\documentclass{amsart} 
\usepackage{amssymb, amscd}
\numberwithin{equation}{section}
\date{January 3 2001}
\def\BB{{\mathbb B}} 
  
\def\CCtimes{{\mathbb C}^\times}
\def\CC{{\mathbb C}}  
\def\DD{{\mathbb D}}  
 
\def\FF{{\mathbb F}} 
\def\GG{{\mathbb G}} 
\def\HH{{\mathbb H}}

\def\LL{{\mathbb L}}
 
\def\PP{{\mathbb P}}
\def\QQ{{\mathbb Q}} 
\def\RR{{\mathbb R}} 
\def\ZZ{{\mathbb Z}}

\def\G{\Gamma}
\def\g{\gamma}

\def\even{{\rm even}} 
 
\def\irr{{\rm irr}} 
 
\def\reg{{\rm reg}}

\def\ss{{\rm ss}}
\def\supp{{\rm supp}} 
\def\st{{\rm st}}
\def\sing{{\rm sing}}

\def\bs{\backslash}
\def\dot{\bullet}

\def\Ccal{{\mathcal C}}
\def\Dcal{{\mathcal D}}
\def\Dcaltilde{\tilde{\Dcal}}
\def\Ecal{{\mathcal E}} 

\def\Hcal{{\mathcal H}}
\def\Ical{{\mathcal I}}

\def\Lcal{{\mathcal L}}

\def\Mcal{{\mathcal M}}
\def\Mcaltilde{\tilde{\Mcal}}
\def\Ocal{{\mathcal O}}
\def\Rcal{{\mathcal R}}
\def\Scal{{\mathcal S}}

\def\Pcal{{\mathcal P}}

\def\Xcal{{\mathcal X}}

\def\That{\hat{T}}
\def\Rhat{\hat{R}}
\def\rhohat{\hat{\rho}}

\newcommand\aut{\operatorname{Aut}}

\newcommand\Gop{\operatorname{G}} 
\newcommand\Hom{\operatorname{Hom}}
\newcommand\Homeo{\operatorname{Homeo}}

\newcommand\pertilde{\operatorname{\widetilde{Per}}}
\newcommand\per{\operatorname{Per}}
\newcommand\pic{\operatorname{Pic}}

\newcommand\proj{\operatorname{Proj}}
\newcommand\sym{\operatorname{Sym}}
\newcommand\GL{\operatorname{GL}}
\newcommand\SL{\operatorname{SL}}
\newcommand\SO{\operatorname{SO}}
\newcommand\PSL{\operatorname{PSL}}
\newcommand\U{\operatorname{U}}
\def\Uu{{\U(\Lambda)}}
\def\lo{{\Lambda_o}}
\def\Ulo{{\U(\Lambda_o)}}
\newcommand\SU{\operatorname{SU}}
\newcommand\PU{\operatorname{PU}}
\newcommand\Symp{\operatorname{Sp}}
\newcommand\braid{\operatorname{Br}}
\newcommand\brcl{\operatorname{BCl}}
\newcommand\brclhat{\operatorname{\widehat{BCl}}}

\newtheorem{theorem}{Theorem}[section]
\newtheorem{lemma}[theorem]{Lemma}
\newtheorem{proposition}[theorem]{Proposition}
\newtheorem{corollary}[theorem]{Corollary}

\theoremstyle{definition}

\newtheorem{example}[theorem]{Example}
\newtheorem{none}[theorem]{}

\theoremstyle{remark} 
\newtheorem{remark}[theorem]{Remark}

\begin{document}

\title{The moduli space of rational elliptic surfaces}

\author{Gert Heckman}
\address{Faculteit Wiskunde en Informatica, Catholic University of
Nijmegen, Postbus 9010, NL-6500 GL Nijmegen, The
Netherlands}
\email{heckman@sci.kun.nl}
\author{Eduard Looijenga}
\address{Faculteit Wiskunde en Informatica, University of
Utrecht, Postbus 80.010, NL-3508 TA Utrecht, The
Netherlands} 
\email{looijeng@math.uu.nl}


\begin{abstract}
We show that the moduli space of rational elliptic surfaces admitting a section
is locally a complex hyperbolic variety of dimension $8$. We compare its 
Satake-Baily-Borel compactification with a compactification
obtained by means of geometric invariant theory, considered by Miranda.
\end{abstract}


\maketitle
\begin{quote}
\center \textit{To Tonny Springer for his $75$th birthday}
\end{quote}

\section*{Introduction}
By a {\it rational elliptic surface} we mean a smooth
complete complex surface that can be obtained from a pencil of
cubic curves in $\PP^2$ with smooth members by successive blowing up 
($9$ times) its base points. A more intrinsic characterization is
to say that the surface is rational and admits a relatively minimal 
elliptic fibration possessing a section. Better yet: it is a smooth 
complete complex surface  whose anticanonical system defines a fibration.  
The description as a blown-up $\PP^2$ is not canonical (in general 
the possible choices are in bijective correspondence with a
weight lattice of an affine root system of type $\hat E_8$), 
but the last characterization makes it plain that the fibration is.
The main goal of this paper is to investigate and describe the
moduli space of these surfaces and certain compactifications
thereof. By assigning to a fiber of a rational elliptic surface its 
Euler characteristic we find a divisor on its base curve, called the
\emph{discriminant divisor}. This discriminant divisor 
is effective and of degree 
$12$. In general it is reduced, meaning that we have $12$ singular fibers,
each of which is a rational curve with a node. It is not difficult to 
show that in that case the discriminant divisor is a complete invariant: 
the projective equivalence class of the discriminant (as a $12$-element 
subset of $\PP^1$) determines the surface up to isomorphism. 
Let us denote by $\Mcal$ the moduli space of rational elliptic surfaces with 
reduced discriminant. One compactification of $\Mcal$ was obtained by 
Miranda \cite{mirandaw} by applying geometric invariant theory to the Weierstra\ss\
models of these surfaces. It gives  a 
projective compactification of $\Mcal$, denoted here by $\Mcal^M$, with
an interpretation of every boundary point as 
corresponding to an isomorphism class of rational elliptic surfaces.

Regarding the discriminant of a rational elliptic surface
as its fundamental invariant leads
to an altogether different compactification of $\Mcal$. 
Let $\Dcal_{12}$ denote 
the space of $\SL (2)$-orbits in the configuration space of
$12$-element subsets of $\PP^1$. A projective compactification 
$\Dcal_{12}^*$ of $\Dcal_{12}$ is obtained by means of 
Geometric Invariant Theory: take the closed
$\SL (2)$-orbits in the space of effective degree $12$ divisors 
on $\PP^1$ that are semistable in the sense that all multiplicities
are $\le\frac{1}{2}.12=6$. There is only one such orbit which is not stable:
it is represented by a divisor which is $6$ times a $2$-element 
subset. 
The variety $\Dcal_{12}^*$ appears in the work of 
Deligne and Mostow \cite{delmost1} as the  Satake-Baily-Borel 
compactification of a $9$-dimensional complex ball $\BB^9$ with just one cusp
(which corresponds to the closed strictly semistable orbit). 
It arises from a period mapping:
for a $12$-element subset $D$ of $\PP^1$, take the cyclic cover 
$C\to \PP^1$ of degree $6$ which totally ramifies in $D$
and then assign to $D$ the abelian variety (of dimension $10$) obtained
from the Jacobian of $C$ by dividing out the Jacobian
of intermediate covers (so that the Galois group acts on this quotient with 
primitive sixth roots of unity). The discriminant gives rise to a closed embedding of $\Mcal$ in $\Dcal_{12}$. Rational elliptic surfaces 
have $8$ moduli, whereas $\dim\Dcal_{12}=9$ and so they define a 
$\SL (2)$-invariant hypersurface in the $12$th symmetric power of 
$\PP^1$. This hypersurface can be 
characterized as defining the $12$-element subsets admitting an equation
that is the sum of a cube and a square.
The compactification of $\Mcal$ we alluded to is the normalization
$\Mcal^*$ of $\Mcal$ in $\Dcal_{12}^*$. A central result of this paper 
is a characterization of the morphism $\Mcal^*\to\Dcal_{12}^*$ in the 
spirit of Deligne and Mostow, namely as a morphism of
Satake-Baily-Borel compactifications defined by an 
`arithmetically defined' hyperball in $\BB^9$. The origin of this 
description is explained by the fact that the degree $6$ cover 
$C\to \PP^1$ naturally comes with a morphism from $C$ to the elliptic curve
of $J$-invarant $0$ when its ramification divisor is the discriminant
of a rational elliptic surface. 

We find that the boundary of $\Mcal$ in 
$\Mcal^*$ is of codimension one and has four irreducible components,
each of which is the closure of a totally geodesic subvariety. Only two
of these irreducible components also appear in Miranda's compactification
and have there the interpretation as parametrizing rational 
elliptic surfaces with a special fiber
(of bicyclic type $I_2$ and of cuspidal type $II$ respectively).

Apart from that, the two compactifications are very much different. The 
natural birational map between $\Mcal^M$ and $\Mcal^*$ is not a morphism 
in either direction and many points of $\Mcal^*$ 
fail to have an interpretation as describing an isomorphism class of a 
rational elliptic 
surface. We therefore consider the closure $\Mcal^{M*}$ of
the diagonally embedded $\Mcal$ in the product of these two 
compactifications. A substantial part of this paper
can be understood as a study of $\Mcal^{M*}$ with 
its projections on $\Mcal^M$ and $\Mcal^*$. 
In the end it turns out that this diagram can be obtained in completely
arithmetic terms (involving a hyperbolic Hermitian lattice over 
the Eisenstein ring); Examples \ref{belangrijkvb1} and \ref{belangrijkvb2}
make this most explicit. The situation is quite similar 
to the relation one of us found between the Baily-Borel 
compactification of the moduli  space of K3 surfaces of degree $2$ 
(resp. $4$) and Shah's GIT compactification  of  the sextic plane 
curves (resp.\ quartic surfaces) 
\cite{looijeng} and the one that Sterk \cite{sterk} found  between the 
moduli space of Enriques surfaces and Shah's GIT compactification of curves on 
$\PP^1\times \PP^1$ of bidegree $(4,4)$ invariant under a certain involution. 

Vakil \cite{vakil} recently showed that some interesting moduli spaces 
define finite monodromy covers of $\Mcal$: the moduli spaces of (1) 
nonhyperelliptic genus $3$ curves endowed with a canonical pencil, 
of (2) genus $4$ curves with an effective even theta characteristic, 
and of (3)  hyperelliptic genus $4$ curves endowed with a noncanonical 
pencil all have this
form. He observes that it then follows from our theorem that these moduli
spaces are locally complex hyperbolic.

As is well-known, $\Mcal$ has also the interpretation as the moduli space
of Del Pezzo surfaces of degree $1$. From this point of view, the above
result neatly fits in a series of similar
characterizations of the moduli spaces of Del Pezzo surfaces 
of given degree: this started with the work of Allcock, Carlson and Toledo 
\cite{all} who associated to a cubic surface $X$ in $\PP^3$ the intermediate 
Jacobian of the cyclic degree $3$ cover of $\PP^3$ ramified along that surface. 
They found that in this way the moduli space of cubic surfaces has the 
structure of a ball quotient. In this case one is so fortunate as to 
have a GIT interpretation of the Baily-Borel compactification so that the 
boundary parametrizes  (mildly) degenerate cubic surfaces. Van Geemen 
(unpublished) and Kondo \cite{kondo1} independently found a similar ball 
quotient description for the moduli space of Del Pezzo surfaces of degree 
two (or equivalently, of  quartic plane curves). It seems that here 
the GIT compactication and the Baily-Borel compactification are related in
a way that is quite similar to the case studied in the present paper.
In particular, neither is a blowup of the other.\par

Let us take the occasion to point out that this is also
the picture for Del Pezzo surfaces of degree four
(in higher degree they are rigid, so this is the remaining case of interest).  
The anticanonical embedding of such a surface realizes that surface in 
$\PP^4$ as the fixed point set of a pencil $P$ of quadrics in $\PP^4$. The
singular quadrics in this pencil define a $5$-element subset $D$ of $P$ 
and the isomorphism type of the pair $(P,D)$ is a complete invariant of 
the surface. The work of Deligne-Mostow identifies the set of such of 
isomorphism types with an open subset of a ball quotient, essentially by 
passing to the Jacobian of the cyclic cover of $P$ of degree $5$ with 
total ramification in $D$.
   
We also mention here that Abramovich and Vistoli \cite{abramvist} defined  
(as a special case of a more general theory of theirs)
a complete Deligne-Mumford stack of moduli of rational elliptic surfaces,  
which is modeled on what we call the \emph{Kontsevich compactification}.  
Here the boundary points label no longer ordinary rational elliptic surfaces, 
but rather surfaces with an `orbispace structure'
\\ 

Let us now describe the contents of the separate sections. 
We begin the paper with a general discussion of the Picard group of a 
rational elliptic surface (Section \ref{ellfibrrev}). Although this
section is brief, we do not use all the material expounded here and we 
therefore advise the reader to consult it as needed.
Section \ref{ellfibrmod1} introduces the main character of this paper,  
the moduli space $\Mcal$. We review Miranda's compactification $\Mcal^M$ of $\Mcal$, which parametrizes elliptic surfaces and we define another 
one, $\Mcal^*$, which does not. 
In Section \ref{meaningful} we define yet another compactification 
that dominates these two and is based on Kontsevich's notion of a stable map. 
This compactification is useful by itself, 
but plays in the present paper only an auxiliary role:
we use it to understand the birational map between $\Mcal^M$ and $\Mcal^*$ in geometric terms. 
In the next two 
Sections \ref{cycliccov} and \ref{centralbraid} we make a careful study 
of the homology of cyclic degree $6$ covers of $\PP^1$ totally ramified 
in $12$ distinct points and the action of a corresponding central extension 
of the braid group of $\PP^1$ with $12$ strands. This discussion belongs to 
algebraic topology rather than to algebraic geometry 
and is independent of the preceding.
Section \ref{bb} recalls the basics of the Satake-Baily-Borel compactification
of a ball quotient and the next section discusses the work of Deligne-Mostow 
for the case that is relevant here. Since this result is a bit hidden 
in their general theory, we outline its proof. 
In passing we obtain a  simple description of the monodromy
group (a unitary group of a rank $10$ lattice over the Eisenstein ring) as
a quotient of the corresponding mapping class group (a centrally extended
braid group). Section \ref{ellfibreis} leads up to the main Theorems
\ref{torelli:3} and \ref{torelli:4} in the next section. 
The final Section \ref{bbmod} is for the most part descriptive. 
It provides what we feel is a natural general context for our results. 
It also suggests an extension of the theory of automorphic forms for ball quotients 
whose geometric counterpart is a theory of compactifications of 
ball quotients with a locally symmetric divisor removed. 
The appendix is devoted to unitary lattices over an Eisenstein ring. 
Part of this is a general discussion, but 
we have also put here the more specific results that we use.\\ 

Some of the initial steps of this work by one of us (GH) were
carried out when he was a visitor of the \'Ecole Normale Superieure at Paris
in May 1998, and he is grateful for the hospitality. 
He also wants to thank Richard Borcherds for an inspiring 
lecture and discussion. We thank Rick Miranda for some helpful correspondence.

Most of the results described here were obtained in the summer of 1999. 

\smallskip
We happily dedicate this paper to our colleague Tonny Springer 
on the occasion of his $75$th birthday.   

\tableofcontents

\section{Rational elliptic surfaces: basic properties}\label{ellfibrrev}
In this section we collect some 
facts---known and perhaps less known---concerning 
rational elliptic surfaces and Del Pezzo surfaces 
of degree one. General references are \cite{manin}, \cite{demazure}, 
\cite{looij:ann}, \cite{morpers} and \cite{dolg_ort}.

By a {\it rational elliptic surface} we shall mean 
a smooth complete rational surface $X$ that admits an
elliptic fibration that is relatively minimal (in the sense that no
exceptional curve is contained in a fiber) and has a section. Then
this fibration is unique since its fibers are the anticanonical
curves on $X$; in particular, its base $P$ is
canonically the  projective line of lines of the plane
$H^0(X,\omega_X^{-1})$. (In fact, any smooth complete surface
whose anticanonical system is a pencil and defines a fibration
is of this form.)
The sections of this fibration are precisely the exceptional 
curves of the first kind of $X$. We can always obtain such a 
surface---though in general in more than one way---as follows: take a
pencil of plane cubic curves  having at least one
smooth member. Its base locus will consist of nine points 
(possibly infinitely near) and blowing these up yields a rational
elliptic surface in our sense (the last blowup giving a section).

It follows from this last description that the Picard lattice of $X$
is isomorphic to the rank $10$ lattice $I_{1,9}$ that
has a basis $\ell ,e_1,\dots ,e_9$ on which the inner product takes
the form $\ell .\ell =1$, $\ell .e_i =0$,  $e_i .e_j=-\delta_{i,j}$.
An isomorphism $I_{1,9}\cong \pic (X)$ can be chosen such that $\ell$
is the class of a line in $\PP^2$ and $e_i$ the class of the exceptional
curve of the $i$th blowup. The class of a fiber of $X\to P$ is the
class of $\omega_X^{-1}$ and is therefore mapped to  $f:=3\ell -e_1-\cdots
-e_9$. 

We first investigate $I_{1,9}$ as an abstract lattice with distinguished 
isotropic vector $f$. A \textit{root} of $I_{1,9}$ is a vector 
$\alpha\in I_{1,9}$ with $\alpha.f=0$
and $\alpha.\alpha =-2$.  The orthogonal reflection with respect to $\alpha$,
\[
s_\alpha : c\mapsto c +(\alpha .c)\alpha
\]
preserves the lattice $I_{1,9}$ and fixes $f$. The set of roots (denoted here 
by $\Rcal$) is an infinite root system; a root basis is 
$\alpha_0:=\ell -e_1-e_2-e_3,\alpha_1:=e_1-e_2,\dots , \alpha_8:=e_8-e_9$, 
which shows that it is of type $\hat E_8$.
The associated Weyl group $W(\Rcal)$ of isometries of $I_{1,9}$ 
generated by the reflections with respect to roots is precisely 
the stabilizer of $f$ in the orthogonal group of $I_{1,9}$ (see 
for instance \cite{vinberg}). We realize $\Rcal$ as an affine root system
(and $W(\Rcal)$ as an affine transformation group) as follows.
The set of vectors $c\in I_{1,9}$ with $c.f= 0$ resp.\ $c.f= 1$
project in $I_{1,9}/\ZZ f$ onto a sublattice $Q$ resp.\ 
an affine lattice $A$ over $Q$. Given a root $\alpha$, then
taking the inner product with that root, makes $\alpha$ appear
as an affine-linear form on $A$. If denote by $\check{\alpha}$ the image 
of $-\alpha$ in $Q$, then the action of $s_\alpha$ in $A$ is given
by $c\mapsto c -(\alpha .c)\check{\alpha}$ and thus $\Rcal$ becomes an 
affine root system on $A$ in the sense of \cite{macdonald}. 
The group $W(\Rcal)$ acts 
faithfully on $A$ and the underlying real affine space
$A_\RR$ receives its standard affine reflection action. 
The image $R$ of $\Rcal$ in $Q$ is a finite root system
of type $E_8$ and spans $Q$. The full translation lattice $Q$ is so 
realized as the translation subgroup of
$W (\Rcal)$. More concretely, the transformation in $I_{1,9}$
associated to $u\in Q$ is the {\it Eichler-Siegel} transformation
\[
T_u: c\mapsto c + (c.f)\hat u - (c.\hat u)f - 
\frac{1}{2}(\hat u.\hat u)(c.f)f,
\]
where $\hat u\in\hat Q$ lifts $u\in Q$. The transformation $T_u$
indeed only depends on $u$ and we have thus defined an 
injective homomorphism $T:  Q\to \SO (I_{1,9})$ of groups.

Let us denote by $\Ecal\subset I_{1,9}$ the set of $e$ with
$e.f=1$ and $e.e=-1$. The natural map $\Ecal\to A$ is a bijection:
if $c\in I_{1,9}$ is such that $c.f=1$, then $(c.c)$ is odd
(this follows from the fact that $Q$ is even and that this is 
true for one such $c$, e.g., $c=e_1$) and so
$e:=c-\frac{1}{2}(1+(c.c))f$ is the
unique element of $c+\ZZ f$ with self-product $-1$. 
So the translation subgroup $T(Q)$ of $W(\hat\Rcal )$ acts simply 
transitively on $\Ecal$.

It is clear that this discussion makes sense in $\pic(X)$ without any
reference to an isomorphism of $(I_{1,9},f)$ onto
$(\pic(X),[\omega^{-1}_X])$. 
We adapt our notation to this situation in an obvious way and write
$f_X, \Rcal_X, \Ecal_X, Q_X, A_X,\dots $. 

An element of $\Rcal_X$ resp.\ $\Ecal_X$ that is the class 
of an irreducible curve is called a \textit{nodal} resp.\ 
\textit{exceptional} class and we denote by
$\Rcal_X^\irr\subset \Rcal_X$ resp.\ $\Ecal_X^\irr\subset\Ecal$ the 
corresponding subset. The following is well-known.

\begin{proposition}
Any irreducible component of a reducible fiber has a 
nodal class and this establishes a bijection between the set of 
irreducible components of reducible fibers and $\Rcal_X^\irr$. 
\end{proposition}

The set $\Rcal_X^\irr$ decomposes according to the set of reducible fibers
$(X_p)_{p\in S}$: 
\[
\Rcal_X^\irr =\sqcup_{p\in S} \Rcal_{X_s}^\irr .
\]
It is convenient to introduce the \textit{closed nodal chamber} as the set of 
$c\in A_\RR$ satisfying $\alpha .c\ge 0$ for all $\alpha\in \Rcal_X^\irr$.
This is a product of closed simplices (a factor for every reducible fiber) 
times an affine space. It is a strict fundamental domain for the action
of the Weyl subgroup $W(\Rcal_X^\irr)\subset W(\Rcal_X)$ in $A_\RR$. 
Let us denote by $Q^\irr_X\subset Q_X$ the image 
of the integral span of $\Rcal_X^\irr$ in $\pic (X)$.

\begin{proposition}\label{exceptional} 
A section of $X\to P$ is an exceptional curve of the first kind and this
identifies the set of sections with $\Ecal_X^\irr$. 
Given $e\in\Ecal_X$, let $e_0\in\Ecal_X$ be the unique element of its
$W(\Rcal^\irr_X)$-orbit mapping to the closed nodal chamber. Then $e_0$ is
the class of a section and $e-e_0$ is a nonnegative linear combination
of nodal classes. The composite map 
$\Ecal_X^\irr\subset \Ecal_X\cong A_X\to A_X/Q^\irr_X$ is a bijection.
\end{proposition}

All of this is known, though perhaps stated somewhat differently in the literature
(see for example \cite{morpers}). 
So $\Ecal_X^\irr$ gets smaller when $\Rcal_X^\irr$ gets bigger.
The generic situation is when $\Rcal_X^\irr=\emptyset$: then 
$\Ecal_X^\irr=\Ecal_X$. The other extreme, $\Ecal_X^\irr$ finite,
happens precisely when $Q_X/Q^\irr_X$ is finite.
The following proposition identifies the rational points of the Picard group
of the generic fiber. 

\begin{proposition}\label{pic}
The group of automorphisms $\aut^0(X/P)$ of $X$ that induce a translation in 
every smooth fiber is faithfully represented in $\pic (X)$. It acts
simply transitively on $\Ecal_X^\irr$ and via the identification
of $\Ecal_X^\irr$ with $A_X/Q^\irr_X$, this group is identified with
the abelian group $Q_X/Q^\irr_X$. It is also the  group of 
automorphisms of $\pic (X)$ that lie in $T(Q).W(\Rcal^\irr_X)$ 
and preserve $\Rcal^\irr_{X_s}$ for every reducible fiber $X_s$.
(This group contains the image of $(Q^\irr_X)^\perp\subset Q$ under $T$ as 
a subgroup of finite index.)
\end{proposition}

This proposition should be known, but since we did not find it 
stated this way, we give a proof. For this we need a property of 
affine Coxeter groups that we recall from \cite{bourbaki}, Ch.~VI, \S~2.
Let $(W, (s_i)_{i\in I})$ be an irreducible Coxeter system of affine type
(with the $s_i$'s distinct) and identify $W$ with its canonical 
representation as an affine transformation group. Denote by $D(I)$ the 
Dynkin diagram on $I$. The normalizer $N(W)$ of $W$ in the affine 
transformation group  acts on $D(I)$ and identifies $N(W)/W$ with 
$\aut (D(I))$. If $I_0\subset I$ is the set \textit{special} vertices
of $D(I)$ (an $i\in I$ is special precisely when every element of
$W$ is the composite of a translation and an element of $W_{I-\{ i\}}$),
then $N(W)/W$ acts faithfully on $I_0$ and the subgroup of $N(W)/W$
induced by translations acts simply transitively on $I_0$. In particular,
if a translation in $N(W)$ fixes a special vertex of $D(I)$, then it lies in
$W$. On the other hand, any element of $N(W)/W$ not coming from a 
translation fixes a special vertex.

We see this illustrated by a Kodaira fibration 
over a smooth curve germ $\Xcal\to \DD$ with special fiber $X_o$ 
(the general fiber is a smooth curve of genus one, the special 
fiber is of Kodaira type). 
If $(C_i)_{i\in I}$ are the distinct irreducible components of $X_o$, then
we have $C_i.C_i=-2$ for all $i$ and if $\sum_i n_iC_i$ is the class 
of the general fiber, then the reflections 
$s_i :c\mapsto c+ (c,C_i)[C_i]$ in $H_2(X_0)$ 
generate an irreducible Coxeter system $(W, (s_i)_{i\in I})$
of affine type acting naturally in the affine 
hyperplane in $\Hom (H_2(X_o),\RR)=H^2(X_o;\RR)$ of forms that 
take the value $1$ on $\sum_i n_iC_i$. Its Dynkin diagram
is just the intersection graph of the $C_i$'s. 
We have $n_i>0$ for all $i$ and $i\in I$ is special precisely when $n_i=1$.
Any automorphism of the general fiber which induces a translation on that fiber
extends to the whole fibration. If it preserves a special component, then it 
preserves every component. So it follows from the preceding that   
that its action on $A_\RR $ is the composite 
of an element of $W$ and a translation. 

\begin{proof}[Proof of \ref{pic}] 
That $\aut (X)$ acts faithfully on $\pic (X)$ is well-known
and easy to prove. 
If $e, e'\in \Ecal_X^\irr$ are represented by sections
$E$, $E'$, then there is a fiberwise translation in the part of 
$X$ that is smooth over $P$ which sends $E$ to $E'$. As recalled above, 
this translation extends as an automorphism $h$ of $X$.
Then $h$ fixes the difference of any two sections, so it 
certainly acts as the identity in $Q_X/Q^\irr_X$. If $E$ and $E'$ meet
a reducible fiber $X_s$ in the same component, then this component
is special. So $h$ fixes every irreducible component of $X_s$. 
The rest of the argument
is now straightforward or follows from the above mentioned property 
of Kodaira fibrations. 
\end{proof}
 
\begin{lemma}\label{eichler}
We have $T(Q_X)\subset \aut^0(X/P).W(\Rcal_X^\irr)$.
\end{lemma}
\begin{proof}
Let $u\in Q_X$. So for every 
$\alpha\in\Rcal_X$ we have $T_u(\alpha )=\alpha -(\alpha .u)f$. 
It follows that for every reducible fiber $X_s$, $T_u$ preserves
the root subsystem $\Rcal_{X_s}$ of $\Rcal_X$ generated by 
$\Rcal_{X_s}^\irr$. So $T_u$ normalizes the associated affine 
Weyl group $W(\Rcal_{X_s}^\irr)$. Choose  
a section $E$. Then $T_u$ sends its class $e\in\Ecal^\irr_X$ to an element 
of the form $w(e')$ with $w\in W(\Rcal_X^\irr)$, where
$e'\in \Ecal^\irr_X$ is the class of a section $E'$. There is a unique 
$h\in\aut^0(X/P)$ that sends $E$ to $E'$. 
We show that $g:=h_*^{-1}w^{-1}T_u$ is in $W(\Rcal_X^\irr)$. It is clear that
$g$ is the identity on the orthogonal complement of 
$\Rcal_X^\irr$ and fixes $e$. Also, for every reducible fiber $X_s$, $g$
normalizes $W(\Rcal_{X_s}^\irr)$ and its image in 
$N(W(\Rcal_{X_s}^\irr))/W(\Rcal_{X_s}^\irr)$ is induced by a translation.
Since $g(e)=e$, it follows that this image is trivial: $g$ acts in
the span of $\Rcal_{X_s}^\irr$ as an element of
$W(\Rcal_{X_s}^\irr)$. This is true for all reducible fibers
and hence $g\in W(\Rcal_X^\irr)$. 
\end{proof}

\begin{remark}
Contraction of an exceptional curve of the first kind with class 
$e\in\Ecal_X^\irr$ produces a smooth rational surface surface $X_e$ with
$\omega_{X_e}.\omega_{X_e}=1$. It follows from Proposition \ref{pic} that its
isomorphism type is independent of the choice of $e$. If all fibers of 
$X\to P$ are irreducible (in other words, $\Rcal_X^\irr=\emptyset$), 
then $\omega_{X_e}^{-1}$ is ample, in other
words, $X_e$ is a Del Pezzo surface of degree one. Conversely, if we
are given a Del Pezzo surface of degree one, then its anticanonical
system consists of irreducible curves and has a unique fixed point. 
Blowing up that point yields an elliptic surface with all its fibers 
irreducible. So the coarse moduli space of Del Pezzo surfaces 
of degree one can be identified with the coarse moduli space of 
smooth rational elliptic surfaces with all its fibers irreducible.
Notice that we have a natural identification of $Q_X$ with
the orthogonal complement of $[\omega_{X_e}]$ in  $\pic (X_e)$.
\end{remark}

\section{Moduli of rational elliptic surfaces I}\label{ellfibrmod1}

\subsection{The Weierstra\ss\ model}
Let $f:X\to P$ be a rational elliptic surface.  
The \emph{discriminant divisor} of $f$ is the divisor on $P$ for which the
multiplicity of $p\in P$ is the Euler characteristic of the fiber $X_p$. 
This is an effective divisor whose degree must be the 
Euler characteristic of $X$, which is $12$.
Assigning to each fiber its modular invariant defines a morphism $J:
P\to\PP^1$. Let us assume that all the singular fibers are of type $I_k$. 
Then $D_\infty:=J^*(\infty )$ is the discriminant divisorof $f$.
In order to understand $J$ over the special
points $0$ and $1$, let us recall that the affine $J$-line is
obtained as the analytic orbifold $\PSL (2,\ZZ)\bs \HH$ with $0$
resp.\ $1$ corresponding to the  singular orbits of
$\omega:=e^{2\pi\sqrt{-1}/6}$ resp. $\sqrt{-1}$. 
The order of ramification 
of the quotient map over such a point is the order of its 
$\PSL (2,\ZZ)$-stabilizer, that is $3$ resp.\ $2$. Since the fibers
of $f$ over  $P-D_\infty$ are smooth, the morphism $J$ is at every point of 
$P-D_\infty$ locally  liftable to a morphism to $\HH$. This implies
that $J^*(0)=3D_0$ and $J^*(1)=2D_1$
with $D_0$ resp.\ $D_1$ a divisor of degree $4$ resp.\ $6$. So $D_\infty$ is 
in the pencil generated by $3D_0$ and $2D_1$. This imposes a
nontrivial condition on $D_\infty$.

To see this, we fix a projective line $P$ and denote by
$H$ the space of sections of $\Ocal_P(1)$. For a nonnegative integer
$k$, $H_k:=\sym^k H$ is then the space of sections of $\Ocal (k)$
and the associated ($k$-dimensional) projective space $P_k$ is  
the linear system of effective degree $k$ divisors on $P$. 
The set of triples $(D_0,D_1,D_\infty)\in P_4\times
P_6\times P_{12}$ with $D_0$ and $D_1$ not a common multiple of an
element of  $P_2$ (to ensure that they generate a pencil)
and $D_\infty$ in the pencil generated by $3D_0$ and $2D_1$ is
an irreducible subvariety of dimension $6+4+1=11$. Denote by
$\tilde\Sigma$ its closure in $P_4\times P_6\times P_{12}$ and by
$\Sigma$ the projection of $\tilde\Sigma$ in $P_{12}$. 
It is clear that $\Sigma$ is irreducible
of dimension $\le 11$. In fact:

\begin{proposition}\label{birat}\label{reducedunique}
The projection $\tilde\Sigma\to\Sigma$ is birational
so that $\Sigma$ is a rational ruled hypersurface in $P_{12}$.
A point $(D_0,D_1,D_\infty)\in\tilde\Sigma$ for which 
$3D_0\not= 2D_1$ comes from a rational elliptic surface. 
It is unique up to $P$-automorphism if
$D_0$ and $D_1$ have disjoint support (a condition fulfilled if 
$D_\infty$ is reduced).  
\end{proposition}

\begin{proof}
Let $W\subset\tilde\Sigma ^2$ be the locus of pairs of distinct
points of $\tilde\Sigma$ with the same image in $P_{12}$. For the
first  assertion it is enough to 
show that $W$ is of dimension $\le 10$. A point of
$(D_0,D_1,D_\infty)\in\tilde\Sigma$ for which $3D_0, 2D_1, D_{12}$ are 
mutually distinct can be represented by a triple 
$(f_0,f_1,f_\infty)\in H_4\oplus H_6\oplus H_{12}$ with $f_\infty
=f_0^3+f_1^2$ so that $D_i$ is the divisor defined by $f_i$. 
Notice that the vector $(f_0^3,f_1^2)\in H_{12}^{\oplus 2}$ is unique
up to a scalar factor. An element
of $W$ is representable by a quadruple $(f_0^3,f_1^2,g_0^3,g_1^2)$ in
$H_{12}$ with $f_0^3+f_1^2=g_0^3+g_1^2$. This identity can also be
written as $(f_1-g_1)(f_1+g_1)=g_0^3-f_0^3$. If the righthand side is
nonzero, then it is factored by the lefthand 
side into two forms of degree six. The family of such factorizations 
(with fixed nonzero righthand side) is of dimension one. 
Since $[f_0:g_0]$ lies in a projective space of dimension $9$, it
follows that $\dim W\le 9+1=10$.

To prove the second assertion we consider the vector 
bundle $\Ecal:=\Ocal_P(2)\oplus\Ocal_P(3)\oplus \Ocal_P$ over $P$.
Denote the projections on its summands by $X,Y,Z$ repectively. 
So for $(f_0,f_1,f_\infty)$ as above, the expression 
$-Y^2Z+ X^3+ 3f_0 XZ^2 + 2f_1 Z^3$ defines a homomorphism
$\Ecal\to \Ocal_P(6)$. Its zero set in the associated projectivized bundle 
$\PP (\Ecal )$ is a Weierstra\ss\ curve over $P$ with modular 
function $J=f_0^3/f_\infty$. If $f_0^3$ and $f_1^2$ are linearly independent, 
then minimal resolution of its singularities gives an elliptic surface 
for which the first summand of $\Ecal$ defines a section. 
This surface is rational.

If $f_0$ and $f_1$ have no nontrivial common zero, then $J$ has degree $12$ and 
$J^*(i)=D_i$ for $i=0,1,\infty$.
Kodaira's theory (see for example \cite{bpv}, Thm.\ 11.1 and 
Subsection \ref{sub:kodaira} below) implies that this elliptic surface 
is unique up to $P$-isomorphism.
\end{proof} 

\begin{remark}
Vakil \cite{vakil} recently proved that the degree of $\Sigma$ is
equal to $3762$. 
In the same paper he also gives several remarkable characterizations
of this hypersurface.
\end{remark}

\begin{example}\label{simplediscr}
This is an example to which we will later return.
Take for $D_\infty$ the $12$th roots of unity in $\CC$, viewed as a
reduced divisor on $\PP^1$. If we take $D_0= 4(0)$ and $D_1 =6(\infty )$,
then clearly $(D_0, D_1, D_\infty)\in\tilde\Sigma$. By the preceding
argument there is a rational elliptic surface 
with $D_\infty$ as discriminant divisor.
\end{example}

We thus recover a result of Dolgachev.

\begin{corollary}[Dolgachev, \cite{dolg}]
The coarse moduli space of rational elliptic surfaces (and hence also
the coarse moduli space of Del Pezzo surfaces of degree one) is
rational.
\end{corollary}
\begin{proof} In view of \ref{birat} we must show that the
$\aut(P)$-orbit space of $\tilde\Sigma$ is rational. 
Generically $\tilde\Sigma $ is fibered
in lines over the product of projective spaces $P_4\times P_6$. 
Let $P'_4\subset P_4$ be the locus where $\aut(P)$ acts freely. Then
$P'_4$ is open-dense in $P_4$, and the orbit space $B:=\aut(P)\bs
P'_4$ is a rational curve. So if $\tilde\Sigma'$ denotes the preimage
of $P'_4$ in $\tilde\Sigma$, then  $\aut(P)\bs\tilde\Sigma'\to B$ is a
morphism to a rational curve whose generic fiber has the structure of a
fibration of  lines over a projective space. This implies that
$\aut(P)\bs\tilde\Sigma'$ is rational.
\end{proof}

\subsection{Miranda's compactification}\label{wstrass}
R.~Miranda gave in his thesis \cite{miranda} a geometric invariant theory
compactification of the space of pencils of cubic plane curves. 
Since pencils with a smooth member define rational elliptic 
surfaces, this leads a compactification of the moduli space of 
(generic) rational elliptic surfaces. Later he found
that the geometric invariant theory of 
Weierstra\ss \ fibrations did that job more 
directly \cite{mirandaw} and so it is this approach that we shall follow. 

Let $U\subset H_4\oplus H_6$ be the open subset of $(f_0,f_1)$ such that 
$f_0^3+f_1^2$ is square free. As was noted in the proof of \ref{reducedunique}, 
the locus $Y^2Z=X^3+ 3f_0 XZ^2 + 2f_1 Z^3$ defines in
$\PP (\Ocal_P(2)\oplus\Ocal_P(3)\oplus \Ocal_P)\times U$ a rational 
elliptic surface $\Xcal_U\to P\times U$ over $U$ with section over $P\times U$. 
The group $\GL(H)$ acts on this fibration. 
Two points of $U$ define isomorphic elliptic
surfaces with section if and only if they are in the same $\GL(H)$-orbit. 
Since the automorphism group of a rational elliptic surface 
acts transisively on 
its sections, it follows that $\GL(H)\bs U$ is the coarse moduli space of 
rational elliptic surfaces with reduced discriminant. We denote that
orbit space by $\Mcal$. A natural projective completion $\Mcal^M$ 
of $\Mcal$ is obtained  
by means of geometric invariant theory applied to the $\SL(H)$-action 
on $H_4\oplus H_6$. With Miranda one easily finds
that  $(f_0,f_1)$ is semistable (resp.\ stable) relative to this action if and only if $f_0^3$ and $f_1^2$ have no nontrivial common zero of order $>6$ 
(resp.\ $\ge 6$).
The proj of the algebra of $\SL(H)$-invariants of the algebra of regular functions on $H_4\oplus H_6$,
\[
\CC [H_4\oplus H_6]^{\SL(H)}=\big( \CC [H_4]\otimes \CC [H_6]\big)^{\SL (H)}
\]
is a projective completion  of $\Mcal$. In more geometric terms:
if $\PP (H_4\oplus H_6)$ stands for the weighted projective space gotten by dividing 
$H_4\oplus H_6-\{ (0,0)\}$ out by the action of the central subgroup  
$\GG_m\subset \GL (H)$, then
\[
\Mcal^M= \SL (H)\bs\bs \PP (H_4\oplus H_6)^\ss .
\]
Here the double backslash indicates that we are forming 
a categorical orbit space. In this case, its closed points are in bijective
correspondence with the closed $\SL(H)$-orbits in $\PP (H_4\oplus H_6)^\ss$. 
We shall refer to $\Mcal^M$  as the {\it Miranda compactification} 
of $\Mcal$. The geometric counterpart of the graded algebra of invariants 
$\CC [H_4\oplus H_6]^{\SL(H)}$ is an orbifold line bundle 
$\Lcal_{\Mcal^M}$ over $\Mcal^M$ such that $\CC [H_4\oplus H_6]^{\SL(H)}$
is the graded algebra of sections of its 
tensor powers with twice the degree. For instance, 
$H^0(\Mcal^M, \Lcal_{\Mcal^M}^{\otimes 2})= H_4^*\otimes 1$ and
$H^0(\Mcal^M, \Lcal_{\Mcal^M}^{\otimes 3})=1\otimes H_6^*$.
 
The minimal strictly semistable orbits in $H_4\oplus H_6$ are 
represented by the pairs of the form $(\lambda f^2,\mu f^3)$ with $f$ a product 
of two distinct linear forms and $\lambda,\mu$ constants that are not both 
zero. In that case the modular function is constant equal to 
$[\lambda^3 :\lambda^3+\mu ^2]\in \PP^1$ and is a complete invariant of the orbit.

A stable orbit can be given more of a geometric content by associating to 
a stable pair $(f_0,f_1)\in H_4\oplus H_6$ the divisor triple 
$(D_0,D_1,D_\infty)\in\tilde\Sigma$ of $(f_0,f_1, f_0^3+f_1^2)$: 
this triple determines the pair $(f_0,f_1)$ up to the action of the 
central subgroup  $\GG_m\subset \GL (H)$. We thus have defined 
an invariant open subset
$\tilde\Sigma^\st$ of $\tilde\Sigma$ characterized by the condition that 
$3D_0$ and $2D_1$ have no point in common of multiplicity $\ge 6$. 

\begin{proposition}[Miranda \cite{mirandaw}]\label{stablepencil}  
A stable orbit defines a rational elliptic surface
all of whose fibers are reduced that is, are of Kodaira type $I_k$ ($k$-gon), 
$II$ (cuspidal curve), $III$ (two rational curves with a common tangent), or $IV$ (three confluent smooth rational curves). Conversely, any 
such rational elliptic surface determines a stable orbit.

A semistable orbit  defines a rational elliptic surface
such that the irreducible components of its fibers have multiplicity
$\le 2$, that is, in addition to the fibers above, we also allow those of type
$I^*_k$. Conversely, such a rational elliptic surface determines a 
semistable orbit in $\tilde\Sigma$. 
The minimal strictly semistable orbits correspond to rational elliptic surfaces
with a $I_4^*$-fiber (such a surface is unique) or with two distinct $I_0^*$-fibers (such a surface has constant modular 
function---see below---and this constant is a complete invariant of the surface).
\end{proposition}

An elliptic surface with two $I_0^*$ fibers is always of the 
following form: start out with a smooth elliptic curve $E$ 
and consider the involution
in $E\times \PP^1$  defined by $(p,[z:1])\mapsto (-p,[-z:1])$. This
involution  has $8$ fixed points that give 
ordinary doubly points on the quotient
surface. A single  blowup resolves these and the resulting smooth
surface $X$ is  rational and fibers over the 
rational curve that is the quotient of
$\PP^1$  by the involution $[z:1])\mapsto [-z:1])$. So to a strictly
semistable orbit of this type is associated a $J$-invariant.

\subsection{Discriminant compactification}\label{discriminant}
We think of $P_{12}=\sym^{12}P$ as the projective space of
effective divisors of degree $12$ on $H$. Let us recall that a $\SL
(H)$-orbit in $P_{12}$ is stable (resp.\ semistable) if and only if
it has no point of multiplicity $\ge 6$ (resp.\ $7$). 
The minimal strictly semistable
elements are of the form $6(a) +6(b)$  with $a$ and $b$ distinct,
hence lie in a single $\SL (H)$-orbit. Let us write 
$\Dcal^\st$ for the ordinary orbit space $\SL(H)\bs P_{12}^\st$ and put 
\[
\Dcal^*:=\SL(H)\bs\bs P_{12}^\ss=\proj (\CC [H_{12}]^{SL(H)}).
\]
So $\Dcal^*$ is a projective one point compactification
of $\Dcal^\st$; the added singleton will be denote $d_\infty$.
The hypersurface $\Dcal^*-\Dcal$ in $\Dcal^*$
parametrizes the nonreduced divisors and is classically called the
{\it discriminant}. There is an orbifold line bundle
$\Lcal_{\Dcal^*}$ on $\Dcal^*$ such that the degree $n$ part of 
$\CC [H_{12}]^{SL(H)}$ is
the space of sections of its $n$th tensor power.
The discriminant is given by the equation $\prod_{1\le i<j\le 12}
(z_i-z_j)^2$ and hence the divisor of a section of $\Lcal_{\Dcal^*}^{11.12}$.

Consider the open part $\Mcal'\subset\Mcal^M$ that parametrizes 
rational elliptic surfaces whose discriminant divisor has no point of
multiplicity $\ge 6$. This means that we discard the surfaces
with a nonreduced fiber or a fiber of type $I_6$ or worse.
So $\Mcal^M-\Mcal'$ is of dimension $\le 3$
and hence everywhere of codimension $\ge 5$ in $\Mcal^M$.
There is an obvious \emph{discriminant morphism} $F :\Mcal'\to \Dcal^*$. 
Assigning to $(f_0,f_1)\in H_4\oplus H_6$ the \emph{discriminant form} 
$f_0^3+f_1^2$ defines an isomorphism 
\[
F^*\Lcal_{\Dcal^*}\cong \Lcal_{\Mcal^M}^{\otimes 6}\vert_{\Mcal'}.
\]
Hence we find:

\begin{corollary}\label{ampelmiranda}\label{bundelvgl}
The algebra of sections 
$\oplus_{k\in\ZZ} H^0(\Mcal',F^*\Lcal_{\Dcal^*}^{\otimes k})$
is zero in negative degrees and of finite type. Its proj  
defines the projective compactification $\Mcal'\subset\Mcal^M$.
\end{corollary}

By Lemma \ref{reducedunique}, $\Mcal$ embeds in $\Dcal$ as a closed hypersurface. 
We denote the normalization of $\Mcal$ in $\Dcal^*$ by $\Mcal^*$
and in $\Mcal^M\times \Dcal^*$ (via the diagonal embedding) by $\Mcal^{M*}$.
The projection $\Mcal^{M*}\to\Mcal^M$ will be special over the 
singleton corresponding to the case where $3D_0=2D_1$ (in other words, 
$(D_0,D_1)=(2(a)+2(b),3(a)+3(b))$ with $a,b\in P$ distinct) and over the locus
where the linear span of $3D_0,2D_1$ has a member with a point of
multiplicity $\ge 7$.
A major goal of this paper is to describe the diagram 
\[
\Mcal^M\leftarrow \Mcal^{M*}\to \Mcal^*
\]
in terms of complex hyperbolic geometry. In particular, we will show
that $\Mcal^*$ is naturally the Baily-Borel compactification of a
ball quotient such that $\Mcal^*-\Mcal$ is the
closure of a union of locally symmetric divisors. This requires a 
better geometric understanding of the above diagram and that is the
topic of the next section.

\section{A geometrically meaningful compactification}\label{meaningful}

We found two compactifications of $\Mcal$ obtained from Geometric
Invariant Theory: one ($\Mcal^M$) based on the Weierstra\ss\ description of a
rational elliptic surface, the other 
($\Mcal^*$) based on the fact that a generic elliptic surface is
defined by its discriminant. It is our goal to define a rather
explicit compactification of $\Mcal$ which dominates both. We also
want it to be geometrically meaningful in the sense that the newly
added points define degenerate elliptic surfaces of some sort.
Together these desiderata imply that the modular function 
of these elliptic surfaces must always be of degree $12$. Since there
exist rational elliptic surfaces whose modular function has lower
degree, there is a price to pay: we must allow the base
to have ordinary double points.

\subsection{Kodaira's theorem}\label{sub:kodaira}
We begin with restating a fundumental
result of Kodaira in more  geometric form. If $P$ is a smooth
complete curve, then
a nonconstant morphism $J: P\to\PP^1$ defines over $P-J^{-1}\{
0,1,\infty\}$ a fibration by elliptic curves given
up to involution. Associated to such a  `Kummer fibration' is
a $\mu_6$-covering of $P$ which will play a central role in this
paper. It is defined as follows. We recall that the abelianization
of $\PSL (2,\ZZ )$ is the cyclic group of
order $6$ with 
\begin{math}
\left(
\begin{smallmatrix}
1 & 1\\
0 & 1\\
\end{smallmatrix}
\right)
\end{math} 
mapping to a generator. We denote that group by $C_6$ and its
generator by $\tau$. So the $\PSL
(2,\ZZ)$ principal bundle over $P-J^{-1}\{ 0,1,\infty\}$ defined by
$J$
determines an unramified $C_6$-covering of $P-J^{-1}\{
0,1,\infty\}$. We extend that covering to a possibly ramified one
over $P$, $C\to P$, by normalizing over $P$. In the case of the
universal example---$J$ is then the identity---this corresponds to
the modular covering $E_o\to \PP^1$ defined by the commutator
subgroup of $\PSL (2,\ZZ)$. The curve $E_o$ is of genus one
and has only one cusp (in other words, it is totally ramified over
$\infty$). If we choose that cusp to be the origin, $E_o$ becomes an
elliptic curve and the fact that it comes with a faithful action of
$\mu_6$ implies that $E_o$ has $J$-invariant $0$. In the general
case, $C\to P$ is simply the normalized pull-back of $E_o\to\PP^1$. 
Here is the list of Kodaira fibers
expressed in terms of the behavior of $J$ at $p$: 
\bigskip

\begin{tabular}{|c|c|c|}
\hline $J(p)$ & $\deg_pJ$ & type\\
\hline\hline $\infty$ & $k(\ge 1)$ & $I_k$ or $I^*_k$\\
\hline $0$ & $0\pmod{3}$ & $I_0$ or $I^*_0$\\
           & $1\pmod{3}$ & $II$ or $IV^*$\\
           & $2\pmod{3}$ & $II^*$ or $IV$\\
\hline $1$ & $0\pmod{2}$ & $I_0$ or $I^*_0$\\
           & $1\pmod{2}$ & $III$ or $III^*$\\
\hline $\notin\{0,1,\infty\}$ & & $I_0$ or $I^*_0$\\
\hline
\end{tabular}
\bigskip

The abelianization of $\SL (2,\ZZ )$ is cyclic of order $12$ with  
\begin{math}
\left(
\begin{smallmatrix}
1 & 1\\
0 & 1\\
\end{smallmatrix}
\right)
\end{math}  
mapping to a generator. We denote group and generator by
$C_{12}$ and $\eta$. So the nontrivial element $-1$ of the
kernel of $\SL (2,\ZZ )\to \PSL (2,\ZZ)$ maps to $\eta^6$.
A relatively minimal elliptic fibration $X\to P$ with $J$ as
modular function determines a $C_{12}$-covering $\tilde C\to P$
which factorizes
over $C\to P$. Thus we associated to every  Kodaira fiber 
an integer modulo $12$, which together with the local behaviour
of $J$ at the corresponding base point determines that fiber. 
Kodaira's basic result says that the lift of the $C_6$-covering to
a $C_{12}$-covering determines $X\to P$ up to $P$-isomorphism
and that any such lift so arises. 
This residue class is in fact the reduction modulo
$12$ of the Euler characteristic of the fiber. So the Euler 
characteristics of the fibers define a further lift to the
integers. (This implies that the Euler
characteristic of $X$ is always divisible by $12$.)
For a fiber with finite $J$-value, its Euler characteristic is the
unique representative of $\ZZ /(12)$ in $\{ 0,1,\dots ,11\}$ (though
$1$, $5$, $7$ and $11$ will not occur), whereas for a fiber $X_p$
with $J(p)=\infty$ it is $\deg_pJ$ (type $I_{\deg_pJ}$) or 
$\deg_pJ+6$ (type $I^*_{\deg_pJ}$). 

A cyclic covering over a smooth rational curve is already given by
the orders of the stabilizers. So if $P$ is rational, then an
elliptic fibration associated to $J$ is already specified 
by a lift of the map $P\to \ZZ /(6)$ defined by $J$ (whose
support will be in $J^{-1}\{ 0,1,\infty\}$) to a finitely supported
map with values in $\ZZ/(12)$. The above receipe defines a lift to
the nonnegative integers and the `integral' of the latter is the 
Euler characteristic of the total space. The total space is
rational precisely when the sum of its fiber Euler characteristics
is equal to $12$. This describes a procedure to obtain all rational
elliptic fibrations and it is the one employed by Miranda in
\cite{miranda2} to recover Persson's classification \cite{persson} of
rational elliptic fibrations up to homeomorphism.

\subsection{Kontsevich compactification}\label{sub:kontcomp} 

Let be given a pair $(J:P\to \PP^1 ,D)$, where
\begin{enumerate}
\item[(a)] $P$ is a complete connected normal crossing curve of
arithmetic genus zero, 
\item[(b)] $J: P\to \PP^1$ a morphism of degree $12$,
\item[(c)] $D$ is a $12$-element subset of the regular part of $P$
contained in $J^{-1}(\infty)$.
\end{enumerate}
For later purposes it will be useful to observe that there then 
exists a $\mu_6$-covering $C\to P$ such that 
\begin{enumerate}
\item[(i)] $C$ is connected normal crossing curve,
\item[(ii)] $C\to P$ is  unramified over $P_\reg-D$ and
\item[(iii)] $C\to P$ is totally ramified over $D$ and
the action of $\mu_6$ in the tangent space of such a ramification
point is the tautological one (i.e., given by scalar
multiplication).
\end{enumerate}
and that this covering is unique up to isomorphism. (The arithmetic
genus of $C$ is easily calculated to be $25$.) 
So to give the pair  $(J:P\to \PP^1 ,D)$ is equivalent to giving a
complete normal crossing curve $C$ with
$C_6$-action as above and a morphism $C\to \PP^1$
constant on orbits of degree $6.12$. The cover $C\to P$ need not be
the pull-back of the modular elliptic curve $E_o\to\PP^1$ for there
may be irreducible components of $P$ in a fiber of $J$
(on which $C\to P$ is necessarily nontrivial). But if we contract all such
components then this is true. In other words, $J$ is covered
by a $C_6$-equivariant morphism $\tilde J:C\to E_o$.

We say that $(J:P\to \PP^1 ,D)$ is {\it Kontsevich stable} if the 
group of its automorphisms that 
induce the identity of $\PP^1$ is finite.
In other words, we require that every connected component 
of $P_\reg -D$ on which $J$ is constant has negative Euler 
characteristic. There is an obvious extension of this notion to
families of such pairs which leads to a well-defined moduli problem. 
Following Kontsevich (\cite{konts} 1.3.2) such pairs have a moduli
stack that is complete, smooth. He also shows that the locus parametrizing 
pairs $(J:P\to \PP^1 ,D)$ with $P$ singular 
defines a normal crossing divisor. His argument shows at the same time  
that the singular points of $P$ are {\it fully smoothable}
in the sense that they are independently smoothable, already at first
order. The underlying variety can be regarded as a coarse moduli space of pairs
$(C, C\to\PP^1)$ obtained as above: here $C$ is a complete connected
normal crossing curve of arithmetic genus $25$ endowed with 
$C_6$-action having in $C_\reg$ exactly $12$ fixed points, each
with tangent character $\chi$ such that the morphism $C\to\PP^1$
is constant on orbits and has degree $6.12$, 
and the group of $\PP^1$-automorphisms of $C$ is finite. 
But the corresponding stack is slightly different.

\begin{remark}
If $(J:P\to \PP^1 ,D)$ is a Kontsevich stable pair, then $(P,D)$ need not be
(Deligne-Knudsen-Mumford) stable as a $12$-punctured curve, but successive 
contraction of its unstable components yields such a curve $(\bar P,\bar D)$
and this curve is unique. There results a morphism from the Kontsevich 
moduli space to the Knudsen-Deligne-Mumford space 
$\Scal_{12}\backslash\overline{\Mcal}_{0,12}$
of stable $12$-punctured rational curves.
\end{remark}

We embed $\Mcal$ in this moduli space by assigning
to a generic rational elliptic fibration $X\to P$ the pair consisting
of its modular function $J:P\to \PP^1$ and the fiber
$J^{-1}(\infty)$. 
The normalization of $\Mcal$ in this moduli space will be called the
{\it Kontsevich compactification} and denoted by $\Mcal^K$.

If $(J:P\to \PP^1, D)$
represents a closed point of $\Mcal^K$, then clearly $D$ will
be contained in $J^{-1}(\infty)$. Specifically, a connected component
of  $J^{-1}(\infty)$ contains as many points of $D$ as the degree of
$J$ on a deleted neighborhood of that component in $P$.
Moreover, every connected component of $J^{-1}(0)$ resp.\
$J^{-1}(1)$ has a basis of deleted neighborhoods in $P$ on which
$J$ has degree divisible by $3$ resp.\ $2$.
The interest of this construction is that such a $J$ 
is still the modular function of an elliptic fibration defined over
the union of the connected components of $P_\reg-D$ on which $J$ is 
nonconstant: if $P'$ is an irreducible component of $P$
on which $J$ is nonconstant, then $J$ determines an elliptic
fibration up to canonical involution. If $p$ is a smooth point of 
$P$, then the fiber over $p$ will be smooth or of type
$I_1$, depending on whether $J(p)$ is finite. In case $p$ is
singular, then we have a singular Kodaira fiber not of type $I_1$.
So if $P'$ has exactly one singular point $p$, then
the Euler characteristic of the fiber over $p$ determined by the
fact that the Euler characteristics of the singular fibers sum up
to $12$. This gives also the answer in the general case since
we can smooth all the singular points of $P$ different from $p$
and do the calculation for this new situation. We thus conclude
that the Euler characteristic of the fiber over $p$ in $P'$ 
must be equal to $12$ minus the degree of $J$ on the connected
component of $P-\{ p\}$ containing $P'-\{ p\}$ plus the multiplicity
of $p$ in $(J|P')^*(\infty)$. 
But beware that in general a singular fiber over a
crossing point will depend on the choice of a branch through
it. For instance, if $P$ has two connected components $P_2$, $P_{10}$ of
degree $2$ and $10$ meeting in a point
$p$ with $J(p)=0$, then the fiber over $p$ in $P_2$ is of
type $II^*$ whereas the fiber over  $p$ in $P_{10}$ is of type
$II$. This issue is adressed and resolved by Abramovich and Vistoli in 
\cite{abramvist} by consistently working in a setting of Deligne-Mumford  
stacks. We shall not go into this here
as it is not needed for what follows. We content ourselves with
observing that $\Mcal^K$ comes as a stack
with a universal morphisms $\Pcal^K\to \Mcal^K\times\PP^1$ of degree
$12$ such that the part of $\Pcal^K$
where this morphism is smooth supports an elliptic fibration for
which $J$ is the modular function. Moreover, $\Pcal^K$ comes (as a
stack) with a $C_6$-covering $\Ccal^K\to \Pcal^K$.

\begin{proposition}
The identity map of $\Mcal$ extends to a morphism from the
Kontsevich compactification $\Mcal^K$ to the Miranda compactification
$\Mcal^M$. Precisely, if $J:P\to \PP^1$
represents a closed point of $\Mcal^K$ and
\begin{enumerate}
\item[(i)] if the fibration has a component $P'$ of $P$ on which
$J$ has degree $>6$, then we assign to $J$ the fibration over this
component (since a nonreduced Kodaira fiber takes off at least $6$
from the degree of modular function, this 
fibration will have only reduced Kodaira fibers);
\item[(ii)] if $P$ has a singular point $p$ with finite $J$-value
such that each component of $P-\{ p\}$ has degree $6$ over $\PP^1$,
then we assign to $J$ the elliptic fibration with constant
modular function $J(p)$ and with two fibers of type $I_0^*$ and
\item[(iii)] if $P$ has a singular point $p$ over $\infty$, such that
each component of $P-\{ p\}$ has  degree $6$ over $\PP^1$, then  
we assign to $J$ the point $m_\infty\in \Mcal^M$ 
(the unique point representing a rational elliptic
surface with a $I_4^*$-fiber).
\end{enumerate}
\end{proposition}
\begin{proof} We begin with proving the first part of (iii).
Suppose that $P$ has a singular point $p$ over
$\infty$. Denote the closures of the connected components of
$P-\{ p\}$ by $P_1$ and $P_2$. Then on $(P_i,p)$ we
have a Kodaira fiber of type $I_{k_i}^*$ for some $k_i\ge 1$. The
Euler characteristic of such a fiber is $6+k_i$ and hence the degree
of $J$ on $P_i$ is $6$.

To see that the birational map from $\Mcal^K$ to $\Mcal^M$
is in fact a morphism, we consider the closure 
$\overline{\Mcal}$ of the diagonal embedding of $\Mcal$ in 
$\Mcal^K\times\Mcal^M$.  Since $\Mcal^K$ is normal it suffices to
prove that the projection $\overline{\Mcal}\to \Mcal^K$ is a
bijection. Or equivalently, that any curve germ in $\Mcal^M$ is
the image of one in $\Mcal^K$. Moreover, we want this lift to be
as prescribed by the proposition. This can be checked in a
straightforward manner.
\end{proof}

\begin{remark}
It can be shown that the natural morphism 
$\Mcal^K\to \Scal_{12}\backslash\overline{\Mcal}_{0,12}$ 
is finite. This implies
that is also possible to define $\Mcal^K$ as the normalization of
$\Mcal$ in $\Scal_{12}\backslash\overline{\Mcal}_{0,12}$. 
Though this avoids appeal to the Kontsevich 
moduli space, we shall need the more powerful interpretation that comes 
with the latter.
\end{remark} 

The identity also extends as a morphism $\Mcal^K\to\Mcal^*$ as
follows. Let be given an allowable pair $(J:P\to \PP^1, D)$
representing  a closed point of $\Mcal^K$. If there
exists an irreducible component $P_c$ of $P$ such that the direct
image of $D$ under the natural retraction $P\to P_c$ is a
stable divisor (all multiplicities $< 6$), then this
irreducible component is unique---we shall call it the 
{\it central component} of $(P,D)$---and we assign to $(J:P\to \PP^1,
D)$ the corresponding point of $\Dcal^*$. 
If no such component exists, then 
there is a unique singular point $p_c$---the {\it central point} of 
$(P,D)$---such that $D$ has $6$ points in each connected 
component of $P-\{ p_c\}$, and we then assign to $(J:P\to \PP^1, D)$
the point of $\Dcal^*$ that corresponds to the unique minimal
semistable orbit (the orbit of divisors that 
have two distinct points, each with multiplicity $6$).
In either case we allow ourselves a mild abuse of 
language by referring to this point of $\Dcal^*$ as the {\it
discriminant} of $(J:P\to \PP^1, D)$.
It is not difficult to verify that this 
defines a morphism $\Mcal^K\to\Dcal^*$. 
Since $\Mcal^K$ is normal this morphism will factorize over
$\Mcal^*$. 

So $\Mcal^K$ dominates $\Mcal^{M*}$.
Understanding of $\Mcal^K$ will help us in understanding
$\Mcal^{M*}$.

\subsection{A partial list of strata} 
Let us describe the generic points of
$\Mcal^K-\Mcal$ (these turn out to be all hypersurfaces).
If $X\to P$ is a generic rational elliptic fibration 
(so with smooth base $P$ and reduced discriminant), then
the modular function $J: P\to \PP^1$ is a degree $12$ covering with
the property that the local degree
of $J$ at a point over $0$ resp.\ $1$ is always equal to $3$ resp.\
$2$. Following Riemann-Hurwitz, the discriminant of $J$ must then
have the form  $8(0)+6(1)+R$, with $R$ of degree $8$. 
This divisor gives us the $8$ moduli parameters. 
Degeneracies will occur when $\supp (R)$ meets $0$, $1$ or $\infty$.
The computation of (co)dimension is based on the full
smoothability property.  

In the list below we make use of a small part of Persson's
classification \cite{persson}. For instance, we use the fact that the
rational elliptic fibrations with a fiber of Kodaira type 
$I_k$ ($8\not= k\le 9$), $II$, $III$, $IV$,
$I_k^*$ ($k\le 4$) respectively are parametrized by an
irreducible variety. We excluded the $I_8$-case since there are two
types of fibrations with an $I_8$ fiber: in one case $(I_8')$
the classes of the irreducible components in the Picard group generate
a primitive sublattice and in the other case $(I_8'')$ the sublattice 
is of index two in a primitive sublattice and either case is
parametrized by an irreducible variety. The two cases can be
distinguished by the fact that in case $(I_8')$
the fiber can degenerate in a $I_9$ fiber, whereas this is
not possible for the $(I_8'')$ case. But either can degenerate 
into a $I_4^*$-fiber and is a degeneration of a $I_7$-fiber.

\begin{enumerate}
\item [($I_{k\ge 2}$)] Then
$P$ has an irreducible component $P_{12}$ of degree $12$
over $\PP^1$ and there is a $z\in P_{12}$ 
where $J|P_{12}$ has local degree $k$. We have an extra component  
$P_0$ in $J^{-1}(\infty)$ which meets
$P_{12}$ in the ramification point. This component contains $k$
points of $D$ and so it is central if and only if $k\ge 7$. Hence
the discriminant has a point of multiplicity
$\min\{ k,12-k\}$. The image of this hypersurface of $\Mcal^K$ 
in $\Mcal^M$ is of dimension $9-k$, whereas its image
in $\Mcal^*$ is of dimension $9-k$ for $k=2,3,4,5$, of dimension $0$
for $k=6$, and of dimension $k-2$ for $k=7,8,9$. The hypersurface in 
question is irreducible
unless $k=8$, in which case there are two irreducible components.
\item[($II$)] Then $P$ has two
irreducible components $P_{10}$, $P_2$ of degrees resp.\ $10$
and $2$ over $\PP^1$ meeting in a point $p$ with $J$-value $0$.
The component $P_2$ ramifies simply over
$0$ and $1$; the component $P_{10}$ has fiber over $0$ resp.\ $1$ of
type $(1,3^3)$ resp.\ $(2^5)$. Over $(P_{10},p)$ we have a
fiber of type $II$ (a cuspidal fiber) and over 
$(P_{2},p)$ a fiber of type $II^*$ (an $\hat E_8$-fiber). 
The central component is $P_{10}$
and the discriminant has a point of multiplicity $2$. 
The images of this subvariety in $\Mcal^*$ and $\Mcal^M$
are hypersurfaces. 
\item[($III$)]  This case and the next are similar to the preceding
case. Here $P$ has two irreducible components $P_9, P_3$ of degrees
$9$ resp.\ $3$ over $\PP^1$ meeting in a point $p$ with with
$J$-value $1$. The component $P_3$ ramifies totally
over $0$ and has a point of simple ramification over $1$; the
component $P_9$ has fiber over $0$ resp.\ $1$ of type $(3^3)$ resp.\
$(2^4,1)$. Over $(P_9,p)$ we have a Kodaira fiber of type
$III$,  and over $(P_3,p)$ one of type
$III^*$ (an $\hat E_7$-fiber). The central component is
$P_9$ and the discriminant has a point of multiplicity $3$. 
The images of this subvariety in $\Mcal^*$ and $\Mcal^M$ are of
codimension two (since we forget $P_3$).
\item[($IV$)] Now $P$ has two irreducible components $P_8, P_4$ 
of degrees $8$ resp.\ $4$ over $\PP^1$ meeting in a point $p$ with
$J$-value $0$.
The component $P_4$ has fiber over $0$ resp.\ $1$ of type 
$(3,1)$ resp.\ $(2^2)$, whereas for $P_8$ 
these data are $(3^2,2)$ resp.\ $(2^4)$. 
They meet in their points of smallest
ramification. Over $(P_8,p)$ we have a 
fiber of type $IV$ and over $(P_4,p)$ a fiber of type
$IV^*$ (an $\hat E_6$-fiber). The central component is
$P_8$ and the discriminant has a point of multiplicity $4$. 
The images of this subvariety in $\Mcal^*$ and $\Mcal^M$ are of
codimension three. 
\end{enumerate}
The following cases involve Kodaira fibers of type $I_k^*$. In all
these cases, $P$ has two irreducible components $P_6$, $P'_6$ 
that are both of degree $6$ over $\PP^1$.

\begin{enumerate}
\item[($I_0^*$)] $P_6$ and $P'_6$ meet in a point $p$ with
$J(p)$ finite. Over $(P_6,p)$ and $(P'_6,p)$ 
we have fibers of type $I_0^*$. The point $p$ is central
and so the discriminant is the orbit of the divisor with two points
of multiplicity $6$. The image of this subvariety in $\Mcal^M$
is a curve and its image in $\Mcal^*$ is a singleton.
\item[($I_{k,k'}^*$)] Here $P_6$ and $P'_6$ are
separated by a central component $P_c$ contained in
$J^{-1}(\infty)$. If $P_c$ meets $P_6$ in $p$, then we have a 
Kodaira fiber of type $I^*_k$ at  $(P,p)$,
where $k=\deg_p(J|P)\in\{ 1,2,3,4\}$. Similarly we find 
a Kodaira fiber of type $I^*_{k'}$ for
$P'_6$. So  $P_c$ meets $D$ in $k+k'$ points.
Hence the discriminant has a point of multiplicity $6-k$ and one 
of multiplicity $6-k'$. This defines a hypersurface in $\Mcal^K$,
whose image in $\Mcal^*$ has dimension $k+k'-1$ (so we get a
hypersurface in $\Mcal^*$ precisely when $k=k'=4$). Its image
in $\Mcal^M$ is a singleton.
\end{enumerate} 

So the boundary of $\Mcal$ in $\Mcal^K$ is a 
union of irreducible hypersurfaces $\Mcal^K(F)$, where $F$ runs over
the Kodaira symbols $I_k$, $k=2,\dots, 7,9$, $I_8'$, $I_8''$, $II$,
$III$, $IV$,
$I^*_0$, $I^*_{k,k'}$ with $k,k'=1,2,3,4$. 
Let us write $\Mcal^{M*}(F)$ for the image of $\Mcal^K(F)$ in 
$\Mcal^{M*}$ and let $\Mcal^M(F)$ and $\Mcal^{M*}(F)$
have a similar meaning. The dimensions of these
subvarieties are listed in the table below. 

It is not hard to check that $\Mcal^{M*}(I_2)$ contains
$\Mcal^{M*}(I_k)$ when $k\le 5$, that
$\Mcal^{M*}(II)\supset \Mcal^{M*}(III)\supset \Mcal^{M*}(IV)$
and that $\Mcal^{M*}(I_{4,4}^*)$ contains $\Mcal^{M*}(I_{k,k'}^*)$.
From these and similar incidence relations we deduce:

\begin{enumerate}
\item[(i)] The irreducible components of the boundary of $\Mcal$
in $\Mcal^{M*}$ are the hypersurfaces
$\Mcal^{M*}(I_2)$, $\Mcal^{M*}(II)$, $\Mcal^{M*}(I_7)$,
$\Mcal^{M*}(I_8')$,
$\Mcal^{M*}(I_8'')$,
$\Mcal^{M*}(I_9)$, $\Mcal^{M*}(I^*_{4,4})$, the curve
$\Mcal^{M*}(I^*_0)$ and the threefold $\Mcal^{M*}(I_6)$.
\item[(ii)] For $k=7,9$ we have
$\Mcal^{M*}(I_k)=\Mcal^*(I_k)\times\Mcal^M(I_k)$ and
$\Mcal^{M*}(I_8)^{(i)}=\Mcal^*(I_8^{(i)})\times \Mcal^M(I_8^{(i)})$ for $i=1,2$.
\item[(iii)] We have inclusions 
\begin{align*}
\Mcal^*(I_9) &\supset \Mcal^*(I_8')\supset \Mcal^*(I_7),\\
\Mcal^M(I_9) &\subset \Mcal^M(I_8')\subset \Mcal^M(I_7),\\
\Mcal^*(I^*_{4,4}) &\supset \Mcal^*(I_8'')\supset \Mcal^*(I_7),\\
\Mcal^M(I^*_{4,4}) &\subset \Mcal^M(I_8'')\supset \Mcal^M(I_7).
\end{align*}
\item[(iv)] The projection of $\Mcal^{M*}(F)\to \Mcal^{*}(F)$ is
birational for $F=I_2$, $I_9$, $II$, $I^*_{4,4}$ and a collapse onto
a point for $F=I_6$, $I^*_0$. 
\item[(v)] The projection of $\Mcal^{M*}(F)\to \Mcal^M(F)$ is
birational for 
$F=I_2$, $II$, $I_6$, $I^*_0$ and and a collapse onto a point
for $F=I_9$, $I^*_{4,4}$. 
\end{enumerate}
The following statements then follow in a straightforward manner:

\begin{corollary}\label{strata}
The boundary of $\Mcal$ in $\Mcal^*$ is the union of
the irreducible hypersurfaces $\Mcal^*(I_2)$, $\Mcal^{*}(II)$,
$\Mcal^{*}(I_9)$ and $\Mcal^{*}(I^*_{4,4})$.
Moreover,
\begin{enumerate}
\item[(i)] $\Mcal^{*}(I_9)\cap\Mcal^{*}(I^*_{4,4})=\Mcal^{*}(I_8')$,
\item[(ii)] $\Mcal^{*}(I_8'')\subset \Mcal^{*}(I^*_{4,4})$, 
\item[(iii)] $\Mcal^{*}(I_8')\cap\Mcal^{*}(I_8'')=\Mcal^{*}(I_7)$,
\item[(iv)]  $\Mcal^{*}(I_6)$ is a singleton contained in
$\Mcal^{*}(I_7)$
and $\Mcal^{*}(I_0^*)$ is a singleton contained in
$\Mcal^{*}(I^*_{4,4})$.
These two make up the preimage of $d_\infty\in \Dcal^*$ in
$\Mcal^*$.  
\end{enumerate}
\end{corollary}

\bigskip
\begin{tabular}{|c|c|c|c|}
\hline $F$   & $\dim \Mcal^*(F)$ & $\dim \Mcal^M(F)$ & $\dim
\Mcal^{M*}(F)$\\
\hline\hline $I_2$ & 7 & 7 & 7\\
\hline $I_{k\le 5}$ & $9-k$ & $9-k$ & $9-k$\\
\hline $I_6$ & 0 & 3 & 3\\
\hline $I_7$ & 5 & 2 & 7\\
\hline $I_8', I_8''$ & 6 & 1 & 7\\
\hline $I_9$ & 7 & 0 & 7\\
\hline $II$  & 7 & 7 & 7\\
\hline $III$ & 6 & 6 & 6\\
\hline $IV$  & 5 & 5 & 5\\
\hline $I^*_0$ & 0 & 1 & 1\\
\hline $I^*_{k,k'}$ & $k+k'-1$ & 0 & $k+k'-1$\\
\hline
\end{tabular}

\section{Homology of a cyclic covering}\label{cycliccov}

\subsection{Symplectic lattices with symmetries}\label{symplsym}
Let be given a finite abelian group $G$ that acts
(morphically) on
a symplectic lattice $L$. We then extend the symplectic form as a
sesquilinear form over $\ZZ [G]$ by 
\[
\langle\; ,\; \rangle :L\times L\to \ZZ[G],\quad 
(a,b)\mapsto \sum_{g\in G} (a\cdot gb) g=\sum_{g\in G} (g^{-1}a\cdot
b)g. 
\]
Indeed, this form is $\ZZ[G]$-linear in the first argument and
$\langle b,a \rangle =-\overline{\langle a,b \rangle}$
(where the overline is the involution which sends each element of $G$
to its inverse). So if we multiply the form by any anti-invariant
element of
$\ZZ[G]$ (such as $g-g^{-1}$ for some $g\in G$), then  we get a
Hermitian form over $\ZZ[G]$. 

We take $G$ to be a cyclic group of order $6$, $C_6$,
with a given generator $\tau\in C_6$. Let
$\chi :\ZZ [C_6]\to \CC$ be the character that sends
$\tau$ to $\omega:=e^{2\pi\sqrt{-1}/6}$. 
The image of this character is the ring of integers $\ZZ
+\ZZ\omega$. We call this ring the {\it Eisenstein ring} and 
denote it by $\Ocal$. For the lattice $L$ as above,
$L_\Ocal :=\Ocal\otimes_{\ZZ[C_6]}L/(\text{torsion})$ is the biggest
torsion free quotient of $L$ on which $C_6$ acts through $\Ocal$.
This quotient of $L$ is realized as the image of $L$ under the 
natural `eigenprojection' 
$\CC\times_\ZZ L\to (\CC\times_\ZZ L)_\chi$.
The composition of the sesquilinear form above with $\chi$
factorizes over a skew-hermitian ($\Ocal$-valued) form:
\[
\phi : L_\Ocal\times L_\Ocal\to\Ocal .
\]
We make this a Hermitian form by multiplying with a square root
of $-3$: we put 
\[
\theta :=\omega -\omega^{-1},
\]
and let our Hermitian form be 
\[
\psi (a,b):=-\theta \phi (a,b).
\] 
As we will show in the Appendix, 
such Hermitian lattices can also be gotten from quadratic forms  
with $C_6$-symmetry.

\begin{example}\label{eisensteincurve}
Let $E_o$ be the elliptic curve of $J$-invariant $0$. It admits a
faithful action of $C_6$ with $\tau$ acting 
on the tangent space at the origin as
multiplication by $\omega$. Note that $H^1(E_o)$
is a free $\Ocal$-module of rank one. The generators
make up a $C_6$-orbit and if $c$ is any one of them, then
$\phi (c,c)= \sum_{i=0}^5 (c,(\tau^*)^{-i}c)\omega^i=2\theta$
and so $\psi (c,c)= -\theta .2\theta=6$.
For later reference we also note that $(\tau^*)^{-1}$ acts on
$H^{1,0}(E_o)$ as multiplication by $\omega^{-1}$.
\end{example}
\begin{example}
Here is another example. Take $L:=\ZZ[C_6]/(\sum_{i=0}^5 \tau^i)$
(which, as a $\ZZ[C_6]$-module, is isomorphic to the augmentation
ideal of $\ZZ[C_6]$). We equip it with the symplectic form 
\[
\tau^i\cdot \tau^j=
\begin{cases}
\pm 1 & \text{if $j=i\pm 1$},\\ 
0 & \text{otherwise.}
\end{cases}
\]
We have $\langle 1,1\rangle = \tau -\tau^{-1}$ and so  
for the image $e$ of $1$ in $L_\Ocal$ we have  
$\psi (e,e)= -\theta .\theta =3$.
\end{example}

\subsection{Cyclic covers}
Let $\pi :C_o\to\PP^1$ be the smooth $C_6$-covering of the projective
line that has total ramification over the 12th roots of unity in the
unit circle and with the generator $\tau$ of $C_6$ acting as
multiplication by 
$e^{2\pi\sqrt{-1}/12}$ on the tangent space of the ramification
points. An affine equation for this curve is $w^6+z^{12}=1$ with
$\tau$
acting as $\tau (z,w)=(z,\omega w)$ and $\pi (z,w)=z$. There is
also $C_{12}$-symmetry, with a generator $\eta$ of $C_{12}$ acting as
$\eta (z,w)=(e^{2\pi\sqrt{-1}/12}z, w)$. So we have an action
$C_6\times C_{12}$ on $C_o$. Our first goal is to describe $H_1(C_o)$
as a module over 
\[
R:=\ZZ[C_6\times C_{12}]= \ZZ[\tau ,\eta ]/(\tau^6-1 ,\eta^{12}-1). 
\]
We make use of F.~Pham's description \cite{pham} of the 
homology (with its intersection form) of the affine piece 
$C_o':= C_o -\pi^{-1}(\infty)$.
Consider the real part of $C_o'$ defined by $x^{12} +u^6=1$ with
$x$ and $u$ in the unit interval. We orient it as going from
$(0,1)$ to $(1,0)$ and denote the singular $1$-simplex thus defined
by  ${\bf e}$. Since ${\bf e}$ is not fixed by any element of
$C_6\times
C_{12}$, ${\bf e}$ generates a free $R$-submodule of the module of
singular
$1$-chains on $C_o'$. Pham observes that
\[
e:= (1-\tau)(1-\eta){\bf e}
\]
is a $1$-cycle with the property that it generates $H_1(C_o)$ as an
$R$-module. Since $R\bf{e}$ does not contain nonzero boundaries,
$H_1(C_o')$ gets
identified (as an $R$-module) with the ideal $(1-\tau)(1-\eta)R$. 
The annihilator of $(1-\tau)(1-\eta)$ in $R$
is the ideal $(\sum_{i=0}^5\tau ^i, \sum_{i=0}^{11}\eta^i)R$ and so
the dual module $H^1(C_o')$ appears naturally as a quotient:
\[
\ZZ[\tau ,\eta ]/(\sum_{i=0}^5\tau ^i, \sum_{i=0}^{11}\eta^i)\cong
H^1(C_o'),
\quad 1\mapsto e^*.
\]
Pham also describes the intersection pairing: the adjoint
homomorphism 
$H_1(C_o)\to H^1(C_o)$ is the antihomomorphism of $R$-modules given
by
\[
e=(1-\tau)(1-\eta){\bf e}\mapsto -(1-\tau)(1-\eta)(1-\tau\eta)e^*.
\]
Notice that the kernel of this map is 
$(1-\tau)(1-\eta)(\sum_{i=0}^{11}(\tau\eta)^i)R{\bf e}$.
The inclusion $C_o'\subset C_o$ induces a surjection on $H_1$; in
fact, $H_1(C_o)$ can be identified with the image of $H_1(C_o)\to
H^1(C_o)$
(compatibly with the intersection pairing). So we find an isomorphism
\[
(1-\tau)(1-\eta)R/(1-\tau)(1-\eta)(\sum_{i=0}^{11}(\tau\eta)^i)\cong
H_1(C_o),\quad (1-\tau)(1-\eta)\mapsto e.
\]
We will identify the lefthand side with the quotient ring
\[
A:=\ZZ[\tau ,\eta ]/(\sum_{i=0}^5\tau ^i, \sum_{i=0}^{11}\eta^i,
\sum_{i=0}^{11}(\tau\eta)^i),
\]
so that $1$ corresponds to $e$. (So as a 
$\ZZ [C_6]$-module, $A$ is generated by $\{\eta^i\}_{i=1}^{10}$.)
The sesquilinear extension of the intersection pairing is given by 
\[
\langle a,b\rangle_R =
(1-\tau)(1-\eta)a.\overline{b(\tau\eta -1)}\in R,\quad 
a,b\in A.
\] 
If we merely regard $H_1(C_o)$ as a $\ZZ [C_6]$-module, then the 
intersection form defines a sesquilinear pairing 
\[
\langle\; ,\;\rangle_{\ZZ [C_6]} :H_1(C_o)\times H_1(C_o)\to \ZZ
[C_6]
\]
that is $\ZZ [C_{12}]$-invariant. The two are of course related by
\[
\langle a,b\rangle_R=\sum_{i=0}^{11} 
\langle ae,\eta^ibe\rangle_{\ZZ [C_6]}\eta^i.
\]
Reducing modulo the ideal generated by $\tau^2-\tau +1$ 
yields sesquilinear pairings 
\[
\langle\; ,\;\rangle _\Ocal : H_1(C_o)_\Ocal\times H_1(C_o)_\Ocal\to
\Ocal\quad\text{and}\quad
\langle\; ,\;\rangle _{\Ocal[C_{12}]} : A_\Ocal\times A_\Ocal \to 
{\Ocal[C_{12}]}
\]
that are related in the same way. The associated Hermitian 
forms are defined by multiplying these by
$-\theta=-\omega(1+\omega)$:
\[
\psi (ae,be):=-\omega(1+\omega)\langle ae,be\rangle_\Ocal\quad\text{
and}\quad
\Psi (a,b):=-\omega(1+\omega)\langle a,b\rangle_{\Ocal[C_{12}]},
\] 
so that
\begin{align*}
\sum_{i=0}^{11}\psi (e,\eta^ie)\eta^i&=\Psi (e,e)\\
&= -\omega(1+\omega)(1-\omega)(1-\eta)(\omega^{-1}\eta^{-1}-1)\\
&= -(1+\omega)((-1-\omega^{-1}) +\eta +\omega^{-1}\eta^{-1})\\
&= 3 -(1+\omega )\eta -(1+\omega^{-1})\eta^{-1}.
\end{align*}
In other words,
\begin{equation}\label{psi}
\psi (e,\eta^ie)=
\begin{cases}
3 &\text{ if } i=0,\\
-1-\omega &\text{ if } i=1,\\
-1-\omega^{-1} &\text{ if } i=-1,\\
0 &\text{ otherwise.}
\end{cases}
\end{equation}
Since $\psi$ is $C_{12}$-invariant,
these formulae completely describe $\psi$ on the generators
$\eta^ie$.

Let us denote by $\Lambda$ the Hermitian $\Ocal$-module underlying
$A_\Ocal$.
So if  $r_i\in\Lambda$ denotes the image of $(\omega\eta)^i$, 
then $(r_1,\dots ,r_{10})$ is a $\Ocal$-basis of $A_\Ocal$ on which
$\psi$ is
given by
\begin{equation}
\psi (r_i,r_j)=
\begin{cases}
3 &\text{ if } j=i,\\
\theta &\text{ if } j=i+1,\\
0 &\text{ if } j>i+1.
\end{cases}
\end{equation}
Notice that for $k\le 10$, the annihilator of the span of $r_1,\dots
,r_{k-1}$ contains the span of $r_{k+1},\dots ,r_{11}$. It is not
hard 
to see that it is in fact equal to it.

\begin{remark}
The homology class of $e$ can be represented more simply as
follows. The closed sector of the (closed) unit disk in the $z$-line 
with $\arg (z)$ between $0$ and $2\pi/12$ has a unique lift to $C$ 
passing through $(0,1)$. If we give this lift its complex
orientation, then it becomes a singular $2$-simplex whose boundary of
is of the form ${\bf e} +{\bf\epsilon}-\eta {\bf e}$, where
${\bf\epsilon}$ is a lift of the arc on the unit circle. So
$(1-\eta){\bf e}$ is
homologous to $-{\bf\epsilon}$. Hence $e=(1-\tau)(1-\eta){\bf e}$ is
homologous
to $(\tau -1){\bf\epsilon}$.
\end{remark}

\begin{remark}\label{isotropic}
It is easy to check that the $\Ocal$-sublattice of $A_\Ocal$ spanned
by $\eta^i$, $i=0,\dots ,k$ is of rank $\min\{ k+1,10\}$ and 
positive definite for $k\le 3$,
positive indefinite for $k=4$, and hyperbolic for $k\ge 5$. Since
multiplication by $\eta$ is a lattice automorphism it follows that
the  $\Ocal$-sublattice spanned by all $\eta^i$ with 
$i\not\equiv 5\pmod{6}$ is a positive (indefinite) sublattice
of rank at least $9$. This is clearly also the maximal rank of a
positive sublattice, so it is of the form $l_0^\perp$ for some
$0$-vector $l_0$. A small calculation shows that we can take 
$l_0= (1+(1+\omega )\eta +2\omega\eta^2 +(2\omega -1)\eta^3 +
(\omega -1)\eta^4)e$.
\end{remark}

\section{A central extension of a braid class group}\label{centralbraid}
\subsection{Braid and braid class groups}\label{braid}
This section reviews some facts concerning the braid groups of
$\CC^\times$ and $\PP^1$. We adhere to the categorical convention 
for the composition law in fundamental groupoids: $\alpha\beta$ 
means that the path $\alpha$ comes after the path $\beta$.

We first establish the terminology.
Fix a positive integer $d$. For any topological surface $X$ we denote
by $X(d)$ the configuration space of $d$-element subsets of $X$. 
The {\it braid group of $X$ with $d$ strands} $\braid_d(X)$ of
$X$ is by definition the fundamental group of $X(d)$. The 
latter requires a choice of base point and so strictly speaking this
group is only defined up to conjugacy. The group $\Homeo(X)$ of 
self-homeomorphisms of $X$ acts on $X(d)$. The image of the 
the fundamental group of the identity component,
$\pi_1(\Homeo(X)^0,1)$
in the  fundamental group of $X(d)$ is normal and we shall refer to
the quotient group as the $d$-pointed {\it braid class group} of $X$,
$\brcl_d(X)$. For $X=\PP^1$ we will often omit $X$ and simply write 
$\braid_d$ and $\brcl_d$.
 
An alternative characterization of $\brcl_d(X)$ is as a mapping class 
group: if we fix a $d$-element subset $S$ of $X$, then $\brcl_d(X)$ 
is the group of isotopy classes of self-homeomorphisms of the pair $(X,S)$ 
that are trivial as an absolute isotopy class of self-homeomorphisms
of $X$. This also gives $\brcl_d$ the interpretation as the orbifold
fundamental group of the moduli space $\Scal_d\bs\Mcal_{0,d}$ of 
smooth rational curves with $d$ punctures.
 
We first consider the case $X=\CC^\times$.
We take as a base point $*$ for $\CC^\times (d)$ the set $\mu_d$ of 
$d$th roots of $1$. We have two special elements $R$ and $T$ of
$\braid_d(\CC^\times )$: $R$ is defined in $\braid_d(\CC^\times )$ by
$t\in [0,1]\mapsto e^{2\pi\sqrt{-1}t/d}.\mu_d $,  
and $T$ is represented by the loop that leaves all elements of
$\mu_d$ in place except $1$ and $e^{2\pi\sqrt{-1}/d}$: these traverse
(in counterclockwise direction) half of the circle that has the segment 
$[1,e^{2\pi\sqrt{-1}/d}]$ as a diameter. 
These two elements generate $\braid_d(\CC^\times )$, but in order to
get a useful presentation of $\braid_d(\CC^\times )$ it is better to
enlarge the number of generators. Let $T_k:=R^kTR^{-k}$ ($k\in\ZZ /d$). 
Clearly, $T_k$ relates to the pair $(e^{2\pi\sqrt{-1}k/d},e^{2\pi\sqrt{-1}(k+1)/d})$
in the same way as $T$ to $(1,e^{2\pi\sqrt{-1}/d})$. 
These elements satisfy:

\begin{align}\label{braidrel1}
\begin{split}
T_kT_{k+1}T_k&=T_{k+1}T_kT_{k+1},\quad  k\in\ZZ/d,\\
T_kT_l&=T_lT_k,\quad  k, l \in\ZZ/d \text{ and }k-l\not=\pm 1.
\end{split}
\end{align}
Together with the obvious relations 
\begin{align}\label{braidrel2}
RT_kR^{-1}=T_{k+1},\quad  k\in\ZZ/d,
\end{align}
these present $\braid_d(\CC^\times )$ in terms of the generators
$R,T_0,\dots ,T_{d-1}$. It is clear that in the braid class group
$R^d$ comes from a loop in $\CC^\times \subset\Homeo^0(\CC^\times)$
(the image of $R$ corresponds to multiplication by
$e^{2\pi\sqrt{-1}/d}$).
So it dies in $\brcl_d(\CC^\times )$, and indeed,
$\brcl_d(\CC^\times )$ is gotten from $\braid_d(\CC^\times )$ by
imposing this extra relation.  

The loop defined by $R^d$ gives the nontrivial element of
$\pi_1(\PSL(2,\CC ))\cong\ZZ /2$. So $R^{2d}$ dies in $\braid_d$.
The reader may check that in $\braid_d$ we also have the relations 
\begin{equation}\label{braidrel3}
R\equiv T_1T_2\cdots T_{d-1},\quad R^{-1}\equiv T_{d-1}T_{d-2}\cdots
T_1.
\end{equation}
One can verify that the relations (\ref{braidrel1}) imply that
$T_1T_2\cdots T_{d-1}$ and $T_{d-1}T_{d-2}\cdots T_1$ have the same  
$d$th power in $\braid_d(\CC^\times )$. So the relations
(\ref{braidrel3}) already imply that 
$R^{2d}$ maps to $1$ in $\braid_d$. 
Conjugating them with $R$ shows that the images of
$T_1T_2\cdots T_{d-1}$ and $T_{d-1}T_{d-2}\cdots T_1$ in $\braid_d$
are invariant under the cyclic permutation $(0,1,\dots ,d-1)$. 
(By suppressing $R$ and adding 
the cyclic invariance we get a presentation of $\braid_d$ in terms
of the $T_i$'s. The cyclic invariance also allows us to eliminate
another generator and this then leads to a presentation due to
Fadell-Van Buskirk in \cite{fadell}.)  
Finally, the braid class group $\brcl_d$ 
is gotten by putting $R^d\equiv 1$.  
 
\subsection{Action of a centrally extended braid class
group}\label{braidaction}
We  continue with the situation of Section \ref{cycliccov}.
We use the presentation of the braid class group $\brcl_{12}$ 
with generators $R=\eta, T_0,\dots ,T_{11}$ subject to the relations
(\ref{braidrel1}), (\ref{braidrel2}), (\ref{braidrel3}) (and
$\eta^{12}=1$).
The loop defining $T=T_0$ can be represented by a homeomorphism
of the pair $(\PP^1,\mu_{12})$ with support in a neighborhood $U$ of
the arc from $1$ to $e^{2\pi\sqrt{-1}/12}$, the loop defining $R$ is
represented by $\eta$. This homeomorphism lifts uniquely over 
$\pi :C_o\to\PP^1$ to a homeomorphism with support in $\pi^{-1}U$.
Let $\That$ denote its isotopy class in the group of homeomorphisms
of $C_o$ that commute with the $C_6$-action. 
(Perhaps we should remark that $\That$ is also
the monodromy that we get from a Milnor fibration: if we let the two
points of ramification  $1,e^{2\pi\sqrt{-1}/12}$ coalesce along the
segment  that connects them, then the $C_6$-covers acquire a
singularity with
local equation $w^6+\zeta^2$ (an $A_5$-singularity) and $\That$ is
the monodromy of this degeneration.)

The action of $\That$ on $H_1(C_o)$ will be a
$\ZZ[C_6]$-linear automorphism that preserves the intersection
pairing. Hence $\That$ will also act on the $\Ocal$-module
$H_1(C_o)_\Ocal$
and preserve the Hermitian form $\psi$ defined in Section
\ref{symplsym}.  
Let us make these actions explicit in terms of Pham's basis.
A suitable representative $\That$ (in the given
isotopy class) will act on $1$-chains on $C_o$ with boundary
supported
by the $\pi$-preimage of $0$ and the 12th roots of unity.
Clearly, $\That$ will not affect the class of $\eta^i{\bf
e}$ if $i\not= 0,1\pmod{12}$.  It is also easily seen that
$\That$ maps the class of $\eta{\bf e}$ to that of ${\bf e}$.
On the other hand $\That{\bf e}$ will be represented by the path
which first  follows ${\bf e}$, stops
just before $1$, makes then a full counterclockwise loop around the
ramification point over $1$, then returns to a point over $0$, and
finally
follows a lift over the segment $[0,e^{2\pi\sqrt{-1}/12}]$. From this
description it follows that this path is as a $1$-chain homologous to
$(1-\tau +\tau\eta){\bf e}$. 

\begin{corollary}\label{monodromy}
The monodromy operator $\That$ acts on $H_1(C_o)$ as follows:
\[
\That(\eta^ie)-\eta^ie=
\begin{cases} 
-(1+\tau) e &\text{if } i=0,\\
e &\text{if } i=1,\\
\tau e &\text{if } i=-1,\\
0 &\text{otherwise.}
\end{cases}
\]
It is in particular of order $6$.
Its action on $H_1(C_o)_\Ocal$ is the given by the complex reflection
\[
\That_\Ocal (x)=x-\tfrac{1}{3}(1+\omega)\psi (x,e)e
=x+\omega^{-1}\theta^{-1}\psi (x,e)e
\]
of order $3$.
\end{corollary}
\begin{proof}
The first statement follows in a straightforward manner from the
fact that $e=(1-\tau )(1-\eta ){\bf e}$, our computation of
$\That(\eta^i{\bf e})$, and the $\ZZ[C_6]$-linearity of $\That$. The
second follows from the first if we bear in mind the Formulae
\ref{psi} for 
$\psi (\eta^ie,e)=\psi (e,\eta^{-i}e)$.
\end{proof}

Consider the mapping class group $\brclhat_{12}$ of 
$C_6$-equivariant isotopy classes generated by $\That$ and $C_6\times
C_{12}$. So $\brclhat_{12}$ is a central extension of 
$\brcl_{12}$ by $C_6$. Let $\That_k:=\eta^k\That\eta^{-k}\in 
\brclhat_{12}$, $k\in\ZZ /12$.
These elements also obey the braid relations 
\begin{align}
\begin{split}
\That_k\That_{k+1}\That_k &= \That_{k+1}\That_k\That_{k+1}, \quad
k\in\ZZ
/12,\\
\That_k\That_l &= \That_l\That_k, \quad k,l\in\ZZ /12, k-l\not=\pm 1.
\end{split}
\end{align}
In view of the relations (\ref{braidrel3}) it is natural to put
\[
\Rhat:=\That_1\That_2\cdots \That_{11},\quad 
\Rhat^*:=\That_{11}\That_{10}\cdots \That_1.
\]
\begin{lemma}
We have $\Rhat =\tau\eta$ and $\Rhat^*=\eta^{-1}$.
\end{lemma}
\begin{proof}
From the definitions we find that
$\Rhat =(\eta \That)^{11}\eta$ and
$\Rhat^*=\eta^{-1}(\That\eta^{-1})^{11}$.
We know a priori that $(\eta T)^{11}$ and $(\That\eta^{-1})^{11}$ are
covering transformations, hence 
it is enough to show that these elements act on $H_1(C_o)$ as 
resp.\ $\tau$ and $1$.
This is verified in a straightforward manner using Corollary
\ref{monodromy}.
\end{proof}

So the $\That_i$'s generate all of $\brclhat_{12}$.
It also follows that $\brclhat_{12}$ 
is a nontrivial central extension of $\brcl_{12}$.

Recall from Section \ref{cycliccov} that we identified
$H_1(C_o)_\Ocal$ 
with the hermitian rank $10$ $\Ocal$-module $\Lambda$.
We noted in Remark \ref{isotropic} (see also the more precise
identification in the Appendix) one finds that the 
form $\psi$ on $\Lambda$ has hyperbolic signature $(9,1)$. 
Since the action of $\brclhat_{12}$ in $H_1(C_o)$ preserves
the $\ZZ [C_6]$-module structure and the sesquilinear form, we have
an induced  monodromy representation $\brclhat_{12}\to\Uu$
with $\Rhat\Rhat^*$ mapping to $\omega$. This drops to a projective
representation $\brcl_{12}\to \PU(\Lambda )$.

\begin{theorem}[Allcock, \cite{allcock}]
The monodromies 
\[
\rho :\brcl_{12}\to \PU(\Lambda )\quad\text{  and   }
\rhohat :\brclhat_{12}\to \Uu
\] 
are surjective.
\end{theorem}

\begin{corollary}\label{surj}
Every unitary automorphism of $\Lambda$ comes from a
$C_6$-equivariant
symplectic automorphism of $H_1(C_o)$.
\end{corollary} 

It follows from \ref{monodromy} that in either case the
image of $\hat T_i$ has order three. So if we define
$\brclhat_{12}[3]$ as the quotient of $\brclhat_{12}$ by the
relations
$\hat T_i^3\equiv 1$ and define $\brcl_{12}[3]$ similarly, then the 
monodromy representations factorize over homomorphisms
$\brclhat_{12}[3]\to \U (\Lambda )$ and $\brcl_{12}[3]\to \PU(\Lambda
)$. We shall see that these are isomorphisms.

\section{Satake-Baily-Borel compactification}\label{bb}

Let $V$ be a complex vector space equipped with a Hermitian form
$\psi :V\times V\to \CC$ of hyperbolic signature $(n,1)$, with $n\ge 
2$. Denote by $\LL =\LL (V)\subset V$ the set of $v\in V$ with $\psi 
(v,v)<0$. Then its projectivization $\BB =\BB (V)\subset\PP (V)$ (a complex 
ball) is a symmetric space for the projective unitary $\PU (V)$.
We regard $\LL$ as an equivariant $\CCtimes$-bundle over $\BB $.
For any integer $k$ we denote by $\LL (k)$ the line bundle defined
by the representation of $\CC^\times$ on $\CC$ given by 
$z\in\CC^\times\mapsto z^k$.
Then $\LL (2)$ is equivariantly isomorphic to the canonical bundle of
$\BB $. (To see this, observe that if $p\in\BB $ is given by the 
negative definite line $L\subset\Lambda_\CC$, then the tangent space
of $\BB $ at $p$ is canonically isomorphic to $\Hom (L,\Lambda_\CC
/L)$ and hence the determinant line of the cotangent space with 
$L^{n+1}\otimes \det (\Lambda_\CC)^{-1}$.) So the canonical bundle of $\BB$
is $\SU (V)$-equivariantly isomorphic to $\LL (n+1)$.

\begin{none}
Suppose $V$ has also the structure of a vector space over an imaginary
quadratic number field $K=\QQ (\sqrt{-d})$ in $\CC$ ($d$ a positive 
square free integer), such that $\psi$ is defined over $K$ and let be
given an arithmetic subgroup $\G$ of $\U (V_K)$. Then $\G$ acts
properly on the $\CCtimes$-bundle $\LL $ and the analytic orbifold
\[
\LL _\G:=\G\bs \LL .
\] 
retains a $\CC^\times$ action.  

The space of {\it $\G$-automorphic forms} of weight $k$ is by
definition
\[
A^k:=H^0(\BB , \LL (k))^{\G}.
\]
Its elements may be thought of as $\G$-invariant functions on 
$\LL$ that are homogeneous of degree $-k$ on every fiber. 
The space $A^k$ is known to be finite
dimensional for all $k\in\ZZ$ and trivial for $k<0$. 
Observe that $A^k=0$ when $k$ is not divisible by the order of
$\G\cap K^\times$. (In the case that interests us this order will
be $6$.) 
Examples of such forms are the Poincar\'e  series: if
$v_0\in\LL $, then
\[
F(a):=\sum_{\gamma\in\G} \psi (a,\gamma v_0)^{-k}
\]
converges uniformly on compact subsets of $\LL $, provided that
$k\ge 2\dim (\Lambda_\CC)=2n+2$. Hence $F$ defines an element of
$A^k$. The direct sum 
\[
A^\dot :=\oplus_{k\ge 0} A^k
\]
is a $\CC$-algebra of regular functions on $\LL _\G$. 
It is an algebra of finite type whose spectrum we 
denote by $\LL ^*_\G$. This is a normal affine variety which contains $\LL _\G$ 
as an open-dense subvariety; we therefore call it the 
{\it automorphic hull} of $\LL _\G$. The group $\CC^\times$ acts 
on $\LL ^*_\G$ with a unique fixed point. The corresponding projective variety at infinity, $\proj (A^\dot )$, will be
denoted by $\BB ^*_\G$. As the notation suggests, the underlying spaces
are in fact orbit spaces of a $\G$-space extensions $\LL^*\supset\LL$
and $\BB^*\supset\BB$. The Satake-Baily-Borel theory constructs these
spaces and we briefly recount how this is done. 

A point of $\partial\BB$ defined over $K$ is called a {\it cusp}
(of the form $\psi_K$). Then the union $\BB^*$ of $\BB$ and the
set of cusps is just the convex hull of the $K$-points of the closure
of $\BB$ in $\PP (V)$. A nonzero isotropic vector $n$ defined over
$K$ defines a cusp $[n]\in \BB^*$ and conversely, a cusp defines an 
isotropic line $I\subset V$ defined over $K$. 
For such a line $I$, let  
\[
\pi_{I^\perp}: V\to V/I^\perp
\]
denote the obvious projection. If $n\in I$ is a generator, then
$\psi (\; ,n)$ defines a coordinate for $V/I^\perp$, so that
$\pi_{I^\perp}$ is basically given by the inner product with $n$.
The image of $\LL\subset V$ is the set of generators
$(V-I^\perp)/I^\perp$. 
Let $\LL^*$ be the disjoint union of $\LL $, the 
punctured lines $\pi_{I^\perp} (\LL )$ and a singleton $V/V=\{ *\}$. Notice
that $U(V_K)$ acts naturally on this union (with $*$ as fixed point). We
give $\LL^* $ the topology generated by 
\begin{enumerate}
\item[(i)] the open subsets of $\LL $, 
\item[(ii)] unions $\Omega_n\cup \pi_{n^\perp}(\Omega_n)$ with 
$n\in V_K-\{ 0\}$ isotropic and 
$\Omega_n\subset\LL$ is the subset defined by $-\psi (z,z)> |\psi
(z,n)|^2$. 
\item[(iii)] unions $\Omega_N\cup 
\cup_{n\in N} \pi_{n^\perp}(\Omega_N)\cup\{ *\}$ with 
$N\subset V_K-\{0\}$ a  finite union of $\G$-orbits of isotropic
vectors and $\Omega_N$ the subset of $\LL$ defined by the 
inequalities  $-\psi (z,z)> |\psi (z,n)|^2$, for all $n\in N$. 
\end{enumerate}
The group $\U (V_K)$ acts on $\LL^* $ as a group of
homeomorphisms. The action of the central subgroup $K^\times$ extends
in an  obvious way to $\CC^\times$ so that in fact $\U
(V_K).\CC^\times$
acts. The orbit space $\G\bs\LL ^*$ is the $\CC^\times$-space
underlying 
the automorphic hull (it is not difficult to verify that the
Poncar\'e series defined above extends continously to $\LL ^* $). The
cuspidal lines define finitely many (regular) $\CC^\times$-orbits in
$\LL^*_\G$, because $\G$ acts with finitely many orbits in the set of
cusps. 

Similarly, the space underlying $\BB ^*_\G$ is the
the $\CCtimes$-orbit space of the $\G$-orbit space of $\BB^*$ endowed
with the {\it horoball} (or {\it Satake}) topology:  
this is the topology of $\BB^*$ generated by 
\begin{enumerate}
\item[(i)] the open subsets of $\BB$, 
\item[(ii)] unions $\PP (\Omega_n) \cup \{ [n] \}$, where 
$\PP (\Omega_n)\subset\BB$ is of the form 
$-\psi (z,z)/|\psi (z,n)|^2>1$, with $n\in V_K$ nonzero isotropic. 
\end{enumerate}
\end{none}

\begin{none}\label{symmhypersurface}
The automorphic hull possesses plenty of totally geodesic
hypersurfaces:
Suppose that $\Hcal$ is a $\G$-invariant collection of
$K$-hyperplanes in
$V$ of hyperbolic signature (that is, orthogonal to a positive
vector).
We assume that $\G$ has finitely many orbits in this collection. 
An example is the case when $\Hcal$ is the set of hyperplanes that
are perpendicular to a vector $v\in V_K$ with $\psi (v,v)=k$ ($k$ a
fixed
positive integer). For every $H\in\Hcal$, $\BB (H)$, is totally
geodesic subball of $\BB$ and the collection of these is 
locally finite on $\BB$. So 
\[
\LL (\Hcal):=\cup_{H\in\Hcal} \LL (H)
\]
is closed in $\LL$ and defines a closed analytic subset $\LL
(\Hcal)_\G$
of $\LL_\G$. If $n>2$, then $\LL^*_\G -\LL_\G$ is of 
codimension $>2$ in $\LL^*_\G$ and an
extension theorem implies that the closure $\LL (\Hcal)_\G^*$ of $\LL
(\Hcal)_\G$ is
analytic in $\LL^*_\G$. 
(This is also true when $n=2$, but that needs an additonal argument.)

This will be a $\CC^\times$-invariant hypersurface, hence algebraic. 
Notice that $\LL (\Hcal)_\G^*$ supports an effective Cartier divisor
if and only if $\LL (\Hcal)$ is defined by a single automorphic form.
(That form then will admit a product expansion.)
\end{none}

\section{The moduli space of rational curves with $12$ punctures}
\label{modulispace}

By a {\it smooth $C_6$-curve} we will mean a complete nonsingular 
complex-projective curve $C$  endowed with
an action of the cyclic group $C_6$ that is
isomorphic to a curve $C_D$ with affine equation
$w^6=\prod_{p\in D} (z-p)$, where $D$ is a $12$-element subset of
$\CC$, with $\tau (w,z)=(\omega w,z)$ (recall that
$\tau$ is a fixed generator of $C_6$ and $\omega = e^{2\pi\sqrt{-1}/6}$).
A more intrinsic characterization is to say that $C$ has genus $25$ and
that the $C_6$-action has $12$ distinct 
fixed points, each with (tangent space) character $\chi$, and 
is free elsewhere. (The Riemann-Hurwitz formula shows that 
its orbit space is then a rational curve.)   

Given such a smooth $C_6$-curve $C$, let $H^{1,0}(C)_\chi$ denote the
space of regular differentials $\alpha$ on $C$ on which $C_6$
acts with character $\chi$, that
is, which satisfy $\tau^*\alpha=\omega^{-1}.\alpha$.
We claim that $H^{1,0}(C)_\chi$ has dimension one. To see this,
represent $C$ by an affine equation $w^6=\prod_{p\in D} (z-p)$
as above. Then $w^{-1}dz$ is a regular differential on $C$ and 
$\tau^*(w^{-1}dz)=\omega^{-1}.w^{-1}dz$. Notice that the only zeroes
of $w^{-1}dz$ are the ramification points and that each 
such point appears with multiplicity $4$. This implies that it is the
only such form up to scalar: any other must be of the form
$f(z)w^{-1}dz$ with $f$ a rational function. In order that it be regular
$f$ should have no poles, so $f$ must be constant. 
(If we let the ramification points move in $\PP^1$, then
a period of such a form is a Lauricella function, see
\cite{delmost1}.)

The coarse moduli space of the  $C_6$-curves under consideration is
the same as the one of $12$ element subsets of a projective line 
(given up to a projective transformation), and so can be identified
with $\Dcal$. This suggests to allow as singular objects
the $C_6$-coverings of a projective line $C_D\to P$ with $D$ a semistable
divisor on $P$ such that over a point of multiplicity $k$ of $D$ we have 
a (plane curve) singularity with local equation $z^k=w^6$ ($k=1,\dots ,6$. 
A good substitute for the sheaf of regular differentials is then
the dualizing sheaf $\omega_C$. 

\begin{lemma}\label{eigenwaardemult}
For a $C_6$-covering $C=C_D\to P$ with $D$ semistable, the 
$\chi$-eigenspace in $H^0(C,\omega_C)$
is one-dimensional. The pull-back of a generator to a normalization
of $C$ is a logarithmic differential whose polar set is the
preimage of multiplicity $6$ locus of the discriminant. 
\end{lemma}
\begin{proof}
Choose an affine equation for $C$ as before.
First note that $w^{-1}dz$ lies in $H^0(C,\omega_C)_\chi$. 
At a point of multiplicity $k$, a local equation of $C$ is 
$z^k=w^6$. A straightforward calculation
shows that the pull-back of $w^{-1}dz$ under normalization has
in each the preimage of this singularity a zero of order $4,1,0,0,0,-1$
for $k=1,2,3,4,5,6$. Any other element of $H^0(C,\omega_C)_\chi$ is of
the form $f(z)w^{-1}dz$ and as in the smooth case we find
that $f$ cannot have any poles, hence must be constant.
\end{proof}

\begin{lemma}\label{eigenlijnbundel}
The orbifold line bundle $\Lcal_{\Dcal^*}$ over $\Dcal^*$ 
is naturally isomorphic to the coarse moduli space of pairs
$(C,\alpha^{\otimes 6})$ with  $C$ a $C_6$-curve with semistable
discriminant divisor and $\alpha\in H^0(C,\omega_C)_\chi$. 
\end{lemma}
\begin{proof}
We use our fixed two dimensional vector space $\Pi$ equipped with a
generator $\zeta$ of $\wedge^2\Pi $. Given a semistable
$F\in \Pi_{12}$, regard $F$ as a
homogeneous function on $\Pi ^*$. Then $w^6=F$ defines a degree $6$ 
covering of $\Pi ^*$. It is an affine
surface with good $\CC^\times$-action (so that $w$ has weight $2$)
whose curve at infinity is a $C_6$-curve
$C$ as above. Then $w^{-1}\zeta$ is a $\CC^\times$-invariant
rational form whose residue at infinity, $\alpha$, is a nonzero
element of $H^0(C,\omega_C)_\chi$. So $\alpha^{\otimes 6}$ is the
residue of $w^{-6}\zeta^{\otimes 6}=F^{-1}\zeta^{\otimes 6}$.
Think of $F^{-1}$ 
as the linear form on the line $\CC F$ in $\Pi_{12}$ spanned by $F$
which takes the value $1$ on $F$. The $\SL (\Pi)$-orbit of such a
linear form defines an element of the complement of the zero section
of $\Lcal_{\Dcal^*}$ and vice versa. Since the constructions are 
$\SL (\Pi )$-equivariant, we thus get a map from the complement of
the zero section of $\Lcal_{\Dcal^*}$ to the moduli space in
question. It is easy to see this this extends to an isomorphism of
$\Lcal_{\Dcal^*}$ to the moduli space.
\end{proof}

\medskip
Let $C$ be a smooth $C_6$-curve as above. The
intersection pairing identifies $H^1(C)$ with $H_1(C)$ as
$\ZZ [C_6]$-modules with symplectic form. Since $H_1(C)$ is 
isomorphic (as a $\ZZ [C_6]$-module with symplectic form) to $A$,
the choice of such an isomorphism induces an
isomorphism of Hermitian $\Ocal$-modules $\Lambda =A_\Ocal\to
H^1(C)_\Ocal$.  We shall refer to a Hermitian isomorphism 
$\Phi :\Lambda\to H^1(C)_\Ocal$ as a {\it $\Lambda$-marking} 
of the $C_6$-curve $C$. By Corollary \ref{surj} such a marking always
comes from a sesquilinear isomorphism $A\to H^1(C)$.

\begin{lemma}
The automorphism group of the $C_6$-curve $C$ acts faithfully on
the quotient $H^1(C)_\Ocal$.
\end{lemma}
\begin{proof}
This is clear for the group of covering transformations.
Any such automorphism that is not a covering transformation must
permute the ramification points nontrivially. It is easy to see that
such an automorphism acts nontrivially on $H^1(C)_\Ocal$.
\end{proof}

This implies that a $\Lambda$-marked $C_6$-curve has no
automorphisms. Hence there is fine moduli space $\Dcaltilde$
in the analytic category of these objects. It is an analytic manifold
of dimension $9$ (use three of the ramification points
as coordinates for the projective line $C_6\bs C$; the other nine
then run over an open subset of $\CC^9$) and
comes with an evident action of the unitary group $\U (\Lambda )$ of 
$\Lambda$: $u\in\U (\Lambda )$ sends $(C,\Phi )$ to $(C, \Phi
u^{-1})$. This action is proper and the orbit space can be identified
with $\Dcal$. Lemma \ref{eigenlijnbundel} suggests we also consider 
the moduli space $\Lcal_{\Dcaltilde}^{1/6}$ of triples $(C,\Phi ,\alpha)$ 
consisting of a $\Lambda$-marked genus $C_6$-curve $(C,\Phi )$ and an 
element $\alpha\in H^0(C,\omega_C)_\chi$.
It is clear that the projection $\Lcal_{\Dcaltilde}^{1/6}\to \Dcaltilde$ is
a $\Uu$-equivariant line bundle.

\begin{lemma}\label{lemma:eindigevert}
The morphism $\Dcaltilde\to\Dcal\subset\Dcal^\st$ extends naturally
to a branched $\Uu$-covering
$\Dcaltilde^\st\to\Dcal^\st$. Moreover, the
$\Uu$-equivariant line bundle $\Lcal_{\Dcaltilde}^{1/6}\to\Dcaltilde$ 
extends  naturally to a $\Uu$-equivariant line bundle
$\Lcal_{\Dcaltilde^\st}^{1/6}\to \Dcaltilde^\st$.
\end{lemma}
\begin{proof}
Let $D$ be a stable effective  degree $12$ divisor in $\CC$ (so
all multiplicities $\le 5$). Given a neighborhood $U$ of 
$D$ in the space of effective  degree
$12$ divisors, denote by $U'\subset U$ the divisors that are
reduced. Then $D'\in U'\mapsto H^1(C_{D'})_\Ocal$ defines a 
locally constant sheaf of $\Ocal$-modules. If
$D$ has multiplicities $5\ge n_1\ge n_2\ge
\cdots n_r\ge 1$ (so that $\sum n_i=12$), and $U$ is sufficiently
small, then the local monodromy
group is isomorphic to subgroup of $\prod_i U(\Lambda^{n_i-1})$.
Since the ranks ${n_i-1}$ are all $\le 4$, the latter is finite
by Subsection \ref{lambda4}, and
hence so is the monodromy group. The assertions of the lemma are a 
formal consequence of this fact.
\end{proof}

\begin{remark}\label{rem:eindigevert}
Closer inspection shows that there is in fact a
moduli interpretation of the added points: an element of
$\Dcaltilde^\st$
is represented by a pair $(C,\Phi)$, where $C$ is a $C_6$-curve with
stable ramification divisor and $\Phi :\Lambda\to H^1(C)_\Ocal$ is
a certain epimorphism of $\Ocal$-modules. The kernel of $\Phi$ is
isomorphic to an orthogonal direct sum of sublattices 
$\Lambda^{n_1-1}\perp \Lambda^{n_2-1}\perp\cdots $ and
$\Phi$ is given up to composition with an element of the local
monodromy group
$\prod_i U(\Lambda^{n_i-1})$. A point of
$\Lcal_{\Dcaltilde^\st}^{1/6}
$ is obtained by also giving an element of $H^0(C,\omega_C)_\chi$. 
\end{remark}

\begin{remark}\label{norm}
If $D$ is stable, then we have a square
norm on $H^0(C,\omega_{C})_\chi$ defined by
\[
\alpha\in H^0(C,\omega_{C})_\chi\mapsto\theta\int_C
\alpha\wedge\bar\alpha.
\]
In case $D$ is reduced, then this is just the restriction of our Hermitian
form $-\psi$ via the embedding  
\[
H^0(C,\omega_C)_\chi\subset H^1(C;\CC )_\chi =\CC\otimes_\Ocal
H^1(C)_\Ocal.
\]
This norm blows up over the point $d_\infty$. To see this, 
use the fact that if $D$ becomes strictly semistable, then
$w^{-1}dz$ becomes a differential on the normalization of $C_D$
with poles of order one. So the integral of the generating section
defined by $|w^{-1}dz|^2$ blows up over $d_\infty$.
\end{remark}

We now define a period mapping. Let 
$(\Lcal_{\Dcaltilde}^{1/6})^\times$ be 
the complement of the zero section of $\Lcal_{\Dcaltilde}^{1/6}$.
Let $\LL$ be as defined in Section \ref{bb} with
$V=\Lambda_\CC$, $K=\QQ (\omega )=\QQ (\sqrt{-3})$ and $\Gamma
=\Uu$. If $(C,\Phi ,\alpha)$ represents 
a point of $(\Lcal_{\Dcaltilde}^{1/6})^\times$, then assign to this triple
the vector $\Phi^{-1}(\alpha )$. This defines the {\it period
mapping}:
\[
\pertilde :(\Lcal_{\Dcaltilde}^{1/6})^\times\to \LL.
\]
This mapping is clearly equivariant with respect to the actions of
$\CC^\times$ and $\Uu$ and both its domain and range are 
analytic manifolds of dimension $10$. This
period mapping extends across the locus with finite
monodromy: we have an extension
\[
\pertilde :(\Lcal_{\Dcaltilde^\st}^{1/6})^\times\to \LL.
\] 
Indeed, if a point of the domain is represented as in 
Remark \ref{rem:eindigevert} by a triple $(C,\Phi ,\alpha)$, then
Lemma \ref{eigenwaardemult} implies that $\alpha$ defines a nonzero
element
of $H^1(C)_\chi$  and the image of $(C,\Phi ,\alpha)$ is the point of
$\LL\cap\ker (\Phi)^\perp$ that is mapped by $\Phi$ to $\alpha$.
For the details we refer to \cite{delmost1}.
The period mapping drops to a morphism
\[
\per :(\Lcal_{\Dcaltilde^\st}^{1/6})\to\LL_\Uu,
\]  
and if we pass to $\CC^\times$-orbit spaces, we also get
\[
\PP (\pertilde ):\Dcaltilde^\st\to\BB \text{  and   } \PP(\per ):
\Dcal^\st\to \BB_\Uu.
\]
The following theorem is a special case of a theorem of 
Deligne-Mostow \cite{delmost1}. 

\begin{theorem}[Deligne-Mostow \cite{delmost1},
see also \cite{couw}]\label{torelli:1}   
The period map $\pertilde$ establishes a $\Uu$-equivariant isomorphism
between the $\CC^\times$-bundle
$(\Lcal_{\Dcaltilde^\st}^{1/6})^\times$ and $\LL$. The induced isomorphism 
$\Dcal^\st\to \BB_\Uu$ extends to an isomorphism between the GIT
compactification
$\Dcal^*\supset\Dcal^\st$ and the Baily-Borel compactification
$\BB_\Uu^*\supset \BB_\Uu$. 
\end{theorem}

Statement and proof are somewhat hidden in the paper and so we give an

\begin{proof}[Outline of proof] Since
$\pertilde$ is $\CC^\times$-equivariant, it is enough to prove that
$\PP (\pertilde ):\Dcaltilde^\st\to\BB$ is an isomorphism.
To this end, one first shows that $\PP (\pertilde )$ is a local
isomorphism in codimension one (this is based on simple type of
local Torelli theorem) and has discrete fibers. This implies
that $\PP (\pertilde )$ has no ramification. So $\PP (\pertilde )$
is a local isomorphism every where. We wish to show
that $\PP (\pertilde )$ is proper; the simple connectivity
of $\BB$ will then imply that $\PP (\pertilde )$ is an isomorphism. 
This will follow if we prove
that $\PP (\per ): \Dcal^\st\to \BB_\Uu$ is proper. In other words,
we want to show that $\PP (\per )$ extends continuously to
the one-point compactifications of its domain and range.

Let $D$ be a strictly semistable divisor of degree $12$ on $P=\PP^1$.
So $D$ has a point of multiplicity $6$. 
Let $\g$ be a small oriented circle around this point.
Then the preimage of $\g$ in $C_D$ consists of $6$ disjoint circles.
If $\tilde \g$ is one of these, then $\int_{\tilde\g } w^{-1}dz$ is 
by \ref{eigenwaardemult} the residue
of a differential with a simple pole
and hence nonzero. The cycle $\tilde \g$ subsists
under small deformations of $D$ and for $D'$ in a neighborhood of $D$
the corresponding integral $\int_{\tilde\g (D')} w^{-1}dz$ is then
analytic in $D'$ and nowhere zero. 
If $D'$ is reduced, then $\tilde\g (D')$ defines an isotropic 
element of $H_1(C_{D'})_\Ocal$. On the other hand, by Remark
\ref{norm},  $\int_{C_{D'}} |w^{-1}dz |^2$ tends to $+\infty$, as
$D'$ approaches $D$. So the same is true for the expression 
\[
\frac{\int_{C_{D'}} |w^{-1}dz |^2}{|\int_{\tilde\g(D')} w^{-1}dz |^2}.
\]
It now follows from our explicit description of the Satake topology 
in Section \ref{bb} that the image of $D'$ under $\PP (\per )$ tends to
the cusp of $\BB_\Uu^*$, as $D'$ tends to $D$.
This proves that $\PP (\per ): \Dcal^\st\to \BB_\Uu$ is proper.

So both $\PP (\pertilde): \Dcaltilde^\st\to \BB$ and
$\PP (\per ): \Dcal^\st\to \BB_\Uu$ are isomorphisms.
Since $\Dcal^*$ and $\BB_\Uu^*$
are normal one point compactifications of $\Dcal^\st$ and
$\BB_\Uu$ respectively, the continuous extension 
$\Dcal^*\to\BB_\Uu^*$  is in fact an isomorphism.
\end{proof}

We can also tell what the image of $\Dcal$ is.
Let us call a hyperplane in $\Lambda$
a {\it mirror} if it is the orthogonal complement of a $3$-vector. 
A mirror has hyperbolic signature and by Lemma \ref{trans6}
any two mirrors are $\Uu$-equivalent. So the collection $\Hcal$ of
mirrors defines an irreducible hypersurface $\BB (\Hcal )^*_\Uu$ in
$\BB^*_\Uu$. If we let of $12$ distinct points in $\PP^1$ two
coalesce, then we get a curve germ in $\Dcal^*$ with generic point
in $\Dcal$ and closed point the generic point $\Dcal^*-\Dcal$.
Associated to this there is a `vanishing $3$-vector' which shows
that $\Dcal^*-\Dcal$ is mapped to $\BB (\Hcal )^*_\Uu$.
Since both $\Dcal^*-\Dcal$ and $\BB (\Hcal )^*_\Uu$ are irreducible
we find:

\begin{theorem}\label{torelli:2}
The period mapping defines an isomorphism 
\[
(\Dcal^*,\Dcal^{st},\Dcal)\cong
(\BB^*_\Uu, \BB_\Uu, \BB_\Uu-\BB(\Hcal)_\Uu ).
\]
\end{theorem}

\begin{remark}
We observed in \ref{discriminant} that the discriminant hypersurface 
$\Dcal^*-\Dcal$ has degree $11.12$ (with respect to the
$\Lcal_{\Dcal^*}$). Hence the locally symmetric hypersurface 
$\BB (\Hcal)^*_\Uu $ is defined by a section of $\LL (6.11.12)$. 
Since $\LL\to\LL_\Uu $ ramifies with order three along $\LL (\Hcal)$, it
follows that the divisorial preimage of $\LL (\Hcal)_\Uu$ 
is $3\LL (\Hcal)$. So $\LL (\Hcal)$ is given by an
automorphic form of weight $2.11.12$ with a character of order $3$.
Since Allcock finds this degree to be $44$ \cite{allcock}, we
assume that his weight  is $1/6$ of ours (the center of $\Uu$ 
consists of the $6$th roots of unity and so the degree of any 
nonzero $\Uu$-automorphic form on $\BB$ is divisible by $6$). 
\end{remark}

\begin{corollary}
The kernel of the monodromy representation $\rho
:\brclhat_{12}\to\Uu$
is the normal subgroup generated by $\That_0^3$ so that $\rho$
induces isomorphisms $\brclhat_{12}[3]\cong \G $
and $\brcl_{12}[3]\cong \PU (\Lambda )$.
\end{corollary}
\begin{proof}
The group $\brcl_{12}$ may be identified with the orbifold
fundamental 
group of $\Dcal$. Via the orbifold isomorphism $\Dcal\cong\BB_\Uu-\BB
(\Hcal)_\Uu$, we 
then get a $\brcl_{12}[3]$-covering. This covering 
factorizes over a covering of $\BB -\BB (\Hcal)$ with  
the kernel of $\brcl_{12}[3]\to \PU (\Lambda )$ as covering group.

Since $\That_0^3$ is trivial in $\brclhat_{12}[3]$, 
a simple loop around a deleted hyperplane has 
monodromy of order three, and so the covering over $\BB -\BB (\Hcal)$
extends as an {\it unramified} covering over the smooth part of 
$\BB (\Hcal)$: we now have a connected  unramified covering over
$\BB -\BB (\Hcal)_\sing$. Since $\BB -\BB (\Hcal)_\sing$ is simply
connected, this covering must be trivial. We conclude that
$\brcl_{12}[3]\to \PU (\Lambda )$ is injective. From this it follows
that $\brclhat_{12}[3]\cong \Uu $ is injective as well.
\end{proof}

\section{Rational elliptic surfaces and the Eisenstein
curve}\label{ellfibreis} 
Recall from our discussion of Kodaira's theorem \ref{sub:kodaira}
that the commutator subgroup of $\PSL (2,\ZZ)$ defines a modular
curve $E_o$ of genus one with a simple cusp. We regard it as an
elliptic curve by taking
the cusp as its origin. It comes with a faithful action the
abelianization $C_6$ of $\PSL (2,\ZZ)$, and so this elliptic curve
has $J$-invariant $0$. In other words, 
it can be analytically obtained as
the quotient $\CC/\Ocal$ with the generator $\tau\in C_6$
acting as complex multiplication by $\omega$. So 
$\ZZ [C_6]$ acts on $H^1(E_o)$ via $\Ocal$. 
We will refer to $E_o$ as the {\it Eisenstein
curve}. Since $\tau$ acts on the tangent space of
the origin with eigenvalue $\omega$, the same is true
for the action of $\tau^*$ on $H^{1,0}(C)$. It follows that
$H^1(C,\CC)_\chi =H^{0,1}(C)$.

The natural map to the $J$-line,
$E_o\to \PP^1$, ramifies over $0$ (two points of order three), $1$
(three points of order two) and $\infty$ (total ramification). 

\begin{lemma}\label{eisensteinquot:1}
Let $X\to P$ be a rational elliptic surface with reduced discriminant
$D_\infty$. Let $J: P\to \PP^1$ be its modular function and 
$C$ be the normalization of $P\times_{\PP^1} E_o$. Then
the $C_6$-covering $C\to P$ is the one considered in Section 
\ref{modulispace}:
it is only ramified over $D_\infty$, the ramification over $D_\infty$
is total and $\tau$ acts in the tangent space of each ramification
point
as multiplication by $\omega$.
\end{lemma}
\begin{proof} It is clear
that the projection $C\to P$ is a $C_6$-covering. There is no
ramification outside the discriminant divisor $J^*(\infty )$ since $J$ is
there locally liftable to a morphism to $E_o$.
The remaining statements follow easily.
\end{proof}

A special feature of this situation is that $C$ comes with a 
$C_6$-equivariant morphism $\tilde J :C\to E_o$. Its degree is
clearly $12$. 

\begin{theorem}\label{eisensteinquot:2}
In the situation of Lemma \ref{eisensteinquot:1} we have:
\begin{enumerate}
\item[(i)] The morphism $\tilde J :C\to E_o$ induces an embedding
$\tilde J^*_\Ocal: H^1(E_o)_\Ocal\to H^1(C)_\Ocal$ of 
$\Ocal$-modules that multiplies the hermitian form by $12$,
\item[(ii)] the line $H^{1,0}(C)_\chi$ is perpendicular to
the image of $\tilde J^*_\Ocal$ and
\item[(iii)] there exists a $6$-vector $z\in H^1(C)_\Ocal$ such
that the image of $\tilde J^*_\Ocal$ is the $\Ocal$-submodule
$H^1(C)_\Ocal$
spanned by $2\theta z$.
\end{enumerate}
\end{theorem}
\begin{proof}
The first assertion  follows from the fact that $\tilde J$ is 
$C_6$-equivariant and of degree $12$ and
the second from the observation that 
$H^1(E_o,\CC)_\chi =H^{0,1}(E_o)$.

The last clause requires more work. 
In view of the connectedness of $\Mcal$, it is enough to prove
that assertion for one particular rational elliptic surface. We take
the case studied in Section \ref{cycliccov}, where
$D_\infty\subset\PP^1$
is the set of $12$th roots of unity
and $C_o\to\PP^1$ is the curve with $C_6\times C_{12}$-action. As
noted in Example \ref{simplediscr}, $D_\infty$ is the discriminant divisor of
an elliptic surface, but we will exhibit such a fibration 
more directly. 
Consider the action of the (order 12)
subgroup $G\subset C_6\times C_{12}$ generated by $\tau^3\eta$. The
orbit space $G\bs C_o$ is a $C_6$-covering of $C_{12}\bs\PP^1$. If we
identify the latter with $\PP^1$ by means of the affine coordinate
$z^{12}$, then we see that $G\bs C_o\to \PP^1$ has total ramification
over $1$, a fiber with two points over $0$ 
and a fiber with three points
over $\infty$. These properties imply that $G\bs C_o$ has genus one
and more than that, namely that $G\bs C_o$ is  
$C_6$-equivariantly isomorphic to the Eisenstein curve $E_o$. The
Eisenstein curve supports a $C_6$-equivariant elliptic fibration.
This pulls back to a $C_6$-equivariant elliptic fibration over $C_o$
and that in turn descends to an elliptic fibration on $\PP^1$.
We therefore denote the resulting $C_6$-morphism 
$\tilde J:C_o\to E_o$. The induced map on the first cohomology
$\tilde J^*: H^1(E_o)\to H^1(C_o)$ is $C_6$-equivariant. We identify
the $\ZZ [\tau, \eta ]$-module $H^1(C_o)$ with the algebra $A$
defined in \ref{cycliccov}. It is clear that the image of $\tilde
J^*$ is the
$\Ocal$-submodule spanned by 
\[
\sum_{i=0}^{11} (\tau^3\eta )^i\in A.
\]
The image $u$ of this element in 
\[
H^1(C_o)_\Ocal\cong A_\Ocal =
\Ocal [\eta ]/(\sum_{i=0}^{11}\eta^i,\sum_{i=0}^{11}(\omega\eta )^i)
\]
is easily calculated to be of the form $2\theta z$, with 
\[
z= \omega ^{-1}(\eta^2 +\eta^8)+ 
(\eta^3 + \eta^4 + \eta^9+\eta^{10}) + \omega (\eta^5 +\eta^{11}).
\]
We claim that $u$ is a $12.6$-vector: this is a straightforward 
computation or one invokes Example \ref{eisensteincurve} and the fact
that the
Hermitian form is multiplied by $12$. So $z$ is a $6$-vector.
\end{proof}

The last assertion of the above proposition implies that 
the condition for a $12$ element subset of $\PP^1$ to be the
discriminant
of a rational elliptic surface {\it imposes a linear constraint on 
the period map} defined in Section \ref{modulispace}. We investigate
this in more detail in the next section.

\section{Moduli of rational elliptic surfaces II}\label{modulispace2}

From now on, we make free use of notions, notation and
results of the theory of $\Ocal$-lattices, as collected and proved 
in the Appendix.

\bigskip
In the Appendix we fix a sublattice
$\Lambda_o$ that is the orthogonal complement of a $6$-vector $z_o\in
\Lambda$. (It is proved in Proposition \ref{trans6} that all such 
sublattices are $\Uu$-equivalent.) According to Proposition 
\ref{unitairiso} the stabilizer of $\Lambda_o$ in $\Uu$ restricts 
isomorphically to the unitary group $\Ulo$ of $\Lambda_o$. 
It follows from Proposition \ref{trans6} that $\Ulo$
has two orbits in the set of primitive $0$-vectors in $\Lambda_o$:
type $(\theta)$ and $(0)$. So the Baily-Borel compactification 
$\BB^*_{o,\Ulo}$ adds two points to $\BB_{o,\Ulo}$. We denote
them  $\infty_\theta$ and $\infty_0$.

We call a hyperplane $H$ of $\Lambda_o$
a {\it mirror trace} if it is the intersection of
a mirror of $\Lambda$ with $\Lambda_o$ and has hyperbolic signature.
This amounts to requiring that the orthogonal complement $H^\perp$ of
$H$ in $\Lambda$ is positive definite and contains the $6$-vector
$z_o$ and a $3$-vector. According to Lemma \ref{class} the 
discriminant of $H^\perp$ then takes the values $6$, $9$, $15$ or
$18$; we denote that number by $d(H)$ and call it the {\it
$d$-invariant} of $H$. 
A special role  will be played by the mirror traces  with
$d$-invariant $6$ or $9$ as in these cases there exist $3$-vectors
$r_1$, $r_2$ in 
$H^\perp$ such that $r_1+r_2$ spans $\Lambda_o^\perp$. Proposition
\ref{dchar} can be  restated as:

\begin{proposition}\label{mirrortrace}
Two mirror traces with the same $d$-invariant are equivalent under
the $\Ulo$-action.
\end{proposition}

We denote the collection of  mirror traces by $\Hcal_o$, and those
with $d$-invariant in a subset $S\subset\{6,9,15,18\}$ by $\Hcal_o(S)$. 
So we get a hypersurface $\Delta:=\BB_o(\Hcal_o)_{\Uu_o}^*$ in $\BB_{o,\Ulo}^*$ 
that has four irreducible components:
$\Delta (d):=\BB_o (\Hcal_o(d))_{\Ulo}^*$, $d=6,9,15,18$. 

The inclusion $\LL_o\subset\LL$ induces a natural map 
\[
\LL_{o,\Ulo}^*\to\LL_\Uu^* 
\]
that is finite and birational onto a hypersurface of $\LL_\Uu^*$ (it
need not be injective though) so that  $\LL_{o,\Ulo}^*$ can be
identified with the normalization of this hypersurface. 
It is clear that $\LL_o (\Hcal_o)_{\Ulo}^*$ is the preimage of
$\LL (\Hcal )_\Uu^*$ under the map displayed above.

\medskip
Let $f:X\to P$ be rational elliptic surface with reduced
discriminant. We have an associated $C_6$-covering $C\to P$ together
with an equivariant morphism $C\to E_o$.  We say that a
$\Lambda$-marking $\Phi :H_1(C)_\Ocal\cong\Lambda$ is {\it adapted}
if $\Phi\tilde J^*$ maps $H^1(E_o)$ to the orthogonal complement of
$\Lambda_o$. Rational elliptic surfaces with adapted markings define
analytic
covers $\Mcaltilde$ and $\Ecal_{\Dcaltilde}|\Mcaltilde$ 
of $\Mcal$ and $\Ecal_{\Dcal}|\Mcal$ respectively, the latter with
Galois group $\Ulo$, the former with Galois group $\Ulo$ modulo its
scalars. The period map induces an equivariant morphism
$\Ecal_{\Dcal}|\Mcaltilde\to \LL_o$.
It follows from the preceding that this morphism is injective; in
fact from Lemma \ref{reducedunique}, Theorem \ref{torelli:2} and
Theorem 
\ref{eisensteinquot:2} we get:

\begin{theorem}\label{torelli:3}
The period mapping induces an isomorphism of arrows:
\[
\begin{array}{ccc}
(\Mcal^*,\Mcal)& \cong &
(\BB_{o,\Ulo}^*,\BB_{o,\Ulo}-\BB_o(\Hcal_o)_{\Uu_o})\\
\Big\downarrow &   &\Big\downarrow \\
(\Dcal^*,\Dcal ) &\cong &(\BB^*_\Uu,\BB_\Uu-\BB(\Hcal)_\Uu ).\\
\end{array}
\]
\end{theorem}

According to \ref{strata}, the boundary of $\Mcal$ in $\Mcal^*$
consists of four irreducible hypersurfaces of $\Mcal^*$:
$\Mcal^*(I_2)$,
$\Mcal^*(II)$, $\Mcal^*(I_9)$ and $\Mcal^*(I^*_{4,4})$,
whereas the irreducible components of $\Delta$ are
$\Delta (18)$, $\Delta (15)$, $\Delta (9)$, $\Delta (6)$.
The period isomorphism \ref{torelli:3} must set up a bijection
between these two sets. Something similar should hold for  
the strata $\Mcal^*(I_6)$ and $\Mcal (I_0^*)$ 
lying over the two cusps $\infty_0$ and $\infty_\theta$ of 
$\BB^*_\Uu$. We complete the picture by determining which goes to
which.

\begin{theorem}\label{torelli:4}
The period isomorphism maps the irreducible
components $\Mcal^*(I_2)$, $\Mcal^*(II)$, $\Mcal^*(I^*_{4,4})$,
$\Mcal^*(I_9)$ 
onto $\Delta (18)$, $\Delta (15)$, $\Delta (9)$, $\Delta (6)$
respectively. Moreover, the singletons $\Mcal^*(I_6)$ and
$\Mcal (I_0^*)$ are mapped to $\{\infty_\theta\}$ and
$\{\infty_0\}$ respectively.
\end{theorem} 

Before we begin the proof, we note that this theorem is
equivalent to the corresponding statements for $\Mcal^K$ (instead of
$\Mcal^*$), for by definition $\Mcal^*(F)$ is the image of
$\Mcal^K(F)$ under the modification $\Mcal^K\to \Mcal^*$. We will
prove the theorem in this form.

Let $(J:P\to\PP^1, D)$ represent a closed point of $\Mcal^K$ and
let $C\to P$ be the corresponding $\mu_6$-covering. 
Consider a deformation of $(J:P\to\PP^1, D)$ over a smooth curve
germ $(\DD ,o)$ with smooth generic fiber. After a finite base change
this is  covered  by a smoothing of $C$:
\[
\Ccal\to \Pcal\to \PP^1\times \DD,
\]
where the first morphism is the quotient by an $\mu_6$-action and
the second is of degree $12$.
We observed in \ref{sub:kontcomp} that there is a natural
$\mu_6$-equivariant morphism $\Ccal\to E_o$.
In a situation like this there is a standard procedure for
comparing the cohomology of the special fiber and the general fiber:
the pull-back of $\Ccal\to \DD$ over the universal cover 
$\tilde \DD^\times$ of $\DD^\times:=\DD-\{ o\}$, 
$\Ccal_{\tilde \DD^\times}\to
\tilde \DD^\times$, is homologically trivial and after a choice of an
adapted $\Lambda$-marking we get an isomorphism of $\Ocal$-modules  
$H^1(\Ccal_{\tilde\DD^\times})_\Ocal\cong \Lambda$ 
such that the image of $H^1(E_o)_\Ocal$
is a multiple of $z_o$. This gives rise to a period morphism
$\tilde\DD^\times\to \BB_o$. The inclusion $C\subset\Ccal$ is
a homological isomorphism, and hence the diagram 
$C\subset\Ccal\leftarrow\Ccal_{\tilde\DD^\times}$
induces a homomorphism of $\Ocal$-modules $\Lambda\to 
H^1(C)_\Ocal$ such the image of $\Ocal 2\theta z_o$ is mapped onto 
$H^1(E_o)_\Ocal$.

\begin{proof}[Proof of \ref{torelli:4}] Consider the case when the
closed
fiber represents a general point of $\Mcal^K(I_2)$, $\Mcal^K(II)$,
$\Mcal^K(I_9)$ or
$\Mcal^K(I^*_{4,4})$. The image of such a point in $\Dcal^*$ is a
semistable orbit of a degree $12$ divisor on $\PP^1$ 
of type $(2,1^{10})$, $(2,1^{10})$, $(3,1^9)$, $(2^2,1^8)$
respectively. So its image under the period isomorphism is
going to be perpendicular to a (primitive) sublattice $L$ of
$\Lambda$ of type $\Lambda^1$, $\Lambda^1$, $\Lambda^2$, 
$\Lambda^1\times\Lambda^1$ respectively. 
In the last two cases, the central component of $P_c$
is in $J^{-1}(\infty)$ and so the morphism $H^1(E_o)\to
H^1(C)\to H^1(C_c)$ will be zero. This implies that in these cases 
$L$ contains $z_o$. This shows that in terms of the notation of
Lemma \ref{class} $L$ is of type $\delta_6$ in 
in the $\Mcal^K(I_9)$-case and of type $\delta_9$ and in
the $\Mcal^K(I^*_{4,4})$-case.
So we then find a point of $\Delta (6)$ and $\Delta (9)$
respectively. 

We now show that for $I_2$ we cannot end up with a point of
$\Delta (15)$. Since we have a period isomorphism, it then will
follow that we must get a point of $\Delta (18)$ and that 
in the remaining case $II$ we get a point of $\Delta (15)$. 
We note that in the $I_2$-case, the lattice $L\cong \Lambda^1$ is
accounted for  by $H^1(C')_\Ocal$, where $C'$ is the 
irreducible component of $C$
that lies over $\infty$. Since the map $C\to E_o$ is constant on  
$C'$, it follows that $L\subset \Lambda_o$. It follows
that $L+\Ocal z_0$ is of type $\delta_{18}$. 
A priori this lattice might be imprimitive, but it certainly does
not contain a lattice of type $\delta_{15}$. 

We know that both $\Mcal^K(I_6)$ and  $\Mcal^K(I_0^*)$ map
to $d_\infty\in\Dcal^*$. So they will map to distinct cusps of 
$\BB_{o,\Ulo}^*$. Hence is enough to show that $\Mcal^K(I_6)$ maps
to  $\infty_\theta$: then $\Mcal^K(I_0^*)$ must necessarily map to 
the other cusp $\infty_0$. A similar argument as used for
$\Mcal^K(I_2)$ shows that a generic point of $\Mcal^K(I_6)$
is mapped to cusp of $\BB_o^*$ that is perpendicular to a sublattice
$L\subset\Lambda_o$ isomorphic to $\Lambda^5$. Then $I:=L\cap
L^\perp$ 
is a primitive isotropic line whose image in $\BB_o^*$ is the cusp in
question.
A primitive isotropic line of type $(0)$ is not perpendicular to a
lattice of type $\Lambda^5$, whereas one of type $(\theta)$ is. So
$\Mcal^K(I_6)$ maps to a cusp of type $\theta$.
\end{proof}

From Corollary \ref{bundelvgl} we deduce a description of the Miranda 
compactification in terms of automorphic forms:

\begin{theorem}\label{torelli:5}
The graded $\CC$-algebra of automorphic forms on $\BB_o$ with values in a 
tensor power $\LL (k)$ with arbitrary poles along the hyperball arrangement
$\BB_o(\Hcal_o(6,9))$ is zero in negative degrees and of finite type.
Its proj reproduces the Miranda compactifiation of 
$\BB_{o,\Ulo}-\BB_o(\Hcal_o(6,9))_\Ulo$.
\end{theorem}

This means that the hypersurface $\Delta (6)\cup \Delta (9)$ in 
$\BB_{o,\Ulo}$ can never be the zero set of an automorphic form, since 
the inverse of such a form would produce an element of the above algebra 
of negative degree. This is in contrast with $\Delta$ itself (see 
\cite{allcock}).

\begin{remark} An intersection of mirror traces in $\BB_o$ of $d$-invariant
$6$ or $9$ is by definition the orthogonal complement of a positive 
definite sublattice $L\subset\Lambda$ spanned by $z_o$ and $3$-vectors
of $d$-invariant $6$ or $9$.
According to Proposition \ref{inbclass} there are, apart from the mirror 
traces themselves, three types:
$(6,9)$, $(9,9)$ and $(6,9,9)$, in which cases $L$ is spanned by
$z_o$ and  $3$-vectors of the indicated $d$-invariant. It also
follows from Proposition \ref{inbclass} that each of these three types
represents a single $\Ulo$-orbit.
So these define irreducible subvarieties $\Delta (6,9)$, $\Delta (9,9)$, 
and $\Delta (6,9,9)$ of $\BB_{o,\Ulo}^*$ of codimension $2$, $2$ and
$3$ respectively. Using  Corollary \ref{strata} one identifies these 
subvarieties in $\Mcal^*$ as $\Mcal^{*}(I_8')$,
$\Mcal^{*}(I_8'')$ and $\Mcal^{*}(I_7)$ respectively.
\end{remark}

\section{Modification of the Baily-Borel compactification.}\label{bbmod}
Although this section is mostly of a descriptive nature, 
it may help to put our results into perspective: we outline an extension
of the Baily-Borel theory which produces the compactifications obtained here in an
algebro-geometrical setting in a canonical fashion. This is closely related to 
the construction described in \cite{looijeng}.

\subsection{Modifications defined by arrangements}\label{arrmod}
Suppose we are given a complex manifold $X$ of dimension $n$ and
a collection $\Hcal$ of smooth hypersurfaces of $X$ that is locally 
finite on $X$
and is \textit{arrangementlike}, in the sense that at each point of
$X$ there exist local analytic coordinates such that each $H
\in\Hcal$ passing through that point is given by a linear equation.
Denote by $D=\cup_{H\in\Hcal}H$ their union. There is a simple
and straightforward way to find a 
modification $\tilde X\to X$ of $X$ such that 
strict transforms of the members of $\Hcal$ get separated: 
if $D^{(k)}$ denotes the union of the codimension 
$k$ intersections of members of the $\Hcal$, then first blow up $D^{(n)}$, 
then the strict transform of $D^{(n-1)}$, and so on, finishing
with blowing up a strict transform of $D^{(2)}$:
\[
X=\tilde X_n \leftarrow \tilde X_{n-1}\leftarrow\cdots\leftarrow 
\tilde X_1=\tilde X.
\]
If we denote the strict transform of $H$ in $\tilde X_k$ by $\tilde H_k$, then
the collection $\{ \tilde H_k\}_{H\in\Hcal}$ is also arrangementlike
and has no intersections of codimension
$>k$. In particular, the $\{ \tilde H_1\}_{H\in\Hcal}$ are disjoint.
It is clear that the blowup is an isomorphism over $\Omega :=X-D$.

\begin{lemma}
The morphism $\tilde X\to X$ is obtained by 
blowing up the fractional ideal $\sum_{H\in\Hcal} \Ocal_X(H)$.  
\end{lemma}
\begin{proof} Let $k$ be the maximal integer for which 
$D^{(k)}$ is nonempty. So $\tilde X_k\to X$ is an isomorphism, but 
$\tilde X_{k-1}\to \tilde X_k$ is not.
So $D^{(k)}$ is locally the intersection of $k$ members
of $\Hcal$ in general position. From this it follows that 
the blowup of $\Ical^D$ factorizes over $\tilde X_{k-1}$. The pull-back
of $\sum_{H\in\Hcal} \Ocal_X(H)$ to $\tilde X_{k-1}$ is up to a twist with a 
principal ideal
equal to $\sum_{H\in\Hcal} \Ocal_{\tilde X_{k-1}}(\tilde H_{k-1})$. 
The lemma now follows with induction.
\end{proof}

A case of interest is when $X$ is the projective space $\PP (V)$
of a complex vector space $V$. If $H\in\Hcal$ is given
by the linear form $\phi_H$ on $V$, then the blowup above is simply obtained 
as follows: consider the morphism $\Omega\to \PP (\CC^\Hcal)$ defined
by $[z] \mapsto [(\phi_H(z)^{-1})_{H\in\Hcal}]$ and take 
the closure of its graph in $\PP (V)\times\PP (\CC^\Hcal)$. 

Assume now that in this situation the collection $\Hcal$ is nonempty and that 
the $H\in\Hcal$ have 
no point in common (in other words, $\Hcal$ contains a set of coordinate 
hyperplanes). Then the projection of 
$\tilde \PP(V)\to \PP (\CC^\Hcal)$ is 
birational onto its image. That image can be regarded 
as a projective completion of the hyperplane complement $\Omega$
and we therefore  denote it by $\hat\Omega$.
(In case $V=\CC^{n+1}$ and $\Hcal$ consists of the set of coordinate hyperplanes, 
then the resulting birational map $\PP^n\dasharrow \PP^n$ is the natural 
$n$-dimensional generalization of the standard Cremona transformation.)
The variety $\hat\Omega$ comes with a natural stratification 
$\{ \Omega (W)\}_W$ into smooth subvarieties. Here the index set runs over
all linear subspaces $W\subset V$ with the property  that
$\PP(W)$ is an intersection of members of $\Hcal$. To be precise: $\Omega (W)$
is the image of $\Omega$ under the projection $\Omega\to \PP (V/W)$.
So it is in fact the hyperplane complement in $\PP (V/W)$ defined by the collection
of $H\in \Hcal$ that pass through $\PP(W)$. 

The variety $\hat\Omega$ defined in the above example 
always exists as a locally compact Hausdorff space. If $X$ is projective,
then conditions can be specified under which $\hat\Omega$ will exist as 
a projective variety. Let us explain briefly how. 

The connected components of the indecomposables of the Boolean algebra
generated by the members of $\Hcal$ define a stratification of $X$.
This stratification is analytically locally trivial.
In a similar fashion, the collection of irreducible components of the preimages
of the members of $\Hcal$ determine a stratification
of $\tilde X$. The preimage of a stratum of $X$ is a union of strata 
of $\tilde X\to X$ and it is easy to see that this preimage is trivial over
the given stratum as a stratified variety. We consider now a somewhat coarser 
partition of $\tilde X$ whose members are 
indexed by the irreducible components of intersections of 
members of $\Hcal$, in which we include the 
empty set as index (this will no longer be a stratification in general:
the closure of a member need not be union of parts):
if $S$ is an irreducible 
component of some $D^{(k)}$ with $k\ge 1$, then let $P_S$ be the closure of 
the preimage of $S-(S\cap D^{(k+1)})$ in $\tilde X$ minus the points that 
lie in the closure of the preimage of $D-D^{(k)}$ and
$P_\emptyset$ will be the preimage of $\Omega$. 
So the open member of this partition 
can be identified with $\Omega$, and the closed members
of this partition are the strict transforms of the members of $\Hcal$.
For $S\not=\emptyset$, the morphism $P_S\to S$ is trivial:  $P_S$ is then
canonically a product $S\times \Omega (S)$, where $\Omega (S)$ is the
complement of a hyperplane configuration in a projective space. 
This structure defines an equivalence relation on $\tilde X$: declare
two points of $\tilde X$ to be equivalent if  they are in the same
member $P_S$ of the partition and have the same image in $\Omega (S)$
(when $S=\emptyset$, read this as: have the same image in $\Omega$).
This equivalence relation is closed and the quotient space $\hat\Omega$
is locally compact Hausdorff. 

If $X$ is projective, and we seek to put a projective structure on $\hat\Omega$, 
then the above example suggests we look for a line bundle $\Lcal$ on $X$ 
with the property that the restriction of $\Lcal$ to $H$ is isomorphic 
to the normal bundle of $H$. Its pull-back to $\tilde X$ will then 
be trivial on the equivalence classes and so 
we would like that $\sum_{H\in\Hcal}\Lcal (-H)$ is 
generated by its sections and that these sections separate the 
equivalence classes on $\tilde X$. In fact, it would be enough 
to know that $\Lcal$ restricted to $H$ is 
isomorphic to a positive power $n_H$ of the
normal bundle of $H$ and then we would ask the corresponding property for
$\sum_{H\in\Hcal}\Lcal (-n_HH)$.

\subsection{Intermediate modification of a cusp}\label{intermod}
We will look at an analogue of this situation
in the case where $X$ is a locally symmetric variety (a quotient of a
a bounded symmetric domain by an arithmetic automorphism group) and the 
hypersurfaces $H$ are totally geodesic. We then also 
wish to understand what happens if we take the closure $D^*$ of $D$ in the 
Baily-Borel compactification $X\subset X^*$ and how the blowup over $X$ extends
across that compactification. The irreducible bounded symmetric domains
admitting totally geodesic complex hypersurfaces are the domains of type
$IV$ (associated to a real orthogonal group of type $\SO (2,n)$) 
and the complex balls. Only the complex balls are relevant here, and as 
they are easier to deal with than the type $IV$ domains, 
we concentrate on them.

So let us take up the situation of Section \ref{bb}. It is known \cite{borel}
that $\G$ has a neat subgroup of finite index (this means that this
subgroup has the property that the subgroup of $\CC^\times$ 
generated by the eigen values of its elements has 
no torsion). For the purposes of this discussion, there is no loss in generality
when  passing to such a subgroup and therefore we assume that $\G$ is neat from 
the start. 

Let us now agree on a bit of notation. If $W\subset V$ is a 
degenerate positive subspace defined over $K$ with radical $I$, then $\BB$
is disjoint with  $\PP(W)$ and so the projection 
$\PP (V)-\PP (W)\to \PP(V/W)$ is defined on $\BB$. We denote the image
by $\BB(W)$ and the projection $\pi_W: \BB\to \BB (W)$. It is
easy to see that $\BB(W)=\PP(V/W)-\PP (I^\perp/W)$. 
So this is an \textit{affine space} over $I^\perp/W$.

There is an evident factorization
\begin{equation*}
\begin{CD}
\pi_W : \BB @>{\pi_I}>>\BB (I) @>>> \BB (W).
\end{CD}
\end{equation*}
The second projection is one of affine spaces.
Let us explicate $\pi_I$.

Suppose $v=(z_0,z_1,z_2,\dots ,z_n)$ are $K$-coordinates for $V$ such
that  $I^\perp$ is defined by $z_0=0$ and $\psi$ assumes the form
\[
\psi (z,w)=z_0\bar w_n+z_n\bar w_0 +\sum_{i=1}^{n-1} z_i\bar w_i.
\]
The intersection of the affine hyperplane defined by $z_0=\psi (z,e_n)=1$
with $\LL$ projects isomorphically onto $\BB$. This intersection
is given by $\Re (-z_n)> \| z'\|^2$, where $z'=(z_1,\dots ,z_{n-1})$.
In terms  of these
coordinates the projection $\pi_I$ is simply $(z',z_n)\mapsto
z'$, and hence a fibration into left half planes, indeed.
The topology near the cusp defined by $I$ is easily described in these
terms also: a neighborhood basis of this cusp intersected with $\BB$ is the
family of shifted fibrations defined by 
$\Re (-z_n)> \| z'\|^2 +a$ with $a$ a positive constant. The boundary
of such subset, in other words a fiber of the function $\Re (z_n)+\| z'\|^2$, 
is an orbit of the unipotent radical of the $\U (V)$-stabilizer of $I$.
This unipotent radical is a Heisenberg group and is described in
\ref{unipotenten}.
Since $\G$ is neat, the $\G$-stabilizer of $I$, $\G_I$, is contained in 
this Heisenberg group and is in fact a cocompact subgroup of it. So the center
of $\G_I$ is infinite cyclic and acts faithfully by purely imaginary 
translations in the fibers of $\pi_I$, whereas the quotient of $\G_I$ by 
its center acts faithfully on the affine space $\BB(I)$ as a lattice of 
maximal rank. Hence 
\[
\G_I\bs \BB\to \G_I\bs \BB (I)
\]
is a punctured disc bundle whose base is a principal homogeneous space for 
the complex torus $\G_I\bs I^\perp /I$. The associated disc bundle
can be understood as the $\G_I$-orbit space of $\BB\sqcup\BB (I)$ endowed 
with a suitable topology with the bundle 
projection given by the obvious retraction
\[
\G_I\bs (\BB\sqcup\BB (I))\to \G_I\bs \BB (I).
\]
The associated line bundle over $\G_I\bs \BB (I)$ has a Riemann form which is 
the negative of the form $\psi$ induced on the translation space $I^\perp /I$.
This implies that the dual of this line bundle is ample. So $\G_I\bs \BB (I)$
can be contracted analytically in $\G_I\bs (\BB\sqcup\BB (I))$.
The result of this contraction is that we added a singleton to $\G_I\bs\BB$.
This is the local model of the Baily-Borel compactification near the cusp 
attached to $I$  
(the added point is that cusp). The contraction mapping itself is the local model
of a well-known (orbifold) resolution of the Baily-Borel compactification, 
one that apparently has the zero section $\G_I\bs \BB(I)$ as exceptional divisor. 

Any $K$-linear subspace $W\subset I^\perp$ which contains 
$I$ defines an intermediate 
contraction and hence an intermediate modication of the cusp as follows.
The image of $W$ in $I^\perp/I$ defines a subtorus 
of $\G_I\bs I^\perp /I$. This subtorus gives rise to a torus fibration: 
\[
\G_I\bs\BB(I)\to\G_I\bs\BB(W).
\]
That fibration is the restriction of a contraction
\[
\G_I\bs (\BB\sqcup\BB (I))\to \G_I\bs (\BB\sqcup\BB (W))
\]
which leaves $\G_I\bs\BB$ unaltered. It can be performed in the 
analytic category for the same reason as for the full
contraction. So $W=I^\perp$ gives the
Baily-Borel model and $W=I$ the natural resolution. 
We still have a natural retraction
\[
\G_I\bs (\BB\sqcup \BB (W))\to \G_I\bs\BB(W)
\]
and this retraction is locally analytically trivial. 

\subsection{Compactifications of arrangement type}\label{comparr}
Now let us return to the more specific situation of \ref{symmhypersurface} 
(we continue to assume that $\G$ is neat). For every $H\in\Hcal$, $\BB (H)$ is 
totally geodesic subball of 
$\BB$ and the collection of these is locally finite on $\BB$. So 
\[
\BB (\Hcal):=\cup_{H\in\Hcal} \BB (H)
\]
is closed in $\BB$ and defines a closed analytic hypersurface $\BB (\Hcal)_\G$
of $\BB_\G$. This hypersurface is arrangementlike in the sense of 
\ref{arrmod} and hence determines a blowup $\tilde\BB\to \BB$.
This blowup is $\G$-invariant and
hence defines a blowup $\tilde\BB_\G\to \BB_\G$ of orbit spaces. 
We explain how this blowup naturally extends across the Baily-Borel 
compactification.
For every isotropic $K$-line $I\subset V$, let us 
denote by $I_\Hcal$ the intersection of $I^\perp$ and the $H\in\Hcal$
containing $I$. So $I\subset I_\Hcal\subset I^\perp$, with
$I_\Hcal= I^\perp$ in case no $H\in\Hcal$ passes through $I$.
The preceding construction attaches to the  collection 
$\{I_\Hcal\}_I$ an intermediate modification of the cusps of $\BB^*_\G$.
Let us denote this blowup
\[
\BB^\Hcal_\G\to \BB^*_\G .
\]  
Each member $H$ of $\Hcal$ passing through $I$ defines an affine 
hyperplane in $\BB(I_\Hcal)$ and hence an orbit in 
$\G_I\bs\BB(I_\Hcal)$ under a complex subtorus of codimension one. 
The closure of the image of $H$ in
$\G_I\bs(\BB\sqcup \BB (I_\Hcal))$ is the preimage of that orbit 
under the retraction of $\G_I\bs(\BB\sqcup \BB (I_\Hcal))$ onto
$\G_I\bs\BB(I_\Hcal)$. In other words the closure of the divisor  
$\BB (\Hcal)_\G$ in $\BB^\Hcal_\G$ is in an obvious sense locally trivial
near the boundary of $\BB_\G$ in $\BB^\Hcal_\G$. This implies that
the normal crossing resolution of this divisor naturally extends across
$\BB^\Hcal_\G$ to give the sought for extension of the blowup:
\[
\tilde\BB^\Hcal_\G\to \BB^\Hcal_\G\to\BB^*_\G .
\]
The closure $\BB (\Hcal )_\G^*$ of $\BB (\Hcal )_\G$ in $\BB^*_\G$
is a hypersurface and the blowup above has the virtue that the 
strict transforms of the irreducible components of this hypersurface
get separated. (This strict transform also supports an effective
Cartier divisor.)

There is a topological contraction of the exceptional locus of 
$\tilde\BB^\Hcal_\G\to \BB^*_\G$ which is of a very similar nature
as our compactification of the hyperplane complement 
$\PP (V)-D$ described in \ref{arrmod} (and is also related 
to the construction described in
\cite{looijeng}): topologically it is gotten as the $\G$-orbit space of a 
stratified extension $\hat\Omega$ of $\Omega:=\BB-\BB(\Hcal)$ as a $\G$-space. 
The strata $\Omega (W)$ of this extension are indexed by certain subspaces $W$ of $V$:
if $\Ical$ denotes the collection of $K$-hyperplanes of $V$ that are isotropic,
then $W$ is an intersection of members of $\Hcal\cup\Ical$.
We require that $W$ is \textit{not} positive definite 
or what amounts to the same, that $\PP (W)\cap\BB^*\not=\emptyset$.
The corresponding stratum $\Omega (W)$ is the image
of $\Omega$ in $\PP (V/W)$. If the algebra of $\G$-automorphic forms
on $\BB$ with arbitrary poles along $\BB (\Hcal)$ is zero in 
negative degrees and of finite type, then we believe that 
the proj of this algebra has $\hat\Omega_\G$ as underlying topological space, thus
endowing the latter with the structure of a projective variety 
that makes the contraction map $\tilde\BB^\Hcal_\G\to \hat\Omega_\G$ a morphism.

\smallskip
Almost all the compactifications we encountered in this paper appear to be 
of this type, as the following examples illustrate (proofs of these statements
are omitted).

\begin{example}\label{belangrijkvb1}
The Knudsen-Deligne-Mumford modification of $\Dcal^*=\Dcal_{12}^*$,
\[
\Scal_{12}\bs \overline{\Mcal}_{0,12}\to \Dcal_{12}^*,
\]
can via the period mapping be identified with
\[
\tilde\BB^{\Hcal}_\Uu\to \BB^*_\Uu ,
\]
where $\Hcal$ is the collection of hyperplanes perpendicular to a $3$-vector. 
\end{example}

\begin{example}\label{belangrijkvb2} 
The sequence of compactifications
\[
\Mcal^K\to \Mcal^{M*}\to \Mcal^*
\]   
is via the period mapping identified with the sequence
\[
\tilde\BB^{\Hcal_o}_{o,\Ulo}\to
\tilde\BB^{\Hcal_o(6,9)}_{o,\Ulo}\to \BB^*_{o,\Ulo},
\]
where $\Hcal_o$ is the restriction of $\Hcal$ above to the complexification of 
$\Lambda_o$ and $\Hcal_o(6,9)\subset \Hcal_o$ the subcollection of
hyperplanes of $d$-invariant $6$ or $9$. 

In this case we have a contraction of the exceptional locus that gives the Miranda
compactification of the $\Ulo$-orbit space of $\Omega:=\BB_o-\BB_o(\Hcal_o(6,9))$. 
The strata $\Omega (W)$ of the extension $\hat\Omega$ for which $W$ has hyperbolic 
signature are listed in Proposition
\ref{inbclass}. Using the obvious notation, we find the following cases:
\begin{enumerate}
\item[(0)] For a hyperplane $W$ of $d$-invariant $6$ resp.\ $9$,  $\Omega (W)$ is a 
singleton. This corresponds in $\Mcal^M$ to the single isomorphism class of a
rational elliptic surface with a $I_9$-fiber, resp.\ a $I_4^*$-fiber.
\item[(1)] For a codimension two intersection $W$ of $d$-invariant $(6,9)$ resp.\ 
$(9,9)$ we get a one-dimensional stratum $\Omega (W)$ parametrizing rational elliptic 
surfaces of type $I_8'$ (resp.\ $I''_8$).
\item[(2)] For a codimension three intersection $W$ of $d$-invariant $(6,9,9)$ we get 
a two-dimensional stratum $\Omega (W)$ parametrizing rational elliptic surfaces
with a $I_7$-fiber.
\end{enumerate}
The maximal strata come from the cases when $W$ is positive degenerate: if
we take for $W$ the intersection of all members of $\Hcal_o(6,9)$ containing an 
isotropic line of type \ $(\theta)$ resp.\ $(0)$, then $\Omega (W)$ is of 
dimension $3$ resp.\ $1$ and parametrizes rational elliptic surfaces with
an $I_6$-fiber resp.\ $I_0^*$-fiber. 
\end{example}

\appendix
\section{Unitary lattices over the Eisenstein ring}
\label{latticelambda}

In this appendix we collect and prove some properties concerning the
lattice $\Lambda$. We advise the reader first to browse through the text 
and then to consult it when the need arises.

The lattice $\Lambda$ is among the lattices considered by Allcock
in \cite{allcock}. Let us begin with an observation implicit in his
paper. Suppose $L$ is a $\ZZ$-lattice equipped with an even symmetric
bilinear form $(\; \cdot \; ):L\times L\to\ZZ$ and an orthogonal
automorphism
$\tau$ of order $6$ that has only {\it primitive} $6$th roots of
unity as
eigen values (in other words, $\tau$ satisfies $\tau^2-\tau +1=0$).
Then $L$ becomes in an obvious manner a torsion free $\Ocal$-module. 
Since $\Ocal$ is a principal ideal domain, this module will be free
also. We shall call the order $3$ automorphism $-\tau$ a {\it
triality} of
$L$ (for this notion naturally extends Cartan's use of that
term---see
below). A skew-hermitian $\Ocal$-valued form $\phi$ on $L$ is then
defined by
\[
\phi (x,y):=\omega (x\cdot y) - (x\cdot \tau y). 
\]
Using
\begin{align*}
2(\tau x\cdot x) &= -((\tau x -1)\cdot (\tau x -1) x)+
(\tau x\cdot \tau x)+(x\cdot x)\\
&=-(\tau^2 x\cdot \tau^2 x)+(\tau x\cdot \tau x)+(x\cdot x)=(x\cdot
x),
\end{align*}
we see that $\phi (x,x)=\frac{1}{2}\theta (x\cdot x)$. So for the
associated Hermitian form $\psi :=-\theta\phi$ on $L$ we have 
$\psi (x,x)=\frac{3}{2}(x\cdot x)$. In other words, 
$(\; \cdot \;)=\frac{1}{3}(\psi +\bar\psi)$.

A remarkable fact is that orthogonal reflections in $L$ (relative 
to $(\; \cdot\: )$)  
determine certain unitary reflections relative to $\psi$:
recall or note that any vector $r\in L$ with $(r\cdot r)=2$ (a
`root') defines an orthogonal $\ZZ$-linear reflection in $L$ that
sends $r$ to $-r$; likewise, the $\Ocal$-linear transformation
\[
s_r(x):= x-\omega^{-1}\phi (x,r)r.
\]
is the identity on the $\psi$-orthogonal complement of $r$
and since $\phi (r,r)=\theta$, it is immediate that $s_r$ multiplies
$r$ by the third root of unity $-\omega$. So $s_r$ is a unitary
reflection in $L$ of order $3$,  which is why such a transformation
is called a {\it triflection}. 
Note that the triflections generate a normal subgroup $\Gop (L)$ of
the unitary group $\U (L)$ of $L$.  

Conversely, every finitely generated torsion free $\Ocal$-module $L$
equipped  with a $\theta\Ocal$-valued Hermitian form $\psi$ 
(or equivalently, a $\Ocal$-valued skew-hermitian form $\phi$) 
so arises, reason for us to call 
such data an {\it $\Ocal$-lattice}. The associated (anti-linear) map 
$x\in L\mapsto \phi (-,x)\in\Hom _\Ocal (L,\Ocal)$ is
bijective precisely when the underlying even symmetric bilinear form 
$(\;\cdot\; )$ is unimodular.

Let us call $x\in L$ an $n$-{\it vector} if $\psi (x,x)=n$
(so then $3$ divides $n$).
If a positive definite $\Ocal$-lattice $L$ is spanned by its
$3$-vectors, then
the underlying even integral lattice decomposes 
canonically into an orthogonal sum of root lattices of type
$A_k$, $D_k$ or $E_k$. This decomposition is unique and hence
respected by $\tau$. Since $\tau$ cannot interchange summands
(otherwise 
it would have eigen values of order $2$ or $3$), this decomposition
is in fact one of $\Ocal$-lattices. So the indecomposable cases must
be of type $A_\even$, $D_{\even\ge 4}$, $E_6$ and $E_8$.
On the other hand, it is easy to see that a triality cannot exist
inside  the Weyl groups $W(A_k)$ or $W(B_k)$ for $k$ even and at
least $4$.
So the possible indecomposable $\ZZ$-lattices with a triality are of
type $A_2$, $D_4$, $E_6$ and $E_8$. For example, a type $D_4$ root
lattice
admits a triality in $W(F_4)$ (which is in fact the automorphism
group of the underlying $\ZZ$-lattice). By inspecting Carter's
description
of conjugacy classes in exceptional
Weyl groups \cite{carter} we find that for a root lattice of type
$A_2$, 
$D_4$, $E_6$ and $E_8$ a triality exists and is unique up to
conjugacy. They can be gotten in a uniform manner as follows: let
$\Lambda^k$
be the $\Ocal$-lattice with basis $r_1,\dots,r_k$, such that each 
$r_i$ is a $3$-vector, $\psi (r_i,r_{i+1})=\theta$ for $i=1,\dots
k-1$ 
and $\psi (r_i,r_j)=0$ when $j>i+1$. So $\Lambda^{10}$ is the
$\Ocal$-lattice 
encountered in Section \ref{cycliccov}. One may verify that
$\Lambda^k$ is positive definite iff $k=1,2,3,4$ and that in these
cases  the underlying root lattice are of type $A_2$, $D_4$, $E_6$,
$E_8$
respectively. (For $k=2$, we get the classical triality on $D_4$.) 
By means of Coxeter \cite{coxeter} we identify $G(\Lambda^k)$ in
Shephard and Todd's 
Table VII in \cite{shep}. The associated triflection group $\Gop
(\Lambda^k)$ 
appears there with number $4)$ for $k=2$, $25)$ for $k=3$ and $32)$
for $k=4$. (The group 
$\Gop (\Lambda^3)$ is the Hesse group of symmetries of the Hesse 
pencil $\lambda(x^3+y^3+z^3)
+\mu(xyz)$; $\Gop (\Lambda^4)$ is sometimes called the {\it Witting
group}.)

\subsection{The lattice $\Lambda^4$}\label{lambda4}
The case $E_8$ is of particular interest: following
\cite{allcock}, $\tau$ is then realizable as the $5$th power of a 
Coxeter transformation. (A Coxeter transformation of such a
root lattice has order $30$ and its eigen values are the eight 
primitive $30$th roots of unity.) 

The $3$-vectors of $\Lambda^4$ are the
roots of the $E_8$-lattice, hence there are $240$ of them. If we 
identify $\Ocal /\theta\Ocal$ with $\FF_3$, then
\[
\Lambda^4_{\FF_3}:=\FF_3\otimes_\Ocal\Lambda^4\cong
\Lambda^4/\theta\Lambda^4
\]
gets the structure of a vector space of dimension $4$ over $\FF_3$.
The skew hermitian $\Ocal$-valued form 
$\phi$ on $\Lambda_4$ induces a symplectic $\FF_3$-valued form 
on $\Lambda^4_{\FF_3}$. It turns out to be nondegenerate. There
results a homomorphism 
\[
\U (\Lambda^4)\to\Symp(\Lambda^4_{\FF_3})\cong \Symp (4,\FF_3)
\]
which Allcock shows to be surjective with kernel the scalar subgroup
$\mu_3$.
(Note that $\omega+1$ is divisible by $\theta$, so that $\omega$ acts
as minus the identity in $\Lambda^4_{\FF_3}$.)  
In particular, $U(\Lambda^4)$ is transitive on 
$\Lambda^4_{\FF_3}-\{ 0\}$. 
He further observes that every nonzero element of $\Lambda^4_{\FF_3}$
has in its preimage precisely three $3$-vectors (a $\mu_3$-orbit). 
Allcock uses this to prove:

\begin{lemma}[\cite{allcock}, Theorem 5.2]\label{simpeltrans}
The group $\U (\Lambda^4)$ acts transitively on the set of
$6$-vectors and
on the set of $3$-vectors in $\Lambda^4$.
\end{lemma}

We shall further exploit this reduction to study $3$- and $6$-vectors
in the $\Lambda^4$-lattice. We begin with noting that 
it remembers the relative position of the $\mu_6$-orbits of two
nonproportional $3$-vectors $r,r'$: the fact that these two span a
positive definite lattice of rank two implies that $|\psi (r,r')|<3$,
and as $\psi (r,r')$ is divisible by $\theta$, we have either $\psi
(r,r')=0$
or $\psi (r,r')\in\mu_6\theta$. This means that their images
in $\Lambda^4_{\FF_3}$ span an isotropic resp.\ nondegenerate rank
two sublattice.

The description of the $6$-vectors in terms of this reduction must be
less straightforward, witness the fact that there are $80.27$
$6$-vectors and $80$ nonzero elements in $\Lambda^4_{\FF_3}$. The
next lemma
offers one such description.

\begin{lemma}\label{zesvektor}
A $6$-vector $z\in\Lambda^4$ can be written in exactly three ways as
the sum of two $3$-vectors $z=r_1+r_2$ with $\psi (r_1,r_2)=\theta$. 
All such  pairs $r_1,r_2$ span the same rank two sublattice 
$L_z$ of $\Lambda^4$. The image of $L_z$ in $\Lambda^4_{\FF_3}$ is a 
nondegenerate plane and assigning to $z$ the mod $\theta$ reduction
of the pair  $(z,L_z)$ defines a bijection between the set of
$\mu_3$-orbits
of $6$-vectors in $\Lambda^4$ and the set of pairs $(v,P)$, where 
$P\subset\Lambda^4_{\FF_3}$ is a nondegenerate plane and $v\in P-\{
0\}$.
\end{lemma}
\begin{proof} 
Consider the set $S$ of pairs of $3$-vectors $(r,r')$ in
$\Lambda^4$ with $\psi (r,r')=\theta$. The mod $\theta$ reduction of
a
pair $(r,r')\in S$
is pair of vectors $(v,v')$ in $\Lambda^4_{\FF_3}$ with symplectic
product $1$. The number of such pairs of vectors is $80.27$. The
$3$-vectors
mapping to $v$ are the elements of the $\mu_3$-orbit of $r$ and
likewise for $r'$. So the preimage of $(v,v')$ in $S$ is the
$\mu_3$-orbit of the {\it pair} $(r, r')$. Hence
$S$ has $80.27.3$ elements. The image of the map $(r,r')\in S\mapsto
r+r'\in\Lambda^4$  consists of $6$-vectors, hence is the set of
all $6$-vectors, since it is $U(\Lambda^4)$-invariant. As there are
$80.27$ $6$-vectors, we see that each $6$-vector occurs precisely
three times.  If $(r,r')\in S$, then $(\omega r', 
r+(1-\omega)r')$ and $((1-\omega)r , \omega r+r')$ are two other
elements of $S$ with the same sum. So there are no more elements in
$S$ with that property. Hence the span  of $r$ and $r'$ only depends
on
$r+r'$. All the assertions of the lemma now have been proved.
\end{proof}

Allcock's result says that mod $\theta$ reduction gives a bijective
correspondence between the $\mu_6$-orbits of
$3$-vectors and the lines $\ell\subset\Lambda^4_{\FF_3}$.  
Lemma \ref{zesvektor} can be understood as
asserting a  similar relationship between the $\mu_6$-orbits  of 
$6$-vectors and the flags $(\ell ,P)$ in
$\Lambda^4_{\FF_3}$, where $\ell$ is a line in a nondegenerate plane
$P$. Since symplectic geometry over a finite field is a priori a lot
simpler than unitary geometry over the
Eisenstein ring, such an interpretation is helpful when determining
the relative position of a $3$-vector and $6$-vector in $\Lambda^4$. 
To see this, note that for a nondegenerate flag $(\ell ,P)$ 
in $\Lambda^4_{\FF_3}$and a line 
$\ell'$ in $\Lambda^4_{\FF_3}$ the following
possibilities present themselves: 
\begin{enumerate}
\item[(a)] $\ell'=\ell$,
\item[(b)] $\ell +\ell' =P$,
\item[(c)] $\ell'\not\subset P$, $\ell'$ not perpendicular to $\ell$,
\item[(d)] $\ell'\not\subset P$, $\ell'$ perpendicular to $\ell$ but
not
to $P$,
\item[(e)] $\ell'$ perpendicular to $P$.
\end{enumerate}
By elementary symplectic geometry, each of these cases represents a
single orbit under the symplectic group. Let us see what this tells
us 
about the relative position of a $3$-vector $r$ and a $6$-vector $z$.
From the preceding it follows that the unitary group of $\Lambda^4$
has precisely five orbits in the set of pairs of $\mu_6$-orbits
$(\mu_6 .r,\mu_6 .z)$.
We give in each of the five cases above a representative example
with $z=r_1+r_2$ (so that $L_z=\Ocal r_1+\Ocal r_2$ and hence $P$ is
the image of $L_z$ in $\Lambda^4_{\FF_3}$). 
\begin{enumerate}
\item[(a)] $r=\omega^2r_1+r_2$ (so $\psi (r,z)=0$), 
\item[(b)] $r=r_1$ (so $\psi (r,z)=3+\theta$), 
\item[(c)] $r=r_3$ (so $\psi (r,z)=-\theta$),
\item[(d)] $r=\omega r_2+r_3$ (so $\psi (r,z)=3$),
\item[(e)] $r=r_4$ (so $\psi (r,z)=0$).
\end{enumerate}
The case (a) is somewhat special: then $r$ and $z$ are perpendicular
and span an imprimitive sublattice. We also see that the orthogonal
complement of $z$ in $L_z$ is spanned by $r$. So any 
$3$-vector with the same mod $\theta$-reduction as $z$ lies in $L_z$
and spans with $z$ a subgroup of finite index in $L_z$.

\begin{corollary}\label{miniclass}
Let $z=r_1+r_2$ be the standard $6$-vector in $\Lambda^4$. Then the
set of those $3$-vectors
in $\Lambda^4$ which have a fixed nonzero Hermitian inner product
with $z$ make up a single $\U (\Lambda^4)_{z}$-orbit.
The $3$-vectors perpendicular to $z$ span a lattice of type
$\Lambda^1\times\Lambda^2$, with basis $(\omega^2r_1+r_2, r_1-\theta
r_2 -2r_3+\theta r_4,r_4)$. A $3$-vector perpendicular to
$z$ spans with $z$ a primitive sublattice if and only if it belongs
to the $\Lambda^2$-summand (hence any two such are in the same
$\U (\Lambda^4)_{z}$-orbit). 
\end{corollary}
\begin{proof} 
Let $r'$ be a $3$-vector in $\Lambda^4$ 
with $\psi (r',z)\not= 0$. It follows from the preceding that
$r'$ is $\U (\Lambda^4)_z$-equivalent to $\omega^ir$, with $i\in \ZZ
/6$ and $r$ a vector mentioned in one of the cases
$\text{(b)},\text{(c)} ,\text{(d})$.
In these cases the exponent $i\in \ZZ /6$ is
determined by the inner product of $r'$ with $z$. 
The last part of the corollary is
straightforward. 
\end{proof}

In case (d), we have that $z-r$ is 
a $3$-vector perpendicular to $r$. We shall need to know 
in how many ways $z$ can be written as a sum of two perpendicular
$3$-vectors. 

\begin{corollary}\label{zesinlambda4}
A $6$-vector $z$ is written in exactly $4$ distinct ways as 
a sum of two perpendicular $3$-vectors in $\Lambda_4$. These vectors 
are orthogonal to the orthogonal complement of $z$ in $L_z$, and so
span with $z$ a rank $3$ sublattice of $\Lambda^4$. 
(For $z=r_1+r_2$ these sum decompositions are $z=(z-r)+r$ with
$r=\omega r_2+r_3$,  $\omega r_2+r_3+\omega^{-1}r_4$,
$\omega r_2+r_3+\omega^{-2}r_4$, $\omega r_1+2r_2-\theta r_3 -r_4$.)
\end{corollary}
\begin{proof}
We begin with noting that a line $\ell'$ in
$\Lambda^4_{\FF_3}$ has property (d) if and only if it is the graph
$\ell_f$ of a nonzero homomorphism $f:\ell\to P^\perp$. It is clear
that there are $8$ such lines. They come in $4$ pairs: we have 
$\ell\subset
\ell_f+\ell_{f'}$ if and only if $f+f'=0$. In that case 
$\ell_f+\ell_{f'}$ is isotropic and
so $\ell_f,\ell_{f'}$ correspond to $\mu_6$-orbits of $3$-vectors 
that are perpendicular. There are unique $3$-vectors  $r_f,r_{f'}$ in
these 
orbits with sum $z$. Let $r$ be a $3$-vector with the same mod 
$\theta$-reduction as $z$. Since $\ell -\ell_f$ is isotropic,  
$r$ is also orthogonal to $r_f$. We noted that $r$ spans the
orthogonal complement of $z$ in $L_z$, and so the second assertion 
of the corollary follows.
\end{proof}

The following is proved in a similar fashion as \ref{miniclass}. The
proof
is in fact easier and so we omit it.

\begin{lemma}
Let $r\in \Lambda^4$ be a $3$-vector. Then the stabilizer group
$\U (\Lambda^4)_r$ acts transitively on the set of those $3$-vectors
in $\Lambda^4$ which have a fixed inner product with $r$.
\end{lemma}

\subsection{The lattice $\Lambda$}
A {\it hyperbolic $\Ocal$-lattice} is
obtained as follows: let $M$ be free finitely generated
$\Ocal$-module. 
Regard $M$ as $\ZZ$-module and consider 
$\Hom_\ZZ (M,\ZZ)\oplus M$. This has the  natural quadratic form  
$q(\xi ,x)=\xi (x)$ for which it is an unimodular $\ZZ$-lattice. Now
let
$\Ocal$ act on $M$ as before and on $\Hom_\ZZ(M,\ZZ)$
contragradiently. 
Then the preceding construction turns $\Hom_\ZZ(M,\ZZ)\oplus M$ into
a
nonsingular $\Ocal$-lattice $H_M$. So $H_\Ocal\cong\Ocal^2$ with
Hermitian form $\psi (z,w)= \theta (z_1\bar w_2 - z_2\bar w_1)$.
Notice that the skew-hermitian form $\phi=-\theta^{-1}\psi$ has
discriminant $1$. We shall denote the given basis of $H_\Ocal$ by
$(e,f)$.
Consider the lattice $\Lambda^4\perp \Lambda^4\perp H_\Ocal$ and 
denote the first two summands
$\Lambda'$ and $\Lambda''$ (with basis $(r_i')_{i=1}^4$ resp.\
$(r_i'')_{i=1}^4$). We shall identify $\Lambda=\Lambda^{10}$ with 
$\Lambda'\perp\Lambda''\perp H_\Ocal$ by 
means of the unitary isomorphism
\[
(r_1,\dots ,r_{10})\mapsto (r_1'',\dots ,r_4'',s+e, \omega e+\theta
f,
\omega^{-1}e+r'_1,r'_2,r'_3,r'_4),
\]
where $s\in\Lambda''$ is characterized by the fact that it is
perpendicular to
$r_1'',r_2'' ,r_3''$ and $\psi (s,r_4'')=\theta$. This shows in
particular that $\Lambda$ has signature $(9,1)$, as asserted earlier.
Notice that this isomorphism also identifies $\Lambda^6$ 
(the span of $r_1,\dots ,r_6$) with $\Lambda''\perp H_\Ocal$.  
 
Allcock \cite{allcock} proves that $\U (\Lambda )$ acts transitivily
on the set primitive $0$-vectors of $\Lambda$, in other words, every
primitive
$0$-vector of $\Lambda^4\perp \Lambda^4\perp H_\Ocal$
can be transformed by a unitary transformation into $e$. 
We derive from this the corresponding statement for the set of
$6$-vectors:

\begin{proposition}\label{trans6}
Each $3$-vector in $\Lambda$ is perpendicular to a primitive null
vector and the group $\U (\Lambda )$ acts transitively on the set of
pairs 
$(r,n)$ with $r$ a $3$-vector and $n$ a primitive $0$-vector
perpendicular to $r$ (in particular, $\U (\Lambda )$ is transitive on
the set of $3$-vectors in $\Lambda$).

Similarly, each $6$-vector in $\Lambda$ is perpendicular to a
primitive null vector. The group $\U (\Lambda )$ acts transitively on
the set of $6$-vectors, but has two orbits in the set of pairs 
$(z,n)$ with $z$ a $6$-vector
and $n$ a primitive $0$-vector perpendicular to $z$. These two orbits
are represented by $(r_1'+r_2',e)$ (type $(\theta)$) and
$(r_1'+r_1'',e)$ (type $(0)$). 
\end{proposition}

\begin{none}\label{unipotenten}
Before we begin the proof it is useful to make a few general observations.
Let $V$ be a finite dimensional complex vector space
equipped with a skew-hermitian form $\phi$. 
Let also be given a nonzero isotropic vector $e\in V$. For every
$v\in V$ with $\phi (v,e)=0$ we define the transformation $T_{e,v}$
in $V$ by 
\[
T_{e,v} (x)= x+ \phi(x,e)v+\phi(x,v)e+\tfrac{1}{2}\psi (v,v)\psi
(x,e)e
\]
One checks that $T_{e,v}$ is unitary and fixes $e$. Its action in
$e^\perp$ is simply given by $x\in e^\perp\mapsto x+\phi(x,v)e$.
Notice that $T_{e,v}$ only depends on the image of $v$ in $e^\perp
/\RR
\sqrt{-1}e$. We have
\[
T_{e,u}T_{e,v}=T_{e,u+v+\frac{1}{2}\phi (v,u)e}.
\]
These transformations make up the unipotent radical of the stabilizer
of $e$ in the unitary group  $\U (V)$. It is a Heisenberg group with
center the transformations $T_{e,\lambda e}$ with $\lambda$ real.
Suppose that $L\subset V$ is a discrete $\Ocal$-submodule in $V$ 
of maximal rank such that $\phi$ takes on $L\times L$ 
values in $\Ocal$. If $e$ and $v$ lie in $L$ and 
$\phi (v,v)$ is even, then clearly $T_{e,v}$ preserves $L$. 
So if $x\in e^\perp\cap L$, then $x+\Ocal e$
is contained in a $U(L)_e$-orbit if $\phi (x,v)=1$ for some
$v\in L\cap e^\perp$ with $\phi (v,v)$ even. Or what amounts to the
same, if $v\in L\cap e^\perp$ with $\psi (x,v)=-\theta $ and $\psi
(v,v)\in 6\ZZ$.
\end{none}

\begin{proof}[Proof of \ref{trans6}]
We only prove the statements involving a $6$-vector, the proof of the
one about a $3$-vector is similar and easier.
We begin with the last clause. Let $(z,n)$ be as in the proposition. 
By Allcock's result, a unitary transformation will map this into a
pair
with second component $e$ and so we may assume that $n=e$. Then
$z$ can be written $x'+x''+\lambda e$ with $x'\in\Lambda'$,
$x''\in\Lambda''$ and $\lambda\in \Ocal$. We must have $\psi
(x',x')+\psi
(x'',x'')=6$. Since the two terms must be nonnegative multiples of
three 
they are $(6,0)$, $(3,3)$ or $(0,6)$. The stabilizer of $e$
contains the interchange of $\Lambda'$ and $\Lambda''$ as well
as the unitary group of each of these summands. So we can
eliminate the last case and by \ref{simpeltrans} assume that
$(x',x'')=(r_1'+r_2',0)$ or $(x',x'')=(r_1',r_1'')$. In either case, 
there exists a $6$-vector $v\in\Lambda'$ with $\psi (z,v)=\theta$ 
and so by the discussion (\ref{unipotenten}) there exists a
unitary transformation fixing $e$ that sends $z$ to $x'+x''$. 
The last assertion follows.

We next show that any $6$-vector $z$ is perpendicular to a
primitive null vector. The orthogonal complement $\Lambda_z$ of $z$
is a free  $\Ocal$-module of signature $(8,1)$. So its
complexification
$\CC\otimes_\Ocal\Lambda_z=\RR\otimes_\ZZ\Lambda_z $ represents zero.

Its real dimension is $\ge 5$ and 
a theorem of Meyer \cite{serre} then implies that $\Lambda_z$ also
represents zero. In other words, there exists a null vector
perpendicular to $z$. 

It remains to see that $r_1'+r_2'$ and $r_1'+r_1''$ are in the same
$U(\Lambda)$-orbit. This is left to the reader.
\end{proof}

\subsection{The lattice $\Lambda_o$}
Let us now fix a sublattice $\Lambda_o\subset\Lambda$ 
that is the orthogonal complement of a $6$-vector $z_o\in\Lambda$. 
In view of \ref{trans6} all such sublattices are unitary equivalent. 

\begin{proposition}\label{unitairiso}
The $\Uu$-stabilizer of $\lo$ maps isomorphically to
the unitary group $\Ulo$ of $\lo$.
\end{proposition}

The proof is a modification of a standard argument in lattice
theory. In order to make it transparent we begin with a
general discussion. Given an $\Ocal$-lattice $L$, let us simply write
$L^*$
for $\Hom_\Ocal (L,\Ocal)$. The skew-hermitian form
$\phi_L:=-\theta^{-1}\psi_L$
on $L$ induces an antilinear map $a_L: L\to L^*$, $x\mapsto \phi (\;
,x)$.
Suppose that $\phi_L$ is nondegenerate (i.e., has nonzero
discriminant). Then $a_L$ maps $L$ bijectively onto a sublattice of
$L^*$ of finite
index, so that $C(L):=L^*/a_L(L)$ is a finite $\Ocal$-module. The
order of
$C(L)$ is then the square absolute value of the discriminant of $L$.
For
instance, if $L$ is spanned by a $3n$-vector, then $C(L)\cong
\Ocal/(n\theta)$,
which has indeed order $3n^2$.
The form $\phi_L$  determines a skew-hermitian form $\phi_{L^*}$ on
$L^*$
such that $\psi_{L^*}(a_L(x),a_L(y))=\psi_L(y,x)$. This form now
takes values in the field $\QQ (\omega)$. If however one of its
arguments
lies in the image of $a_L$, then it takes values in $\Ocal$. So
$\psi_{L^*}$ induces a skew-hermitian form $\phi_{C(L)}: C(L)\times
C(L)\to
\QQ (\omega)/\Ocal$. It is clear that every unitary
transformation of $L$ induces a unitary transformation in $C(L)$.

Suppose now $L$ of discriminant $\pm 1$ and let $M\subset L$ be a
primitive
nondegenerate submodule with orthogonal complement $N$. 
So $M\perp N$ sits in $L$ as a submodule of finite index.
Composing  $a_L$ with restriction to $M\perp N$ gives an embedding
of $L/(M+N)$ in $C(M\perp N)=C(M)\perp C(N)$. This image is isotropic
for the skew-hermitian $\QQ (\omega)/\Ocal$-valued form on $C(M)\perp
C(N)$.
Since $L$ has discriminant $\pm 1$, it is a maximal sublattice in
$\QQ\otimes_\ZZ L$ on which $\phi$ is $\Ocal$-valued, an so its image
in $C(M)\perp C(N)$ is maximally isotropic. It is clear that the
projection of this image in either summand is a bijcetion. In other
words, the
image is the graph of an isomorphism $\alpha : C(M)\cong C(N)$ which
changes the sign of the forms.

It is clear that an automorphism of $M\perp N$ preserves $L$ if and
only if
is preserves the image of $L$ in $C(M)\perp C(N)$.
So a pair of unitary transformations of $M\perp N$ of the form
$(u_M,u_N)$  preserves $L$ if and only if $\alpha$ commutes with the
unitary transformations in $C(M)$ and $C(N)$ induced by $u_M$ and
$u_N$.

\begin{proof}[Proof of \ref{unitairiso}]
We apply this to the case at hand:
$L=\Lambda$, $M=\lo$ and $N$ spanned by the $6$-vector $z_o$. 
Then $C(\lo)\cong C(\Ocal z_o)\cong
\Ocal/(2\theta)$, where in the latter case the skew form takes the
value $\frac{1}{2}\theta^{-1}$ on a generator. One easily verifies
that the
group of unitary transformations of  $\Ocal/(2\theta)$ is $\mu_6$. As
this
is also the group of unitary transformations of  $\Ocal z_o$, it
follows that
every  unitary transformation of $\lo$ extends uniquely to unitary
transformation of $\Lambda$.
\end{proof}

In order to classify the $3$-vectors in $\Lambda$ relative to 
$z_o$, we first consider the abstract
$\Ocal$-lattices spanned by a $3$-vector and a $6$-vector.

\begin{lemma}\label{class}
Let $L$ be a positive definite $\Ocal$-lattice of rank two spanned by
a $6$-vector $z$ and $3$-vectors. Then we are in one of the following
four cases: $L$ has a basis $(e_1,e_2)$ such that
\begin{itemize}
\item[($\delta_6$)] $z=e_1+e_2$ and $\psi$ has the matrix
\begin{math}
\left(
\begin{smallmatrix}
3 & \theta\\
\bar\theta& 3\\
\end{smallmatrix}
\right)
\end{math} 
so that $L$ has discriminant $6$ or
\item[($\delta_9$)] $z=e_1+e_2$ and $\psi$ has the matrix
\begin{math}
\left(
\begin{smallmatrix}
3 & 0\\
0& 3\\
\end{smallmatrix}
\right)
\end{math}
so that $L$ has discriminant $9$ or
\item[($\delta_{15}$)] $z=e_1$ and $\psi$ has the matrix
\begin{math}
\left(
\begin{smallmatrix}
6 & \theta\\
\bar\theta& 3\\
\end{smallmatrix}
\right)
\end{math}
so that $L$ has discriminant $15$ or
\item[($\delta_{18}$)] $z=e_1$ and $\psi$ has the matrix
\begin{math}
\left(
\begin{smallmatrix}6 & 0\\
0 & 3\\
\end{smallmatrix}
\right)
\end{math}
so that $L$ has discriminant $18$.
\end{itemize}
Moreover, is $M\supset L$ a rank two $\Ocal$-lattice that strictly 
contains $L$, then we are in case $\delta_{18}$ and $M$ is isomorphic
to the lattice of case $\delta_6$.
\end{lemma}
\begin{proof}
Suppose first $L$ spanned by the $6$-vector $z$ and a $3$-vector $r$.
We have $\psi (z,r)=\theta u$ for some $u\in\Ocal$.  
Since $L$ is positive definite, we must have $| u|^2\le 6$. Since
$u\in\Ocal$, this implies that up to a unit $u$ equals $0$, $1$, $2$
or $\theta$. By multiplying $r$ with a unit we may assume that
$u$ acually equals one of these values. For $u=0$ we get case
$\delta_{18}$,
and for $u=1$ we get case $\delta_{15}$. For $u=2$ we get case 
$\delta_9$ by taking $(e_1,e_2)=(z,\omega^{-2}r+\omega^{-1}z)$ and
for $u=3$
we get case $\delta_6$ by taking $(e_1,e_2)=(z-r,r)$.

For the last part of the lemma, we observe that for an overlattice
$M\supset L$ we must have that the quotient of the discriminant of
$M$ by the
discriminant of $L$ must be the norm of an element of $\Ocal$.
Since the discriminant of $M$ is also divisible by $3$, this implies
that
$L$ is of type $\delta_9$ or $\delta_{18}$. The 
case $\delta_9$ has as underlying integral lattice a root lattice of
type 
$A_1\perp A_1$. This admits no even overlattice and hence cannot
occur. 
There remains the case that $L$ is of type $\delta_{18}$ with $M$ of
discriminant $6$. It is then not hard to see that $M$ is as asserted.
 
\end{proof}

If $r\in\Lambda$ is a $3$-vector, which together with 
$z_o$ spans a primitive positive definite sublattice of $\Lambda$,
then according to Lemma \ref{class}, the discriminant of 
this sublattice can take $4$ values : $6$, $9$, $15$ or $18$. We call
this value the {\it $d$-invariant} of $r$. Proposition \ref{trans6}
shows
that the primitive isotropic lines $I\subset \Lambda_o$ come in
two types (types $(0)$ and $(\theta)$) and that each type is 
represented by a single $\U (\Lambda_o)$-orbit.

\begin{proposition}\label{dclass}
Let $I\subset\Lambda_o$ be a primitive isotropic line
and denote by $I(6)$ resp.\ $I(9)$ the span of 
$I$ and the $3$-vectors $r\in I^\perp$ with 
$d$-invariant $6$ resp.\ $9$. Then:
\begin{enumerate}
\item[($\theta$)] If $I$ is type $(\theta )$, then  $I(6)/I$ and
$I(9)/I$ 
are perpendicular sublattices of $I^\perp/I$ of 
rank $1$ and $2$ respectively. Moreover, there are precisely $4$
rank one sublattices of $I^\perp/I$ spanned by the image of 
a $3$-vector in $I^\perp$ of $d$-invariant $9$.
\item[(0)] If $I$ is type $(0)$, then $I(6)=I$ and $I(9)/I$ is of
rank $1$.
\end{enumerate}
\end{proposition}
\begin{proof}
By Proposition \ref{trans6} we may assume that $I$ is spanned by $e$
and that $z_o=r_1'+r_2'$ in case $(\theta)$ and $z_o=r_1'+r_1''$ 
in case $(0)$. This identifies $I^\perp/I$ with
$\Lambda'\perp\Lambda''$.
A $3$-vector in $I^\perp$ maps to a $3$-vector in $I^\perp/I\cong
\Lambda'\perp\Lambda''$ of the same $d$-invariant and the $3$-vector
of $\Lambda'\perp\Lambda''$ lies in $\Lambda'$ or in $\Lambda''$.

In case $(\theta)$ it is clear that any $3$-vector in $\Lambda''$ has
$d$-invariant $18$, so if we are after the $3$-vectors of
$d$-invariant $6$ or $9$, then we only have to deal with 
$\Lambda'$. Assertion ($\theta$) then follows from 
Lemma \ref{zesvektor} and Corollary \ref{zesinlambda4}.

Case $(0)$ follows from the simple observation that $r_1'+r_1''$
cannot be written in any other way as a sum of two $3$-vectors in
$\Lambda'\perp\Lambda''$.
\end{proof}

\begin{proposition}\label{inbclass} 
Let $L\subset \Lambda$ be a primitive sublattice containing
$z_o$. Then $L^\perp$ is isomorphic to an orthogonal product of
lattices
$\Lambda^{k_1}\perp\Lambda^{k_2}\perp\cdots $ with $k_1\ge 6$ 
if and only if $L$ is spanned by $z_o$ and
$3$-vectors of $d$-invariant $6$ and $9$ and we are then
in one of the following cases: 
\begin{enumerate}
\item[$(6)$]  $(L,z_o)\cong(\Lambda^2,r_1+r_2)$,
$L^\perp\cong \Lambda^8$ and $L$ is spanned by $z_o$ and a
$3$-vector of $d$-invariant $6$, 
\item[$(9)$] $(L,z_o)\cong (\Lambda^1\perp\Lambda^1,r+r')$,
$L^\perp\cong \Lambda^7\perp\Lambda^1$ and is spanned by $z_o$ and a
$3$-vector of $d$-invariant $9$,
\item[$(6,9)$] $(L,z_o)\cong (\Lambda^3,r_1+r_2)$, $L^\perp\cong
\Lambda^7$  and $L$ is spanned by $z_o$ and two $3$-vectors of
$d$-invariant $6,9$. 
\item[$(9,9)$] $(L,z_o)\cong (\Lambda^3,r_1+r_3)$,
$L^\perp\cong \Lambda^7$ and $L$ is spanned by $z_o$ and two
$3$-vectors of $d$-invariant $9,9$. 
\item[$(6,9,9)$] $(L,z_o)\cong (\Lambda^4,r_1+r_2)$,
$L^\perp\cong \Lambda^6$ and $L$ is $L$ is spanned by $z_o$ and
three $3$-vectors of $d$-invariant $6,9,9$.
\end{enumerate}
Each of these possibilities respresents a single 
$\U(\Lambda_o)$-equivalence class and this is also the complete
list of $\U(\Lambda_o)$-equivalence classes of
positive definite sublattices of $\Lambda$ spanned
by $z_o$ and $3$-vectors of $d$-invariant $6$ and $9$.
\end{proposition}
\begin{proof} 
Let us first assume that $L^\perp$ is isomorphic to an orthogonal
product 
$\Lambda^{k_1}\perp\Lambda^{k_2}\perp\cdots $ with $k_1\ge 6$.
Since $\Lambda\cong \Lambda^6\perp\Lambda^4$ and
$\Lambda^{k_1}\cong\Lambda^6\perp\Lambda^{k_1-6}$, 
we see that it is enough to investigate
the corresponding issue in $\Lambda^4$. The $6$-vectors in
$\Lambda^4$ are all unitary equivalent, and so we can assume that
$z=r_1+r_2$. The assertions regarding the classification now follow 
from Corollary \ref{miniclass}. 

Assume now that $L\subset\Lambda$ is a positive definite sublattice
and spanned by $z_o$ and $3$-vectors of $d$-invariant $6$ and $9$.
Assume also that its rank is $\le 5$. Then the 
orthogonal complement of the lattice $L$ is hyperbolic of
sufficiently high rank and so by Meyer's theorem contains a 
primitive null vector. We may assume that this null vector is $e$ and
that $z_o$ is either $r_1'+r_2'$ or $r_1'+r_1''$. So $L$ projects
isomorphically  to a sublattice $\bar L\subset \Lambda'\perp
\Lambda''$ spanned by
$3$-vectors. Since the $3$-vectors helping to span $L$ are of 
$d$-invariant $6$ or $9$, \ref{dclass} implies that
$\bar L\subset \Lambda'$ when $z_o=r_1'+r_2'$ and 
$\bar L\subset \Ocal r_1'+\Ocal r_1'' $ when $z_o=r_1'+r_1''$.
In particular, $L$ is of rank $\le 4$. All the assertions 
now follow in a straightforward manner from \ref{miniclass},
\ref{dclass} and  \ref{unipotenten}.
\end{proof}

\begin{proposition}\label{dchar}
The $3$-vectors in $\Lambda$ of fixed $d$-invariant form a single
$U(\Lambda_o)$-orbit.
\end{proposition}
\begin{proof}
For $d= 6$ or $9$ this is part of the statement of the previous
proposition.
The cases $d=15$ and $d=18$ are handled in a similar way.
\end{proof}


\begin{thebibliography}{99}

\bibitem{abramvist}
D.~Abramovich, A.~Vistoli: \textit{Complete moduli for 
fibered surfaces}, in: Recent Progress in Intersection Theory,  
Proc. Intern.\ Conf.\ on Intersection 
Theory, Bologna 1997, G. Ellingsrud, W. Fulton, A. Vistoli (eds.),
Birkh\"auser 2000, see also \texttt{arXiv.org/abs/math.AG/9804097},\par
\texttt{~/math.AG/9811059}, \texttt{~/math.AG/9908167}. 
	
\bibitem{allcock} 
D.~Allcock:
{\it The Leech Lattice and Complex Hyperbolic Reflections},
Invent.\ Math.\ 140 (2000), 283--301. 

\bibitem{all} 
D.~Allcock, J.A.~Carlson, D.~Toledo:
{\it The complex hyperbolic geometry for moduli of cubic surfaces}, 
55 pp., \texttt{arXiv.org/abs/math/0007048}, see also   
{\it A complex hyperbolic structure for moduli of cubic surfaces},
C.~R.\ Acad.\ Sci.\ Paris, t.\ 326, S\'erie I (1998), 49--54. 

\bibitem{bpv} 
W.~Barth, C.~Peters, A.~Van de Ven:
{\it Compact Complex Surfaces},
Erg.\ der Math.\ u.\ i.\ Grenzgebiete, 3e Folge, 4,
Springer-Verlag, Berlin etc.\ 1984.

\bibitem{borel} 
A.~Borel:
\textit{Introduction aux groupes Arithm\'etiques}
Hermann, Paris 1969.


\bibitem{bourbaki}
N.~Bourbaki: 
{\it Groupes et alg\`ebres de Lie, Ch.\ 4,5,6},
Masson, Paris 1981.

\bibitem{carter}
R.W.~Carter: 
{\it Conjugacy classes in Weyl groups,} 
in: Seminar on algebraic groups
and related finite groups, Lect. Notes in Math. 131, 
Springer-Verlag, Berlin etc.\ 1970.

\bibitem{couw}
W.~Couwenberg:
{\it Complex Reflection Groups and Hypergeometric Functions},
Thesis (123 p.), Katholieke Universiteit Nijmegen, 1994.

\bibitem{coxeter}
H.S.M.~Coxeter:
{\it Finite groups generated by unitary reflections},
Abh.\ Math.\ Sem.\ Hamburg 31 (1967), 125--135.


\bibitem{delmost1}
P.~Deligne, G.D.~Mostow:
{\it Monodromy of hypergeometric functions and non-lattice integral
monodromy},
Publ.\ Math.\ IHES 63 (1986), 58--89.

\bibitem{delmost2}
P.~Deligne, G.D.~Mostow:
{\it Commensurabilities among lattices in $\PU (1,n)$},
Ann.\ of Math.\ Study 132, Princeton U.P., Princeton 1993.

\bibitem{demazure}
M.~Demazure:
{\it Surfaces de Del Pezzo, I-V.}
In: S\'eminaire sur les Singularit\'es des Surfaces,
Lecture Notes in Math.\ 777,
Springer-Verlag Berlin etc.\ 1980.

\bibitem{dolg}
I.~Dolgachev:
{\it Rationality of fields of invariants},
in: Algebraic Geometry Bowdoin 1985
Proc.\ Symp.\ Pure Math.\ 46 Part 2, 3--16.
AMS, Providence (RI) 1987. 

\bibitem{dolg_ort}
I.~Dolgachev, D.~Ortland:
{\it Point sets in projective spaces and theta functions}
Ast\'erisque 165 (1988),
Soc.\ Math.\ de France.

\bibitem{fadell}
E.~Fadell, J.~van Buskirk: 
{\it The braid groups of $E^2$ and $S^2$},
Duke Math.\ J.\ 29 (1962), 243--258.

\bibitem{kondo1}
S.~Kondo:
\textit{A complex hyperbolic structure on the moduli space 
of curves of genus three},
J.~reine u.\ angew.\ Math.\ 525 (2000) 219--232. 

\bibitem{kondo2}
S.~Kondo:
\textit{The moduli space of curves of genus 4 and Deligne-Mostow's
complex reflection groups}, preliminary version (2000), 10 pp.


\bibitem{konts}
M.~Kontsevich:
{\it Enumeration of rational curves via torus actions},
in: The Moduli of Curves, 335--368, 
R.~Dijkgraaf, C.~Faber, G.~van der Geer 
eds., Progress in Math.\ 129, Birkha\"user (1995).
 
\bibitem{manin}
Yu.I.~Manin: 
{\it Cubic Forms}, 2nd ed.,
North-Holland Math.\ Library,
Elsevier Sc.\ Pub.\ 1986.

\bibitem{looij:ann}
E.~Looijenga: 
\textit{Rational surfaces with an anti-canonical cycle},
Ann.\ of Math.\ 114 (1981), 267--322.

\bibitem{looijeng}
E.~Looijenga: 
\textit{New compactifications of locally symmetric varieties},
in: Proceedings of the 1984 Conference in Algebraic Geometry, 341--364,
J.~Carrell, A.V.~Geramita, P.~Russell eds. 
CMS Conference Proceedings, vol.\ 6, Amer.\ Math.\ Soc.\ Providence RI
(1984).

\bibitem{macdonald}
I.~Macdonald:
\textit{Affine root systems and Dedekind's $\eta$-function},
Invent.\ Math.\ 15 (1972), 91--143.

\bibitem{miranda}
R.~Miranda: 
{\it On the stability of pencils of cubic curves},
Amer.\ J.\ Math.\ 102 (1980), 1172--1202.

\bibitem{mirandaw}
R.~Miranda: 
{\it The Moduli of Weierstra\ss\ Fibrations over $\PP^1$},
Math.\ Ann.\ 255 (1981), 379--394.

\bibitem{miranda2}
R.~Miranda: 
{\it Persson's list of singular fibers for a rational elliptic
surface},
Math.\ Z.\ 205 (1990), 191--211.

\bibitem{morpers}
I.~Morrison, U.~Persson:
{\it Numerical sections on elliptic surfaces},
Comp.\ Math.\ 59 (1986), 323--337.

\bibitem{most}
G.D.\ Mostow:
{\it Generalized Picard lattices arising from half-integral
conditions},
Publ.\ Math.\ IHES 63 (1986), 91--106.

\bibitem{mumford}
D.~Mumford, J.~Fogarty, F.~Kirwan:
{\it Geometric Invariant Theory}, 3rd ed.,
Erg.\ der Math.\ u.\ i.\ Grenzgebiete 34,
Springer-Verlag Berlin etc.\ 1994. 

\bibitem{persson}
U. Persson: 
{\it Configurations of Kodaira fibers on rational elliptic surfaces},
Math.\ Z.\ 205 (1990), 1--47.


\bibitem{pham} F.~Pham:
{\it Formules de Picard-Lefschetz g\'en\'eralis\'ees et ramification
des int\'egrales},
Bull.\ Soc.\ Math.\ de France 93 (1965), 333--367.

\bibitem{serre} J.-P.~Serre:
{\it Cours d'Arithm\'etique},
Presses Univ.\ de France, Paris 1970.

\bibitem{sterk} 
H.~Sterk:
\textit{Compactifications of the period space of Enriques surfaces I, II},
Math.\ Z.\ 207 (1991), 1--36 and Math.\ Z.\ 220 (1995), 427--444.

\bibitem{shep}
G.C.~Shephard, J.A.~Todd:
{\it Finite unitary reflection groups},
Canadian J.\ Math.\ 6 (1954), 274--304.

\bibitem{vakil}
R.~Vakil:
{\it Twelve points on the projective line, branched covers,
and rational elliptic fibrations.} To appear.

\bibitem{vinberg}
E.B.~Vinberg:
{\it Hyperbolic reflection groups},
Russian Math.\ Surveys 40 (1985), 31--75.

\end{thebibliography}
\end{document}